\documentclass[10pt]{article}
\usepackage[utf8]{inputenc}
\usepackage{subfig}
\usepackage{amssymb}
\usepackage{bigints}
\numberwithin{equation}{section}
\usepackage{amsmath}
\usepackage{graphicx}
\usepackage[colorinlistoftodos]{todonotes}
\usepackage[colorlinks=true, allcolors=blue]{hyperref}
\usepackage{bm}
\usepackage{amsmath,amsthm,amssymb,stmaryrd}
\usepackage{enumerate}
\usepackage{enumitem}
\usepackage{relsize}
\usepackage{authblk}
\newtheorem{theorem}{Theorem}[section]

\newtheorem{corollary}{Corollary}[section]

\theoremstyle{remark}
\newtheorem*{remark}{Remark}
\newcommand{\bbint}[2]{\ensuremath{\;\backslash\!\!\!\!\backslash\!\!\!\!\!\int_{#1}^{#2}}}
\newcommand{\reglim}[2]{\ensuremath{\overset{\times}{\underset{#1\rightarrow #2}{\lim}}}}

\ExplSyntaxOn

\NewDocumentCommand{\pFq}{O{}mmmmm}
 {
  \group_begin:
  \keys_set:nn { hypergeometric } { #1 }
  \hypergeometric_print:nnnnn { #2 } { #3 } { #4 } { #5 } { #6 }
  \group_end:
 }
\NewDocumentCommand{\hypergeometricsetup}{m}
 {
  \keys_set:nn { hypergeometric } { #1 }
 }

\tl_new:N \l_hypergeometric_divider_tl
\tl_new:N \l_hypergeometric_left_tl
\tl_new:N \l_hypergeometric_right_tl

\keys_define:nn { hypergeometric }
 {
  symbol .tl_set:N = \l_hypergeometric_symbol_tl,
  symbol .initial:n = F,
  separator .tl_set:N = \l_hypergeometric_separator_tl,
  separator .initial:n = {},
  skip .tl_set:N = \l_hypergeometric_skip_tl,
  skip .initial:n = 8,
  divider .choice:,
  divider/semicolon .code:n = \tl_set:Nn \l_hypergeometric_divider_tl { \;; },
  divider/bar .code:n = \tl_set:Nn \l_hypergeometric_divider_tl { \;\middle|\; },
  divider .initial:n = semicolon,
  fences .choice:,
  fences/brack .code:n = 
   \tl_set:Nn \l_hypergeometric_left_tl {[}
   \tl_set:Nn \l_hypergeometric_right_tl {]},
  fences/parens .code:n = 
   \tl_set:Nn \l_hypergeometric_left_tl {(}
   \tl_set:Nn \l_hypergeometric_right_tl {)},
  fences .initial:n = brack,
 }

\cs_new_protected:Nn \hypergeometric_print:nnnnn
 {
  {} \sb {#1} \l_hypergeometric_symbol_tl \sb { #2 }
  \left\l_hypergeometric_left_tl
  \genfrac .. 
           {0pt} 
           {} 
           { \__hypergeometric_process:n { #3 } } 
           { \__hypergeometric_process:n { #4 } } 
  \l_hypergeometric_divider_tl
  #5
  \right\l_hypergeometric_right_tl
 }

\cs_new_protected:Nn \__hypergeometric_process:n
 {
  \clist_use:nn { #1 }
   {
    {\l_hypergeometric_separator_tl}
    \mspace { \l_hypergeometric_skip_tl mu }
   }
 }

\ExplSyntaxOff

\hypergeometricsetup{
  fences=parens,
  separator={,},
  divider=bar,
}

\title{\bf Finite-Part Integration of the Hilbert Transform}
\author[1,2]{Philip Jordan D. Blancas}
\author[1]{Eric A. Galapon}
\affil[1]{Theoretical Physics Group, National Institute of Physics, University of the Philippines, Diliman, Quezon City, 1101 Philippines}
\affil[2]{Research on Optical Science, Engineering, and Systems Laboratory, Department of Physics, Ateneo de Manila University, Loyola Heights, Quezon City, 1108 Philippines}

\date{\today}

\begin{document}

\maketitle
\begin{abstract}
The one-sided and full Hilbert transforms are evaluated exactly by means of the method of finite-part integration [E.A. Galapon, \textit{Proc. Roy. Soc. A} \textbf{473}, 20160567 (2017)]. In general, the result consists of two terms---the first is an infinite series of finite-part of divergent integrals, and the second is a contribution arising from the singularity of the kernel of transformation. The first term is precisely the result obtained when the kernel of transformation is binomially expanded in positive powers of the parameter of transformation, followed by term-by-term integration, and the resulting divergent integrals assigned values equal to their finite-parts. In all cases, the finite-part contribution is present while the presence or absence of the singular contribution depends on the interval of integration and on the parity of the function under transformation about the origin. From the exact evaluation of the Hilbert transform, the dominant asymptotic behavior for arbitrarily small parameter is obtained.
\end{abstract}


\section{Introduction}
\label{intro}

\hspace{\parindent} A Hilbert integral transform is characterized by a simple pole singularity in the interior of its contour of integration arising from its kernel of transformation, and is executed by interpreting the integral as a principal-value integral. Typically the integration ranges along the entire real line  and the transform is given by
\begin{equation}\label{ht}
    \hat{f}(\omega)=\frac{1}{\pi}\operatorname{PV}\!\!\int_{-\infty}^{\infty}  \frac{f(x)}{\omega-x}\,\mathrm{d}x ,
\end{equation}
where PV denotes the principal value. When the function $f(x)$ possesses parity symmetry, the Hilbert transform \eqref{ht} reduces into either of the forms
\begin{equation}\label{htsr}
    \frac{2}{\pi}\operatorname{PV}\!\!\int_0^{\infty} \frac{x f(x)}{\omega^2-x^2}\,\mathrm{d}x,\;\;\; \frac{2}{\pi}\operatorname{PV}\!\!\int_0^{\infty} \frac{\omega f(x)}{\omega^2-x^2}\,\mathrm{d}x  ,
\end{equation}
depending on whether $f(x)$ is odd or even, respectively.  The Hilbert transform is widely used in signal processing in diverse fields, such as in acoustics \cite{causalityapp2}, high-energy physics \cite{exptkkhighenergy}, quantum scattering theory \cite{NUSSENZVEIG1960209}, material characterization \cite{HalfHilberKramersKronig1}, optics \cite{hilbertbookking_2009}; moreover, the causality between the dispersion and attenuation functions in optical systems are related via the Hilbert transform \cite{hilbertbookking_2009b}.

The Hilbert transform \eqref{ht} and its special reductions \eqref{htsr} are routinely analyzed and evaluated by means of contour integration. In this paper we introduce the method of finite-part integration \cite{CauchyGalapon, TermByTermpaper11, doi:10.1063/1.5003479, StieltjesFinitePartpapermain13, doi:10.1063/5.0038274, FPIRegularized, StieltjesSeriesLarge} in evaluating exactly the Hilbert transform and its generalizations. Finite-part integration is a method of evaluating a well-defined integral by means of the finite-part of divergent integrals induced from the given integral itself. The method was introduced to solve the problem of missing terms that arise in evaluating the Stieltjes transform 
\begin{equation}\label{stransform}
    \int_0^a \frac{u(x)}{x^{\nu}(\omega+x)}\,\mathrm{d}x,\;\;\; 0\leq \nu<1, \;\;\; \omega>0,
\end{equation}
by term-by-term integration that leads to an infinite series of divergent integrals, with the divergent integrals assigned finite-values via analytic continuation \cite{TermByTermpaper11}. It was established in \cite{TermByTermpaper11, FPIRegularized} that the Stieltjes transform assumes the following evaluation in terms of finite-part integrals:
\begin{theorem}
Let $u(x)$ be analytic at $x=0$ and let $\rho_0$ be the radius of convergence of its Taylor expansion there. If the Stieltjes transform \eqref{stransform} exists for a given positive $a\leq\infty$ and $u(x)$ is analytic in the interval $[0,a]$, then
\begin{equation}
\label{stieltjesandhilberttransformpair2}
	\begin{split}
		\displaystyle\int_0^a
		\dfrac{u(x)}{\omega + x} {\rm d}x
		= \sum_{k=0}^\infty (-1)^k \omega^k
		\bbint{0}{a} \dfrac{u(x)}{x^{k+1}} {\rm d}x
		- u(-\omega) \ln(\omega) ,
	\end{split}
\end{equation}
\begin{equation}
	\label{stieltjesandhilberttransformpair2b}
	\begin{split}
		\displaystyle\int_0^a
		\dfrac{u(x)}{x^\nu (\omega + x)} {\rm d}x
		= \sum_{k=0}^\infty (-\omega)^k
		\bbint{0}{a} \dfrac{u(x)}{x^{k+\nu+1}} {\rm d}x
		+ \dfrac{\pi}{\sin(\pi\nu)} 
		\dfrac{u(-\omega)}{\omega^\nu} ,\;\; 0< \nu<1 ,
	\end{split}
\end{equation}
for all $\omega<\mathrm{min}(a,\rho_0)$
\end{theorem}

\noindent In both expressions, the integral $\bbint{0}{a} x^{-k-\nu-1} u(x) \; {\rm d}x$ ($0 \leq \nu < 1$) denotes the finite-part of the divergent integral $\int_0^a x^{-k-\nu-1} u(x) \; {\rm d}x$. (See Section-\ref{introtofinitepartintegrationsection}.)

It can be discerned that the summations in equations \eqref{stieltjesandhilberttransformpair2} and \eqref{stieltjesandhilberttransformpair2b} arise from the term-by-term integration of the binomial expansion of the kernel $(\omega+x)^{-1}$ about $\omega=0$ with the resulting divergent integrals assigned the values equal to their finite-parts. The second terms in \eqref{stieltjesandhilberttransformpair2} and \eqref{stieltjesandhilberttransformpair2b} are the missing terms, which are referred to as the singular contribution, when naive term-by-term integration is performed and the divergent integrals are replaced with their finite-parts. The singular contributions were recovered from the fundamental contour integral representation of the finite-part integral \cite{CauchyGalapon, TermByTermpaper11}. Here, finite-part integration of the Hilbert transform will result in expressions similar to those given by equations \eqref{stieltjesandhilberttransformpair2} and \eqref{stieltjesandhilberttransformpair2b}. But under some circumstances, the singular contribution of the Hilbert transform integral vanishes.  

Our results here do not only provide an exact evaluation of the Hilbert transform but, more importantly, lay the necessary groundwork for the application of finite-part integration in the resummation problem of divergent series appearing in perturbation theory in many areas of physics. In \cite{doi:10.1063/1.5003479}, Tica and Galapon devised a prescription to use the result \eqref{stieltjesandhilberttransformpair2} for the Stieltjes transform in the strong asymptotic regime for physical quantities assuming a Stieltjes integral representation. The prescription yields a more accurate resummation scheme than the standard scheme by Pade approximants in the non-perturbative regime. However, the resummation there can only treat alternating divergent series. But non-alternating series arise also in many contexts, such as in effective action for the vacuum polarization by a uniform electric field \cite{schwinger}, in the partition function for the self-interacting QFT \cite{mera}, and in QED effective action in time-dependent electric backgrounds \cite{dunne}, to mention a few. For this case, the Hilbert transform is expected to play the role of the Stieltjes transform in the alternating case \cite{doi:10.1063/1.5003479}. A necessary component of the resummation by finite-parts for the Stieltjes transform is knowledge of the asymptotic behavior of the Stietljes integral for arbitrarily small values of the parameter $\omega$. This information dictates the appropriate kernel of the Steiltjes transform for the resummation. It is also expected that the same asymptotic information is necessary in a resummation involving the Hilbert transform. Therefore, we do not only give an exact evaluation here but also obtain the explicit dominant behavior of the Hilbert transform for arbitrarily small values of $\omega$.

In application, we do not expect that the relevant Hilbert transform is restricted to \eqref{ht} and \eqref{htsr} with $f(x)$ analytic at the origin as is commonly assumed. Here, we extend the analysis in the presence of branch point singularity at the origin. In particular, for $0\leq\nu<1$, we evaluate the one-sided Hilbert transforms 
\begin{equation}
\label{introhilberteq3}
\operatorname{PV}\!\!\!
\displaystyle\int_{0}^a \dfrac{x^{-\nu} f(x)}{\omega-x} {\rm d}x , \;
\operatorname{PV}\!\!\!\int_0^a
\dfrac{x^{-\nu} f(x)}{\omega^2 - x^2} {\rm d}x , \;
\operatorname{PV}\!\!\!\int_0^a
\dfrac{x^{1-\nu} f(x)}{\omega^2 - x^2} {\rm d}x 
\end{equation}
for any positive $\omega<a\leq\infty$, and the full transforms
\begin{equation}
\label{introhilberteq4}
\operatorname{PV}\!\!\!
\displaystyle\int_{-a}^a \dfrac{x^{-\nu} f(x)}{\omega-x} {\rm d}x, \;
\operatorname{PV}\!\!\!
\displaystyle\int_{-a}^a \dfrac{|x|^{-\nu} f(x)}{\omega-x} {\rm d}x,\;
\operatorname{PV}\!\!\!
\displaystyle\int_{-a}^a \dfrac{|x|^{-\nu}\operatorname{sgn}(x)f(x)}{\omega-x} {\rm d}x 
\end{equation}
for any real $|\omega|<a\leq\infty$. These Hilbert transforms can serve as starting points in evaluating more Hilbert transforms. Using the results in \cite{FPIRegularized}, the results here can be extended to cover cases in the presence of logarithmic singularities at the origin. The tabulation of Hilbert transforms in terms of finite-part integrals is important as they may provide guidance in choosing the appropriate resummation scheme for non-alternating divergent series.

The rest of the paper is organized as follows. In Section-\ref{introtofinitepartintegrationsection}, we outline the method of finite-part integration. In Sections-\ref{onesidedhilbert}, \ref{fullhilbert}, and \ref{otherformhilbert}, we present various theorems involving the one-sided Hilbert transform, full Hilbert transforms, and their special reductions, respectively. In Section-\ref{examplesection}, we demonstrate some examples of using finite-part integration in evaluating Hilbert transform integrals. Finally, in Section-\ref{mellintransformmethod}, we show the various methods for evaluating finite-part integrals. In Appendix-\ref{appendix3}, we tabulate a new set of Hilbert transform integrals. And in Appendix-\ref{appendix4}, we list down the finite-part integrals used to derive the Hilbert transform integrals in Appendix-\ref{appendix3}.


\section{Finite-part Integrals}
\label{introtofinitepartintegrationsection}

\hspace{\parindent} To apply the method of finite-part integration, we will cast the Hilbert transform such that the induced divergent integrals are of the form 
\begin{equation}
\label{introfpirelevantconcept1}
\displaystyle\int_0^a 
\dfrac{f(x)}{x^{k+\nu}} {\rm d}x \;\;\;
\text{for}\;\;\; 
k = 1, 2, 3,\dots,\;\; 0\leq\nu<1,\;\; 0 < a \leq \infty ,
\end{equation}
where the divergence arises from the non-integrable singularity at the origin. The finite-part is obtained by replacing the lower limit of integration with some positive $\varepsilon<a$ and the resulting integral decomposed in the form 
\begin{equation}
\label{introfpirelevantconcept2}
\displaystyle\int_\varepsilon^a 
\dfrac{f(x)}{x^{k+\nu}} {\rm d}x 
= C_\varepsilon + D_\varepsilon  ,
\end{equation}
where $C_{\varepsilon}$ ($D_{\epsilon}$) constitutes all terms that converge (diverge) in the limit as $\varepsilon$ approaches 0. Then the finite-part is given by the limit
\begin{equation}
\label{canonical}
    \bbint{0}{a}\frac{f(x)}{x^{k+\nu}}\,{\rm d}x = \lim_{\varepsilon\rightarrow 0} C_{\varepsilon}.
\end{equation}
If the upper limit of integration happens to be infinite, the finite-part is given by
\begin{equation}
\label{introfpirelevantconcept4}
\bbint{0}{\infty}
\dfrac{f(x)}{x^{k+\nu}} {\rm d}x 
= \lim_{a \rightarrow \infty}
\displaystyle\bbint{0}{a}
\dfrac{f(x)}{x^{k+\nu}} {\rm d}x ,
\end{equation}
which we assume to exist in this paper. To uniquely define the finite part, the diverging part $D_{\varepsilon}$ must contain only diverging algebraic powers of $\varepsilon$ and $\ln\varepsilon$ (see Section-\ref{section71}). By definition, the finite-part integral always exists and is unique.

\begin{figure*}[t]
\centering\includegraphics[width=0.68\linewidth]{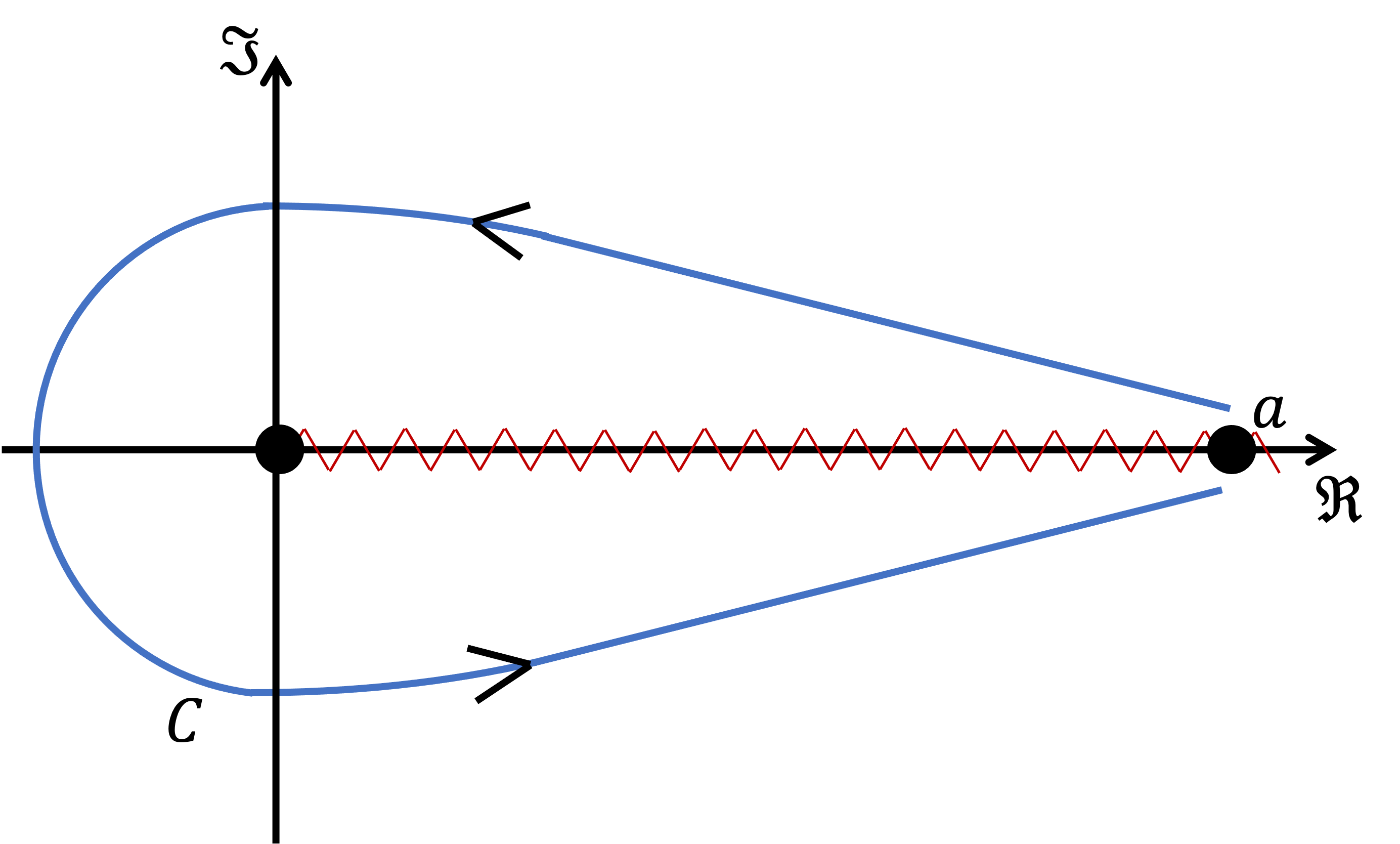}
\caption{The contour $C$ of integration in the contour integral representation of the finite-part integrals (\ref{principalvaluegeneralresult3b}) and (\ref{principalvaluegeneralresult3withbc}). The contour starts at $a$ above the branch cut of $\log z$ and $z^{-\nu}$, goes around the origin and ends at $a$ below the cut without enclosing any singularity of $f(z)$. The same contour is employed in the contour integral representation of the Hilbert transform.}
\label{contourforpositvecasefpi}
\end{figure*}

In this paper, we assume that $f(x)$ is a function of the real variable $x$ over the interval $[0,a]$ ($0<a\leq \infty$) and that it has a complex extension $f(z)$ of the complex variable $z$ that is analytic in the interval $[0,a]$. This means that $f(x)$ is the restriction of $f(z)$ in the interval $[0,a]$. By the principle of analytic continuation, the function $f(z)$ is unique and is completely determined by $f(x)$. We will refer to a function $f(x)$ over some interval $[0,a]$ having such property as complex analytic in the given interval. We denote the distance of the singularity of $f(z)$ closest to the origin by $\rho_0$. If $f(z)$ is entire, then $\rho_0=\infty$. Equivalently, $\rho_0$ is the radius of convergence of the Taylor series expansion of $f(z)$ about $z=0$.

Central to the method of finite-part integration is the fact that the finite-part possesses a contour integral representation in the complex plane \cite{TermByTermpaper11}. The representation depends on the singularity of the divergent integral at the origin: a pole singularity (corresponding to $\nu=0$) or a branch point singularity (corresponding to $\nu\neq 0$).  
\begin{theorem}\label{fpicontourrep} Let $f(x)$ be complex analytic in the interval $[0,a]$ with $f(0)\neq 0$. Then 
\begin{equation}
\label{principalvaluegeneralresult3b}
\bbint{0}{a}\dfrac{f(x)}{x^m}{\rm d}x
=\dfrac{1}{2\pi i}
\displaystyle\int_C\dfrac{f(z)}{z^m}
\left[\log(z)-i\pi\right]{\rm d}z 
\;\;\;\text{   for   }m=1, 2, 3, \cdots ,
\end{equation}
\begin{equation}
\label{principalvaluegeneralresult3withbc}
\begin{split}
\bbint{0}{a} \dfrac{f(x)}{x^{m+\nu}} {\rm d}x
= \dfrac{1}{e^{-2\pi i\nu}-1}
\displaystyle\int_C \dfrac{f(z)}{z^{m+\nu}} {\rm d}z 
\;\;\;\text{   for   } 0 < \nu < 1 ,\; m=1, 2, 3, \cdots ,
\end{split}
\end{equation}
where $f(z)$ is the complex extension of $f(x)$ in the complex plane, $\log z$ and $z^{\nu}$ take the positive real line as their branch cuts with their values coinciding with the real-valued functions $\ln x$ and $x^{\nu}$ above the branch cut, respectively, and the contour $C$ is as shown in Figure-\ref{contourforpositvecasefpi}.  
\end{theorem}


\section{One-sided Hilbert Transforms}
\label{onesidedhilbert}

\hspace{\parindent}In this section, we execute the finite-part integration of the  one-sided Hilbert transform of $x^{-\nu} f(x)$, 
\begin{equation}
\label{onesidedintro31}
\operatorname{PV}\!\!\!\displaystyle\int_0^a 
\dfrac{f(x)}{x^\nu (\omega - x)} {\rm d}x,\;\; 0\leq \nu<1, \;\; \omega>0,\;\; a\leq \infty,
\end{equation}
where $f(x)$ is complex analytic in the interval $[0,a]$. 


\begin{figure*}[ht]
\centering
\includegraphics[width=0.68\linewidth]{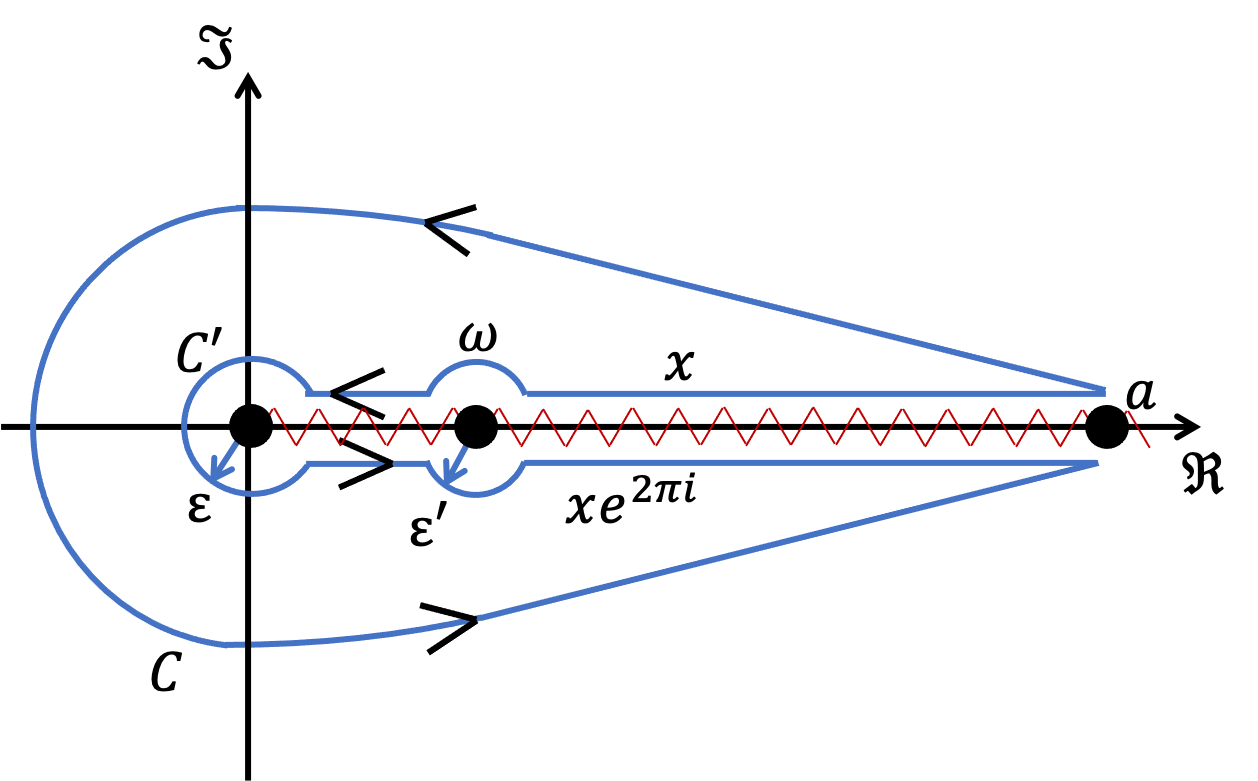}
\caption{The deformation of the contour $C$ into the contour $C'$. The deformed contour avoids the simple pole of the kernel of transformation at $x=\omega$ by semicircles above and below the cut and centered at the pole with equal radii $\varepsilon'$. }
\label{contour}
\end{figure*}

\subsection{Case $\nu=0$}

\begin{theorem}
\label{theorem1} Let $f(x)$ be complex analytic in the 
interval $[0,a]$ for some fixed positive $a<\infty$. 
If $f(0)\neq 0$, then for all $\omega\in(0,a)$ satisfying $\omega<\rho_0$, 
\begin{equation}
\label{mainequationprincipalvalueorder1onepole}
\operatorname{PV}\!\!\!\displaystyle\int_0^a \dfrac{f(x)}{\omega-x}{\rm d}x
= - \sum_{k=0}^{\infty}\omega^k
\bbint{0}{a}\dfrac{f(x)}{x^{k+1}}{\rm d}x
+ f(\omega)\ln (\omega) . 
\end{equation}
If $f(x)$ happens to have a zero at the origin of order $m$ such that $f(x)=x^m g(x)$ where $g(0)\neq 0$, then
\begin{equation} 
\label{mainequationprincipalvalueorder1onepolexmgx}
\begin{split}
\operatorname{PV}\!\!\!\displaystyle\int_0^a \dfrac{f(x)}{\omega-x}
&{\rm d}x
=-\sum_{k=0}^{m-1} \omega^k \displaystyle\int_0^a x^{m-k-1} g(x) \; {\rm d}x \\
&- \sum_{k=0}^\infty 
\omega^{m+k} \bbint{0}{a} \dfrac{g(x)}{x^{k+1}} {\rm d}x
+ \omega^m g(\omega) \ln(\omega) .
\end{split}
\end{equation}
If the principal-value integral \eqref{onesidedintro31} exists in the limit as 
$a\rightarrow\infty$, then equations \eqref{mainequationprincipalvalueorder1onepole} and \eqref{mainequationprincipalvalueorder1onepolexmgx} hold for 
$a=\infty$ as well. Furthermore, for $a=\infty$, equations \eqref{mainequationprincipalvalueorder1onepole} and \eqref{mainequationprincipalvalueorder1onepolexmgx} hold
for all positive $\omega<\rho_0$; in particular, if $f(z)$ is entire, they hold for all $\omega>0$.
\end{theorem}

\begin{proof} Let $\log(z)$ be the complex logarithmic function in Theorem-\ref{fpicontourrep} and $f(z)$ the complex extension of $f(x)$. Consider the contour integral
\begin{equation}
\label{mainequationprincipalvalueorder1onepolecomplex}
\displaystyle\int_C
\dfrac{f(z)\log(z)}{\omega-z}{\rm d}z,
\end{equation}
where the contour $C$ is as shown in Figure-\ref{contourforpositvecasefpi} and does not enclose any singularity of $f(z)$. Deforming the contour $C$ into $C'$ as indicated in Figure-\ref{contour}, with the intention to eventually take the limits $\varepsilon,\varepsilon'\to 0$, we obtain 
\begin{equation}
\label{mainequationprincipalvalueorder1onepolecontourintegration}
\begin{split}
\displaystyle\int_C
\dfrac{f(z)\log(z)}{\omega-z}
& {\rm d}z
= \displaystyle\int_a^{\omega+\varepsilon'}\dfrac{f(x)\log(x)}{\omega-x}{\rm d}x
+\displaystyle\int_{\omega-\varepsilon'}^{\varepsilon}\dfrac{f(x)\log(x)}{\omega-x}{\rm d}x \\
& +\displaystyle\int_{\varepsilon}^{\omega-\varepsilon'}\dfrac{f(x)\log(xe^{2\pi i})}{\omega-xe^{2\pi i}}{\rm d}x 
+\displaystyle\int_{\omega+\varepsilon'}^{a}\dfrac{f(x)\log(xe^{2\pi i})}{\omega-xe^{2\pi i}}{\rm d}x \\
& +\displaystyle\int_0^{\pi}\dfrac{f(\omega+\varepsilon' e^{i\phi})\log(\omega+\varepsilon' e^{i\phi})}{\omega-(\omega+\varepsilon' e^{i\phi})}i\varepsilon' e^{i\phi}{\rm d}\phi \\
&+\displaystyle\int_{\pi}^{2\pi}
\dfrac{f(\omega +\varepsilon' e^{i\phi})\log((\omega +\varepsilon' e^{i\phi}) e^{2\pi i})}{\omega-(\omega +\varepsilon' e^{i\phi})}
i\varepsilon' e^{i\phi} {\rm d}\phi \\
&+\displaystyle\int_0^{2\pi}\dfrac{f(\varepsilon e^{i\theta})\log(\varepsilon e^{i\theta})}{\omega-\varepsilon e^{i\theta}}i\varepsilon e^{i\theta}{\rm d}\theta .
\end{split}
\end{equation}
The last term of \eqref{mainequationprincipalvalueorder1onepolecontourintegration} vanishes in the limit $\varepsilon\rightarrow 0$. On the other hand, the  
fifth and sixth terms combine to take the limit
\begin{equation}
\begin{split}
&\lim_{\varepsilon'\rightarrow 0} \displaystyle\int_0^{\pi}\dfrac{f(\omega+\varepsilon' e^{i\phi})\log(\omega+\varepsilon' e^{i\phi})}{\omega-(\omega+\varepsilon' e^{i\phi})}i\varepsilon' e^{i\phi}{\rm d}\phi \\
&+ \lim_{\varepsilon'\rightarrow 0} \displaystyle\int_{\pi}^{2\pi}
\dfrac{f(\omega +\varepsilon' e^{i\phi})\log((\omega +\varepsilon' e^{i\phi}) e^{2\pi i})}{\omega-(\omega +\varepsilon' e^{i\phi})}
i\varepsilon' e^{i\phi} {\rm d}\phi
= - 2\pi i f(\omega) [ \ln(\omega) 
+ \pi i ].
\end{split}
\end{equation}
The first four terms of \eqref{mainequationprincipalvalueorder1onepolecontourintegration} constitute the desired principal-value integral in the limit $\varepsilon,\varepsilon'\rightarrow 0$. Returning all the limits back into \eqref{mainequationprincipalvalueorder1onepolecontourintegration}, the principal-value integral assumes the form
\begin{equation}\label{principalvaluegeneralresult1}
\begin{split}
\operatorname{PV}\!\!\!\displaystyle\int_0^a
\dfrac{f(x)}{\omega-x}{\rm d}x
= &\; \dfrac{1}{2\pi i}\displaystyle\int_C \dfrac{f(z)\log(z)}{\omega-z}{\rm d}z 
+ \pi i f(\omega) + f(\omega)\ln (\omega) .
\end{split}
\end{equation}

We rewrite equation \eqref{principalvaluegeneralresult1} by replacing $f(\omega)$ in the second term with its contour integral representation 
\begin{equation}
\label{cauchyprincipalvalueuniversal}
f(\omega)
=\dfrac{1}{2\pi i}\displaystyle\int_C
\dfrac{f(z)}{z-\omega}{\rm d}z 
\end{equation}
in accordance with the Cauchy integral formula, where the contour $C$ is the same contour as in \eqref{principalvaluegeneralresult1}. The equality \eqref{cauchyprincipalvalueuniversal} holds because $C$ does not enclose any singularity of $f(z)$. Then the principal-value integral becomes
\begin{equation}
\label{principalvaluegeneralresult2}
\begin{split}
\operatorname{PV}\!\!\!\displaystyle\int_0^a
\dfrac{f(x)}{\omega-x}{\rm d}x
=&\; \dfrac{1}{2\pi i}\displaystyle\int_C \dfrac{f(z)\left[\log(z)-\pi i\right]}{\omega-z}{\rm d}z 
+ f(\omega)\ln (\omega) .
\end{split}
\end{equation}
Next, we replace the kernel in the integrand with the expansion
\begin{equation}
\label{binomialexpansiontheorem1}
\dfrac{1}{\omega - z}
= -\sum_{k=0}^{n-1} \dfrac{\omega^k}{z^{k+1}}
- \frac{\omega^n}{z^n (z-\omega)} ,\;\;\;
\text{ for }\;\;n=1, 2, 3, \cdots ,
\end{equation}
to obtain 
\begin{equation}
\label{principalvaluegeneralresult2c}
\begin{split}
\operatorname{PV}\!\!\!\displaystyle\int_0^a
\dfrac{f(x)}{\omega-x}{\rm d}x
=&\; f(\omega)\ln (\omega)
- \sum_{k=0}^{n-1}\dfrac{\omega^k}{z^{k+1}}\frac{1}{2\pi i} \int_C
\frac{f(z)}{z^{k+1}}\left[\log(z)-\pi i\right] {\rm d}z 
+  R_n ,
\end{split}
\end{equation}
where the remainder term is given by
\begin{equation}\label{remainder}
R_n
= -\dfrac{\omega^n}{2\pi i}
\displaystyle\int_C
\dfrac{f(z)}{z^n (z - \omega)} 
\left[\log(z)-\pi i\right] {\rm d}z .
\end{equation}

Because $C$ does not enclose any of the singularities of $f(z)$, we recognize that the contour integrals in \eqref{principalvaluegeneralresult2c} are just the finite-parts of the divergent integrals $\int_0^a x^{-k-1} f(x) \;{\rm d}x$ in accordance with equation \eqref{principalvaluegeneralresult3b}. Now, if we choose the contour $C$ such that $\omega/|z|<1$ for all $z$ in $C$, the remainder term \eqref{remainder} vanishes as $n\rightarrow\infty$. If $f(z)$ happens to be entire, the choice is possible for all positive $\omega<a$ for we can make $C$ as large as we please. On the other hand, if $f(z)$ is not entire and the distance from the singularity closest to the origin is $\rho_0 < \infty$, the choice is possible for all positive $\omega\in(0,a)$ that satisfy $\omega<\rho_0$. Taking into account all conditions just mentioned and taking the limit as $n\rightarrow\infty$ in \eqref{principalvaluegeneralresult2c}, we arrive at the result \eqref{mainequationprincipalvalueorder1onepole} for finite $a$.

In order to obtain equation \eqref{mainequationprincipalvalueorder1onepolexmgx}, we utilize equation \eqref{binomialexpansiontheorem1} where $m=n$. In this way, we write the one-sided Hilbert transform of $f(x)=x^m g(x)$ as
\begin{equation} 
\label{proofeqnprincipalvalueorder1onepolexmgx}
\begin{split}
\operatorname{PV}\!\!\!\displaystyle\int_0^a \dfrac{f(x)}{\omega-x}
&{\rm d}x
=-\sum_{k=0}^{m-1} \omega^k \displaystyle\int_0^a x^{m-k-1} g(x) \; {\rm d}x 
+ \omega^m \operatorname{PV}\!\!\!\displaystyle\int_0^a \dfrac{g(x)}{\omega-x} {\rm d}x.
\end{split}
\end{equation}
The integrals within the sum of equation \eqref{proofeqnprincipalvalueorder1onepolexmgx} are convergent. However, the second integral on the right-hand side of the same equation is simply the one-sided Hilbert transform of $g(x)$. If we apply equation \eqref{mainequationprincipalvalueorder1onepole} to the second integral of the right-hand side, then we recover equation \eqref{mainequationprincipalvalueorder1onepolexmgx}.

It can be established that the series in equations \eqref{mainequationprincipalvalueorder1onepole} and \eqref{mainequationprincipalvalueorder1onepolexmgx} uniformly converge for all finite $a$ under the stated conditions. This allows us to perform the limit as $a\rightarrow\infty$ inside the summation so that equations \eqref{mainequationprincipalvalueorder1onepole} and \eqref{mainequationprincipalvalueorder1onepolexmgx} hold for $a=\infty$ if the Hilbert transform exists in the limit $a\rightarrow\infty$. 
\end{proof}

\begin{remark}
If we substitute $f(x) = x^m g(x)$ directly into equation \eqref{mainequationprincipalvalueorder1onepole}, we can recover equation \eqref{mainequationprincipalvalueorder1onepolexmgx} from the fact that the finite-part integrals for the first $(m-1)$ terms in the infinite series reduce to convergent integrals, which are precisely the first sum in equation \eqref{mainequationprincipalvalueorder1onepolexmgx}. Hence, in the succeeding Theorems, we will simply substitute $f(x)=x^m g(x)$ to obtain various forms of Hilbert transforms for $x^{-\nu} f(x)$ when $f(0)=0$ of order $m$.
\end{remark}

\begin{corollary}
\label{corolary31}
Under the same relevant conditions as in Theorem-\ref{theorem1}, equation \eqref{mainequationprincipalvalueorder1onepole} has the dominant behavior
\begin{equation}
\label{mainequationprincipalvalueorder1onepolexx}
\operatorname{PV}\!\!\!\displaystyle\int_0^a \dfrac{f(x)}{\omega-x}{\rm d}x
\sim f(0)\ln\omega  
\end{equation}
as $\omega \rightarrow 0$. Correspondingly, the asymptotic behavior for equation \eqref{mainequationprincipalvalueorder1onepolexmgx} is
\begin{equation} 
\label{mainequationprincipalvalueorder1onepolexx2}
\operatorname{PV}\!\!\!\displaystyle\int_0^a \dfrac{f(x)}{\omega-x}{\rm d}x 
\sim -\displaystyle\int_0^a x^{m-1} g(x) \; {\rm d}x 
\end{equation}
as $\omega\rightarrow 0$, provided that the integral does not vanish.
\end{corollary}

\begin{proof}
In equation \eqref{mainequationprincipalvalueorder1onepole}, $\ln(\omega)$ dominates over the rest of the terms, which leads to equation \eqref{mainequationprincipalvalueorder1onepolexx} provided that $f(0)\neq 0$. On the other hand, the first term within the summation in equation \eqref{mainequationprincipalvalueorder1onepolexmgx}
dominates in the same limit, which leads to \eqref{mainequationprincipalvalueorder1onepolexx2}.
\end{proof}

\begin{remark}
Similar proof applies to subsequent Corollaries to obtain the dominant terms in the asymptotic limit $\omega\rightarrow 0$. Hence, no proof be given in the subsequent Corollaries unless necessary. 
\end{remark}


\subsection{Case $\nu\neq 0$}

\begin{theorem}
\label{theorem2}
Under all relevant conditions given in Theorem-\ref{theorem1}, the principal-value integral 
\begin{equation}
\label{mainequationprincipalvalueorder1onepolewithbc}
\operatorname{PV}\!\!\!\int_{0}^{a}
\dfrac{f(x)}{x^\nu(\omega-x)}{\rm d}x
= - \sum_{k=0}^\infty \omega^k
\bbint{0}{a} \dfrac{f(x)}{x^{\nu+k+1}}{\rm d}x 
- \dfrac{\pi}{\tan(\pi \nu)} 
\dfrac{f(\omega)}{\omega^\nu} 
\end{equation}
holds for $0<\nu<1$. If $f(x)=x^m g(x)$, then
\begin{equation}
\label{mainequationprincipalvalueorder1onepolewithbcxmgx}
\begin{split}
\operatorname{PV}\!\!\!\int_{0}^{a}
\dfrac{f(x)}{x^\nu(\omega-x)}{\rm d}x
= - \sum_{k=0}^{m-1} \omega^k \displaystyle\int_0^a x^{m-k-\nu-1} g(x) \; {\rm d}x \\
- \sum_{k=0}^\infty 
\omega^{k+m} \bbint{0}{a} \dfrac{g(x)}{x^{k+\nu+1}} {\rm d}x
- \pi\cot(\pi\nu) \; \omega^{m-\nu} g(\omega) .
\end{split}
\end{equation}
\end{theorem}

\begin{proof}
The proof for this theorem is almost similar to the proof of Theorem-\ref{theorem1}. Evaluating the contour integral
\begin{equation}
\label{mainequationprincipalvalueorder1onepolecomplexwithbc}
\displaystyle\int_C
\dfrac{f(z)}{z^\nu (\omega-z)}{\rm d}z
\end{equation}
along contour $C$ as shown in Figure-\ref{contour} yields
\begin{equation}
\label{principalvaluegeneralresult1withbc}
\begin{split}
\operatorname{PV}\!\!\!\displaystyle\int_0^a
\dfrac{f(x)}{x^\nu(\omega-x)}{\rm d}x
= &\; 
\dfrac{1}{(e^{-2\pi i\nu} - 1)}\displaystyle\int_C \dfrac{f(z)}{z^\nu(\omega-z)}{\rm d}z 
- \dfrac{\pi}{\tan(\pi \nu)} 
\dfrac{f(\omega)}{\omega^\nu}.
\end{split}
\end{equation}
We then introduce the expansion \eqref{binomialexpansiontheorem1} back into the first term of \eqref{principalvaluegeneralresult1withbc}, followed by term-by-term integration over the sum. Under the same conditions as those in Theorem-\ref{theorem1}, the remainder term vanishes and the contour integrals in the sum are identified as finite-part integrals according to  
\eqref{principalvaluegeneralresult3withbc}. Then equation \eqref{mainequationprincipalvalueorder1onepolewithbc} follows.

In deriving equation \eqref{mainequationprincipalvalueorder1onepolewithbcxmgx}, we substitute $f(x) = x^m g(x)$ back into equation \eqref{mainequationprincipalvalueorder1onepolewithbc} and obtain
\begin{equation}
\label{principalvaluegeneralresult2bwithbc}
\operatorname{PV}\!\!\!\int_{0}^{a}
\dfrac{f(x)}{x^\nu(\omega-x)}{\rm d}x
= - \sum_{k=0}^\infty \omega^k
\bbint{0}{a} \dfrac{x^m g(x)}{x^{\nu+k+1}}{\rm d}x 
- \pi\cot(\pi \nu) \; \omega^{m-\nu} 
g(\omega) . 
\end{equation}
The first $(m-1)$ finite-part integrals inside the summation in equation \eqref{principalvaluegeneralresult2bwithbc} reduce to convergent integrals. As a consequence, the summation in equation \eqref{principalvaluegeneralresult2bwithbc} splits into two parts, resulting in equation \eqref{mainequationprincipalvalueorder1onepolewithbcxmgx}.
\end{proof}

\begin{corollary}
\label{corolary32}
Under the same relevant conditions as in Theorem-\ref{theorem2}, equation \eqref{mainequationprincipalvalueorder1onepolewithbc} has the dominant behavior
\begin{equation}
\label{mainequationprincipalvalueorder1onepolewithbcx}
\operatorname{PV}\!\!\!\int_{0}^{a}
\dfrac{f(x)}{x^\nu(\omega-x)}{\rm d}x
\sim 
- \dfrac{\pi}{\tan(\pi \nu)} 
\dfrac{f(0)}{\omega^\nu} 
\end{equation}
as $\omega \rightarrow 0$. And equation \eqref{mainequationprincipalvalueorder1onepolewithbcxmgx} has the dominant behavior 
\begin{equation}
\operatorname{PV}\!\!\!\int_{0}^{a}
\dfrac{f(x)}{x^\nu(\omega-x)}{\rm d}x
\sim - \displaystyle\int_0^a x^{m-\nu-1} g(x) \; {\rm d}x 
\end{equation}
in the same asymptotic limit, provided that the integral does not vanish.
\end{corollary}


\section{Full Hilbert Transform}
\label{fullhilbert}

\hspace{\parindent}In this section, we perform finite-part integration of the full Hilbert transform over a symmetric interval $[-a,a]$ and in the entire real line $\mathbb{R}$ by taking the limit $a\rightarrow\infty$. We will proceed by splitting the symmetric integral into separate integrals over the intervals $[-a,0]$ and $[0,a]$ to turn it into an integral over the interval $[0,a]$. This will lead to direct evaluation of the full Hilbert transform using the already known Stieltjes transform and the one-sided Hilbert transform just established above. Here, when we say that $f(x)$ is complex analytic in the interval $[-a,a]$, we mean that $f(x)$ has a complex extension $f(z)$, which is analytic and coincides with $f(x)$ in the interval $[-a,a]$. Again, by the uniqueness of the analytic extension, the function $f(z)$ is uniquely determined by $f(x)$.


\subsection{Case $\nu=0$}

\begin{theorem}
\label{theorem3} 
Let $f(x)$ be complex analytic in the interval $[-a,a]$ for some finite positive $a$. If $f(0)\neq 0$, then for all $\omega\in(-a,a)$ that satisfy $|\omega|<\rho_0$, 
\begin{equation}
\label{splitthewhole6maintheorem}
\operatorname{PV}\!\!\!\int_{-a}^a
\dfrac{f(x)}{\omega - x}{\rm d}x 
= \sum_{k=0}^\infty \omega^k
\bbint{0}{a}
\dfrac{{\rm d}x}{x^{k+1}}
\left[
(-1)^k f(-x)
- f(x)
\right] .
\end{equation}
When $f(x)$ is even, then the principal-value integral in \eqref{splitthewhole6maintheorem} simplifies into
\begin{equation}
\label{splitthewhole6maintheoremeven}
\operatorname{PV}\!\!\!\int_{-a}^a
\dfrac{f(x)}{\omega - x}{\rm d}x 
= -2\sum_{k=0}^\infty \omega^{2k+1}
\bbint{0}{a}
\dfrac{f(x)}{x^{2k+2}} {\rm d}x .
\end{equation}
On the other hand, if $f(x)$ has a zero at the origin of order $m$ such that $f(x)=x^m g(x)$ where $g(0)\neq 0$, then
\begin{equation}
\label{splitthewhole6maintheoremgofx}
\begin{split}
\operatorname{PV}\!\!\!\int_{-a}^a
\dfrac{f(x)}{\omega - x}{\rm d}x 
= \sum_{k=0}^{m-1} \omega^k 
\displaystyle\int_0^a x^{m-k-1} 
[(-1)^{k+m} g(-x) - g(x)] \; {\rm d}x \\
+ \sum_{k=0}^\infty 
\omega^{k+m}
\bbint{0}{a}
\dfrac{{\rm d}x}{x^{k+1}} [(-1)^k g(-x) - g(x)] .
\end{split}
\end{equation}
And when $g(x)$ is an even function, the principal-value integral in \eqref{splitthewhole6maintheoremgofx} reduces to
\begin{equation}
\label{splitthewhole6maintheoremgofxeven}
\begin{split}
\operatorname{PV}\!\!\!\int_{-a}^a
& \dfrac{f(x)}{\omega - x}{\rm d}x 
= -2 \sum_{k=0}^\infty 
\omega^{2k+m+1}
\bbint{0}{a}
\dfrac{g(x)}{x^{2k+2}} {\rm d}x \\
&- 2\sum_{k=0}^{\lfloor (m-1)/2\rfloor} \omega^{2k+1+2\lfloor m/2 \rfloor - m} 
\displaystyle\int_0^a x^{2\lfloor (m-1)/2\rfloor -2k} 
g(x) \; {\rm d}x ,
\end{split}
\end{equation}
where $\lfloor r \rfloor$ is the floor function of $r$. The results in equations \eqref{splitthewhole6maintheorem}-\eqref{splitthewhole6maintheoremgofxeven} also hold for $a=\infty$ if the principal-value integrals exist in the limit $a\rightarrow \infty$. If $a=\infty$, then the results hold for $|\omega| < \rho_0$; otherwise, if $f(x)$ is entire, then the results hold for all real $\omega$.
\end{theorem}
\begin{proof}
Let us prove equation \eqref{splitthewhole6maintheorem} for the two possible signs of $\omega$ by transforming the integration from the interval $[-a,a]$ to the interval $[0,a]$. For the case $\omega > 0$, we write the left-hand side of equation (\ref{splitthewhole6maintheorem}) as
\begin{equation}
\label{splitthewhole2}
\begin{split}
\operatorname{PV}\!\!\!\int_{-a}^a
\dfrac{f(x)}{\omega - x}{\rm d}x 
= \displaystyle\int_0^a 
\dfrac{f(-x)}{\omega + x} {\rm d}x
+ \operatorname{PV}\!\!\!\int_0^a
\dfrac{f(x)}{\omega - x}{\rm d}x .
\end{split}
\end{equation}
We recognize that the first term is the Stieltjes transform of the function $u(x)=f(-x)$, and the second term is the one-sided Hilbert transform of $f(x)$. Since $f(x)$ is complex analytic in the symmetric interval $[-a,a]$, $u(x)=f(-x)$ is necessarily complex analytic in the interval $[0,a]$ so that the result in (\ref{stieltjesandhilberttransformpair2}) holds for the first term; likewise, $f(x)$ is complex analytic in $[0,a]$ so that the result in \eqref{mainequationprincipalvalueorder1onepole} holds for the second term. Substituting the results (\ref{stieltjesandhilberttransformpair2}) and \eqref{mainequationprincipalvalueorder1onepole} into \eqref{splitthewhole2} confirms that \eqref{splitthewhole6maintheorem} holds for $\omega>0$. For the case $\omega<0$, we express the left-hand side of equation (\ref{splitthewhole6maintheorem}) as
\begin{equation}
\label{splitthewhole7}
\begin{split}
\operatorname{PV}\!\!\!\int_{-a}^a
\dfrac{f(x)}{\omega - x}{\rm d}x 
&= - \operatorname{PV}\!\!\!\int_0^a 
\dfrac{f(-x)}{(-\omega) - x} {\rm d}x
- \displaystyle\int_0^a
\dfrac{f(x)}{(-\omega) + x}{\rm d}x .
\end{split}
\end{equation}
This time the first term is the one-sided Hilbert transform of $f(-x)$ in the parameter $|\omega|$, and the second term is the Stieltjes transform of $f(x)$ in the same parameter. Substituting (\ref{stieltjesandhilberttransformpair2}) and \eqref{mainequationprincipalvalueorder1onepole} into \eqref{splitthewhole7} likewise confirms that \eqref{splitthewhole6maintheorem} holds for $\omega<0$. 

Let us move on to the case when $f(x)$ is an even function. Here, we must impose $f(x)=f(-x)$ on equation \eqref{splitthewhole6maintheorem}. In doing so, all even orders of $k$ vanish. Hence, we shift $k \rightarrow 2k+1$ to recover equation \eqref{splitthewhole6maintheoremeven}. 

Next, we prove equation \eqref{splitthewhole6maintheoremgofx} by substituting $f(x) = x^m g(x)$ into equation \eqref{splitthewhole6maintheorem}. The substitution yields
\begin{equation}
\label{splitthewhole7b}
\operatorname{PV}\!\!\!\int_{-a}^a
\dfrac{f(x)}{\omega - x}{\rm d}x 
= \sum_{k=0}^\infty \omega^k
\bbint{0}{a}
\dfrac{{\rm d}x}{x^{k+1}}
\left[
(-1)^k (-x)^m g(-x)
- x^m g(x)
\right] .
\end{equation}
Again the finite-part integrals up to the $k=m-1$ term in the infinite series reduce to convergent integrals, leading to the splitting of the summation into sums of convergent integrals and sums of finite-part integrals. The result is equation \eqref{splitthewhole6maintheoremgofx}.

Finally, we derive equation \eqref{splitthewhole6maintheoremgofxeven} by imposing $g(x)=g(-x)$ on equation \eqref{splitthewhole6maintheoremgofx}. The result is
\begin{equation} 
\label{splitthewhole7c}
\begin{split}
\operatorname{PV}\!\!\!\int_{-a}^a
\dfrac{f(x)}{\omega - x}{\rm d}x 
=& \sum_{k=0}^\infty 
\omega^{k+m}
[(-1)^k - 1]
\bbint{0}{a}
\dfrac{g(x)}{x^{k+1}} {\rm d}x \\
&+ \sum_{k=0}^{m-1} \omega^k 
[(-1)^{k+m} - 1]
\displaystyle\int_0^a x^{m-k-1} 
g(x) \; {\rm d}x  .
\end{split}
\end{equation}
The terms corresponding to even $k$ in the first sum vanish. Shifting $k\rightarrow 2k+1$ verifies the first summation in equation \eqref{splitthewhole6maintheoremgofxeven}. On the other hand, some terms in the second summation of \eqref{splitthewhole7c} vanish when $k+m$ is an even number. Hence, when $m$ is even, the non-vanishing terms appear if $k$ is odd, leading to
\begin{equation} 
\label{splitthewhole8a}
-2 \sum_{k=1, 3, 5, \cdots}^{m-1} \omega^k 
\displaystyle\int_0^a x^{m-k-1} 
g(x) \; {\rm d}x
= -2 \sum_{k'=0}^{(m-2)/2} \omega^{2k'+1} \displaystyle\int_0^a x^{m-2k'-2} g(x) \; {\rm d}x ,
\end{equation}
and when $m$ is odd, the non-vanishing terms appear if $k$ is even such as
\begin{equation} 
\label{splitthewhole8b}
-2 \sum_{k=0, 2, 4, \cdots}^{m-1} \omega^k 
\displaystyle\int_0^a x^{m-k-1} 
g(x) \; {\rm d}x
= -2 \sum_{k'=0}^{(m-1)/2} \omega^{2k'} \displaystyle\int_0^a x^{m-2k'-1} g(x) \; {\rm d}x .
\end{equation} 
In order to express equations \eqref{splitthewhole8a} and \eqref{splitthewhole8b} as a single equation, we use the floor function $\lfloor r \rfloor$ wherein we consider the largest integer that is less than $r$. Through floor function, we define $\lfloor m/2 \rfloor$ as
\begin{equation} 
\label{splitthewhole8c}
\lfloor m/2 \rfloor
=\begin{cases}
m/2 &\text{for} \;\;\;\;\; m=\text{even} , \\
\\
(m-1)/2 &\text{for} \;\;\;\;\; m=\text{odd} . 
\end{cases}
\end{equation}
Hence, equation \eqref{splitthewhole8c} recovers the second summation in equation \eqref{splitthewhole6maintheoremgofxeven}.
\end{proof}

\begin{corollary}
Under the same relevant conditions as in Theorem-\ref{splitthewhole6maintheorem}, equation \eqref{splitthewhole6maintheorem} has the dominant behavior
\begin{equation}
\label{splitthewhole6maintheoremapprox}
\operatorname{PV}\!\!\!\int_{-a}^a
\dfrac{f(x)}{\omega - x}{\rm d}x 
\sim \bbint{0}{a}
\dfrac{{\rm d}x}{x}
\left[ f(-x) - f(x) \right] \end{equation}
as $\omega\rightarrow 0$, provided that the finite-part integral does not vanish, a condition implied in the subsequent results where the leading contribution involves finite-part or convergent integral. If $f(x)$ happens to be even with $f(0)\neq 0$, then 
\begin{equation}
\label{splitthewhole6maintheoremevenapprox}
\operatorname{PV}\!\!\!\int_{-a}^a
\dfrac{f(x)}{\omega - x}{\rm d}x 
\sim -2 \omega
\bbint{0}{a}
\dfrac{f(x)}{x^{2}} {\rm d}x, \;\;\;\omega \rightarrow 0 .
\end{equation}
When $f(x)=x^m g(x)$ with $g(0)\neq 0$, equation \eqref{splitthewhole6maintheoremgofx} has the dominant behavior
\begin{equation}
\label{splitthewhole6maintheoremgofxapprox}
\begin{split}
\operatorname{PV}\!\!\!\int_{-a}^a
\dfrac{f(x)}{\omega - x}{\rm d}x 
\sim \displaystyle\int_0^a x^{m-1} 
[(-1)^{m} g(-x) - g(x)] \; {\rm d}x, \;\;\;\omega \rightarrow 0 ,
\end{split}
\end{equation}
and when $g(x)$ is even, 
\begin{equation}
\begin{split}
\operatorname{PV}\!\!\!\int_{-a}^a
\dfrac{f(x)}{\omega - x}{\rm d}x 
\sim -2
\omega^{2\lfloor m/2 \rfloor - m + 1} 
\displaystyle\int_0^a x^{2\lfloor (m-1)/2\rfloor} 
g(x) \; {\rm d}x
, \;\;\;\omega \rightarrow 0 .
\end{split}
\end{equation}
\end{corollary}


The inclusion of the full Hilbert transform for $|x|^{-\nu} \textup{sgn}(x) f(x)$ is necessary because the Hilbert transform integral 
\begin{equation}
\label{HTnotKKrelations}
\operatorname{PV}\!\!\!\int_{-\infty}^\infty \dfrac{\operatorname{sgn}(x)  f(x)}{ |x|^\nu (\omega-x)} {\rm d}x
\end{equation}
reduces into
\begin{equation}
\label{HTKKrelations}
\operatorname{PV}\!\!\!\int_{0}^\infty \dfrac{ f(x)}{x^\nu (\omega^2-x^2)} {\rm d}x \;\;\;
\text{or}
\;\;\;
\operatorname{PV}\!\!\!\int_{0}^\infty \dfrac{ x f(x)}{x^\nu (\omega^2-x^2)} {\rm d}x
\end{equation}
for $0\leq \nu<1$ if $f(x)$ is a symmetric function. In fact, some Hilbert transforms of symmetric $f(x)$ in Appendix-\ref{appendix3} possess the congruence between equations \eqref{HTnotKKrelations} and \eqref{HTKKrelations}. 

\begin{theorem}
\label{theorem3b} 
Under all relevant conditions given in Theorem-\ref{theorem3}, the principal-value integral
\begin{equation}
\label{splitthewhole6maintheoremwsignum}
\operatorname{PV}\!\!\!\int_{-a}^a
\dfrac{\textup{sgn}(x) f(x)}{\omega - x}{\rm d}x 
= -\sum_{k=0}^\infty \omega^k
\bbint{0}{a}
\dfrac{{\rm d}x}{x^{k+1}}
\left[
(-1)^k f(-x)
+ f(x)
\right] + 2f(\omega) \ln|\omega| 
\end{equation}
holds for $f(0)\neq 0$. If $f(x)$ is even, then equation \eqref{splitthewhole6maintheoremwsignum} simplifies into
\begin{equation}
\label{splitthewhole6maintheoremwsignumeven}
\operatorname{PV}\!\!\!\int_{-a}^a
\dfrac{\textup{sgn}(x) f(x)}{\omega - x}{\rm d}x 
= -2\sum_{k=0}^\infty \omega^{2k}
\bbint{0}{a}
\dfrac{f(x)}{x^{2k+1}} {\rm d}x
+ 2f(\omega) \ln|\omega| .
\end{equation}
When it happens that $f(x) = x^m g(x)$, then
\begin{equation}
\label{splitthewhole6maintheoremwsignumxmgx}
\begin{split}
\operatorname{PV}\!\!\!\int_{-a}^a
\dfrac{\textup{sgn}(x) f(x)}{\omega - x}{\rm d}x 
= -\sum_{k=0}^{m-1} 
\omega^k \displaystyle\int_0^a
x^{m-k-1} 
[(-1)^{k+m} g(-x) + g(x)] \; {\rm d}x \\
- \sum_{k=0}^\infty
\omega^{k + m} \bbint{0}{a}
\dfrac{{\rm d}x}{x^{k+1}} 
[(-1)^k g(-x) + g(x)]
+ 2 \omega^m g(\omega) \ln|\omega| ,
\end{split}
\end{equation}
and when $g(x)$ is an even function, the principal-value integral in \eqref{splitthewhole6maintheoremwsignumxmgx} reduces to
\begin{equation}
\label{splitthewhole6maintheoremwsignumxmgxeven}
\begin{split}
\operatorname{PV}\!\!\!\int_{-a}^a
\dfrac{\textup{sgn}(x) f(x)}{\omega - x}{\rm d}x 
= - 2\sum_{k=0}^\infty
\omega^{2k + m} \bbint{0}{a}
\dfrac{g(x)}{x^{2k+1}} {\rm d}x 
+ 2 \omega^m g(\omega) \ln|\omega| \\
-2\sum_{k=0}^{\lfloor (m-2)/2 \rfloor } 
\omega^{2k + 2\lfloor (m+1)/2 \rfloor - m} \displaystyle\int_0^a
x^{2\lfloor (m-2)/2 \rfloor -2k+1} 
g(x) \; {\rm d}x  .
\end{split}
\end{equation}
\end{theorem}

\begin{proof}
The proof for this theorem is almost the same as illustrated in Theorem-\ref{theorem3}. However, we utilize first the definition of $\text{sgn}(x)$, which is given by
\begin{equation}
\label{signum}
\text{sgn}(x)
= \begin{cases}
+ 1 &\text{for} \;\;\;\; x > 0 , \\
\\
- 1 &\text{for} \;\;\;\; x < 0 ,
\end{cases}
\end{equation}
before we split the principal-value integral in equation \eqref{splitthewhole6maintheoremwsignum} into two integrals involving a Stieltjes transform and a one-sided Hilbert transform.
\end{proof}

\begin{corollary}
Under the same relevant conditions as in Theorem-\ref{theorem3b}, equations \eqref{splitthewhole6maintheoremwsignum} and 
\eqref{splitthewhole6maintheoremwsignumeven} have the dominant behavior
\begin{equation}
\operatorname{PV}\!\!\!\int_{-a}^a
\dfrac{\textup{sgn}(x) f(x)}{\omega - x}{\rm d}x 
\sim 2f(0) \ln|\omega| 
\end{equation}
as $\omega\rightarrow 0$. When $f(x)=x^m g(x)$, equation \eqref{splitthewhole6maintheoremwsignumxmgx} has the dominant behavior
\begin{equation}
\begin{split}
\operatorname{PV}\!\!\!\int_{-a}^a
\dfrac{\textup{sgn}(x) f(x)}{\omega - x}{\rm d}x 
\sim - \displaystyle\int_0^a
x^{m-1} 
[(-1)^{m} g(-x) + g(x)] \; {\rm d}x 
\end{split}
\end{equation}
as $\omega \rightarrow 0$, provided that the integral does not vanish. And when $g(x)$ is even symmetric, then equation \eqref{splitthewhole6maintheoremwsignumxmgxeven} has the dominant behavior
\begin{equation}
\operatorname{PV}\!\!\!\int_{-a}^a
\dfrac{\textup{sgn}(x) f(x)}{\omega - x}{\rm d}x 
\sim
\begin{cases} 
-2\omega\displaystyle\bbint{0}{a} x^{-1}  g(x) \; {\rm d}x, & m=1 , \\
\\
-2 \dfrac{\omega^{ 2\lfloor (m+1)/2 \rfloor} }{\omega^m} \displaystyle\int_0^a
x^{2\lfloor (m-2)/2 \rfloor +1} g(x) \; {\rm d}x,  
& m \geq 2 
\end{cases}
\end{equation}
in the same asymptotic limit, as long as the finite-part or the convergent integral does not vanish.
\end{corollary}


\subsection{Case $\nu\neq 0$: Not in Absolute Value}

\hspace{\parindent} In this Section we evaluate the full Hilbert transform of $x^{-\nu}f(x)$ for $0<\nu<1$. The integration is to be performed above the negative real axis where $x^{-\nu}=e^{-\pi\nu i} |x|^{-\nu}$.

\begin{theorem} \label{theorem43}
Under all relevant conditions given in Theorem-\ref{theorem3}, the principal-value integral
\begin{equation}
\label{splitthewhole1withbcmain}
\begin{split}
\operatorname{PV}\!\!\!\int_{-a}^a
& \dfrac{f(x)}{x^\nu(\omega - x)} {\rm d}x 
= \sum_{k=0}^\infty
\omega^k \bbint{0}{a}
\dfrac{{\rm d}x}{x^{k+\nu+1}}
\left[
e^{- i\pi\nu} 
(-1)^{k}  f(-x) 
- f(x)
\right] 
- i\pi \dfrac{f(\omega)}{\omega^\nu} 
\end{split}
\end{equation}
holds for $0 < \nu < 1$. If $f(x)$ is even, then equation \eqref{splitthewhole1withbcmain} reduces to 
\begin{equation}
\label{splitthewhole1withbcmainevenmain}
\begin{split}
\operatorname{PV}\!\!\!\int_{-a}^a
& \dfrac{f(x)}{x^\nu(\omega - x)} {\rm d}x 
= - 2i\sin(\pi\nu/2) e^{- i\pi\nu/2}
\sum_{k=0}^\infty \omega^{2k}
\bbint{0}{a} \dfrac{f(x)}{x^{2k+\nu+1}} {\rm d}x \\
&- 2\cos(\pi\nu/2)  e^{- i\pi\nu/2}
\sum_{k=0}^\infty \omega^{2k+1}
\bbint{0}{a} \dfrac{f(x)}{x^{2k+\nu+2}} {\rm d}x
- i\pi \dfrac{f(\omega)}{\omega^\nu} .
\end{split}
\end{equation}
If $f(x) = x^m g(x)$, then
\begin{equation}
\label{splitthewhole1withbcmainxgxmain}
\begin{split}
\operatorname{PV}\!\!\!\int_{-a}^a
\dfrac{f(x)}{x^\nu(\omega - x)} {\rm d}x
= \sum_{k=0}^{m-1}
\omega^k \displaystyle\int_0^a
\dfrac{\left[(-1)^{k+m}    e^{- i\pi\nu} g(-x) - g(x) \right] {\rm d}x }{x^{k-m+\nu+1}} 
 \\
+\sum_{k=0}^\infty \omega^{k+m} \bbint{0}{a} \dfrac{{\rm d}x}{x^{k+\nu+1}} [(-1)^k e^{-i\pi\nu} g(-x) - g(x)]
- i\pi \omega^{m-\nu} g(\omega) ,
\end{split}
\end{equation}
and if $g(x)$ is an even function, equation \eqref{splitthewhole1withbcmainxgxmain} simplifies to
\begin{equation}
\label{splitthewhole1withbcmainxgxmaineven}
\begin{split}
\operatorname{PV}
&\!\!\!\int_{-a}^a
\dfrac{f(x)}{x^\nu(\omega - x)} {\rm d}x
= - i\pi \omega^{m-\nu} g(\omega) 
-2i e^{-i\pi\nu/2} \sin(\pi\nu/2) \sum_{k=0}^\infty \omega^{2k+m} \bbint{0}{a} \dfrac{g(x) \; {\rm d}x}{x^{2k+\nu+1}} \\
&\hspace{0.6cm}
-2 e^{-i\pi\nu/2} \cos(\pi\nu/2) \sum_{k=0}^\infty \omega^{2k+m+1} \bbint{0}{a} \dfrac{g(x) \; {\rm d}x}{x^{2k+\nu+2}}\\
&-2i e^{-i\pi\nu/2}\sin(\pi\nu/2) \sum_{k=0}^{\lfloor (m-2)/2 \rfloor}
\omega^{2k+2\lfloor (m+1)/2 \rfloor - m} \displaystyle\int_0^a
\dfrac{g(x) \; {\rm d}x }{x^{2k+\nu-1-2\lfloor (m-2)/2 \rfloor}}
 \\
&\hspace{0.6cm}
-2 e^{-i\pi\nu/2}\cos(\pi\nu/2) \sum_{k=0}^{\lfloor (m-1)/2 \rfloor}
\omega^{2k+1+2\lfloor m/2 \rfloor - m} \displaystyle\int_0^a
\dfrac{g(x) \; {\rm d}x }{x^{2k+\nu-2\lfloor (m-1)/2 \rfloor}}  .
\end{split}
\end{equation}
\end{theorem}

\begin{proof}
The key step in proving this theorem is by utilizing $x^{-\nu} = e^{-i\pi\nu} |x|^{-\nu}$ along the negative real axis. Hence, if we split the left-hand side of equation \eqref{splitthewhole1withbcmain} as two integrals, then
\[
\begin{split}
\operatorname{PV}\!\!\!\int_{-a}^a
\dfrac{f(x) \; {\rm d}x}{x^\nu(\omega - x)}
= 
 e^{-i\pi\nu} \displaystyle\int_0^{a}
\dfrac{f(-x) \; {\rm d}x}{x^\nu(\omega + x)}
+ \operatorname{PV}\!\!\!\int_0^a
\dfrac{f(x) \; {\rm d}x}{x^\nu(\omega - x)} 
\end{split}
\]
for $\omega > 0$ and 
\[
\begin{split}
\operatorname{PV}\!\!\!\int_{-a}^a
\dfrac{f(x) \;  {\rm d}x }{x^\nu(\omega - x)} 
= - e^{-i\pi\nu} \operatorname{PV}\!\!\!\int_0^{a}
 \dfrac{f(-x) \; {\rm d}x}
{x^\nu((-\omega) - x)} 
- \displaystyle\int_0^a
\dfrac{f(x) \; {\rm d}x}{x^\nu((-\omega) + x)}
\end{split}
\]
for $\omega < 0$. 
The subsequent steps require the application of equations \eqref{stieltjesandhilberttransformpair2b} and \eqref{mainequationprincipalvalueorder1onepolewithbc}, which leads to the results given in equation \eqref{splitthewhole1withbcmain}. We also perform the same procedures as executed in the previous theorems to derive equations \eqref{splitthewhole1withbcmainevenmain}-\eqref{splitthewhole1withbcmainxgxmaineven}.
\end{proof}

\begin{corollary}
Under the same relevant conditions as in Theorem-\ref{theorem43}, equations \eqref{splitthewhole1withbcmain} and \eqref{splitthewhole1withbcmainevenmain} have the dominant behavior 
\begin{equation}
\label{splitthewhole6maintheoremwsignumapprox}
\operatorname{PV}\!\!\!\int_{-a}^a
\dfrac{f(x)}{x^\nu(\omega - x)}{\rm d}x 
\sim -i\pi 
\dfrac{f(0)}{\omega^\nu} .
\end{equation}
as $\omega\rightarrow 0$. When $f(x)=x^m g(x)$, the principal-value integral in \eqref{splitthewhole1withbcmainxgxmain} has the dominant behavior
\begin{equation}
\label{splitthewhole6maintheoremwsignumapproxb}
\operatorname{PV}\!\!\!\int_{-a}^a
\dfrac{f(x)}{x^\nu(\omega - x)}{\rm d}x \sim 
\displaystyle\int_0^a 
\dfrac{(-1)^m e^{-i\pi\nu} g(-x) - g(x)}{x^{\nu-m+1}} {\rm d}x 
\end{equation}
in the same asymptotic limit, provided that the integral does not vanish. And if $g(x)$ is even symmetric, the dominant term of equation \eqref{splitthewhole1withbcmainxgxmaineven} is
\begin{equation}
\label{splitthewhole6maintheoremwsignumapproxc}
\operatorname{PV}\!\!\!\int_{-a}^a
\dfrac{f(x) \; {\rm d}x}{x^\nu(\omega - x)} \sim 
\begin{cases}
-2 e^{-i\pi\nu/2} \cos(\pi\nu/2) \displaystyle\int_0^a \dfrac{g(x) \; {\rm d}x}{x^{\nu-m+1}}  
&\text{for} \;\; 
m = \text{odd} , \\
\\
-2i e^{-i\pi\nu/2} \sin(\pi\nu/2)
\displaystyle\int_0^a \dfrac{g(x) \; {\rm d}x}{x^{\nu-m+1}}  
&\text{for} \;\;
m = \text{even} 
\end{cases}
\end{equation}
in the same asymptotic limit, provided that the convergent integrals do not vanish as well.
\end{corollary}


\subsection{Case $\nu\neq 0$: In Absolute Value}

\begin{theorem}
\label{theorem44}
Under all relevant conditions given in Theorem-\ref{theorem3}, the principal-value integral 
\begin{equation}
\label{splitthewhole1withbcmainabspositive}
\begin{split}
\operatorname{PV}\!\!\!\int_{-a}^a
\dfrac{f(x)}{|x|^\nu(\omega - x)} {\rm d}x 
= &\; \sum_{k=0}^\infty 
\omega^k \bbint{0}{a}
\dfrac{{\rm d}x}{x^{k+\nu+1}}
\left[
(-1)^{k}  f(-x) 
- f(x)
\right] \\
&\; + \pi \tan(\pi\nu/2) \; \textup{sgn}(\omega) \dfrac{f(\omega)}{|\omega|^\nu} 
\end{split}
\end{equation}
holds for $0 < \nu < 1$. If $f(x)$ is even, then
\begin{equation}
\label{splitthewhole1withbcmainabspositiveeven}
\begin{split}
\operatorname{PV}\!\!\!\int_{-a}^a
\dfrac{f(x)}{|x|^\nu(\omega - x)} {\rm d}x 
= &\; -2 \sum_{k=0}^\infty 
\omega^{2k+1} \bbint{0}{a}
\dfrac{f(x)}{x^{2k+\nu+2}} {\rm d}x 
+ \pi \tan(\pi\nu/2) \; \textup{sgn}(\omega) \dfrac{f(\omega)}{|\omega|^\nu} .
\end{split}
\end{equation}
If $f(x) = x^m g(x)$, then equation \eqref{splitthewhole1withbcmainabspositive} becomes
\begin{equation}
\label{splitthewhole1withbcmainabspositivexgx}
\begin{split}
\operatorname{PV}\!\!\!\int_{-a}^a
& \dfrac{f(x)}{|x|^\nu(\omega - x)} {\rm d}x 
= \sum_{k=0}^{m-1} 
\omega^k \displaystyle\int_0^a
\dfrac{\left[ 
(-1)^{k+m} g(-x) - g(x)
\right] {\rm d}x }{x^{k -m+\nu+1} }
\\
&+ \sum_{k=0}^\infty \omega^{k+m} \bbint{0}{a} \dfrac{\left[(-1)^k g(-x) - g(x)\right] {\rm d}x}{x^{k+\nu+1}} 
+ \pi  \dfrac{\textup{sgn}(\omega)}{\cot(\pi\nu/2)} \dfrac{\omega^m g(\omega)}{|\omega|^\nu} .
\end{split}
\end{equation}
And if $g(x)$ is even, then equation \eqref{splitthewhole1withbcmainabspositivexgx} reduces to
\begin{equation}
\label{splitthewhole1withbcmainabspositivexgxeven}
\begin{split}
\operatorname{PV}\!\!\!\int_{-a}^a
& \dfrac{f(x)}{|x|^\nu(\omega - x)} {\rm d}x 
= -2 \sum_{k=0}^{\lfloor (m-1)/2 \rfloor} 
\omega^{2k+1+2\lfloor m/2 \rfloor - m} \displaystyle\int_0^a
\dfrac{ g(x) \; {\rm d}x }{x^{2k -2 \lfloor (m-1)/2 \rfloor +\nu} }
\\
&- 2 \sum_{k=0}^\infty \omega^{2k+m+1} \bbint{0}{a} \dfrac{g(x) }{x^{2k+\nu+2}} {\rm d}x
+ \pi  \dfrac{\textup{sgn}(\omega)}{\cot(\pi\nu/2)} \dfrac{\omega^m g(\omega)}{|\omega|^\nu} .
\end{split}
\end{equation}
\end{theorem}
\begin{proof} We proceed in the same manner as in the Theorem-\ref{theorem43} to establish the results. But this time, we need to utilize the definition for $|x|$, which is given by
\begin{equation}
\label{absolutevalue}
|x|
=\begin{cases}
x &\text{if}\;\;\;\;\; x>0, \\
\\
-x &\text{if}\;\;\;\;\; x<0 .
\end{cases}
\end{equation}
\end{proof}

\begin{corollary}
Under the same relevant conditions as in Theorem-\ref{theorem44}, equations \eqref{splitthewhole1withbcmainabspositive} and \eqref{splitthewhole1withbcmainabspositiveeven} have the dominant behavior 
\begin{equation}
\label{splitthewhole1withbcmainabspositivecorollary1}
\operatorname{PV}\!\!\!\int_{-a}^a
\dfrac{f(x)}{|x|^\nu(\omega - x)}{\rm d}x 
\sim \pi \tan(\pi\nu/2)  \; \textup{sgn}(\omega) \dfrac{f(0)} {|\omega|^\nu} 
\end{equation}
as $\omega \rightarrow 0$. When $f(x)=x^m g(x)$, equation \eqref{splitthewhole1withbcmainabspositivexgx} has the dominant behavior
\begin{equation}
\begin{split}
\operatorname{PV}\!\!\!\int_{-a}^a
& \dfrac{f(x)}{|x|^\nu(\omega - x)} {\rm d}x \sim \displaystyle\int_0^a
x^{m - \nu - 1}
\left[ (-1)^{m} g(-x) - g(x) \right]
{\rm d}x 
\end{split}
\end{equation}
in the same asymptotic limit provided that the integral does not vanish. And if $g(x)$ is even symmetric, equation \eqref{splitthewhole1withbcmainabspositivexgxeven} has the dominant behavior
\begin{equation}
\begin{split}
\operatorname{PV}\!\!\!\int_{-a}^a
\dfrac{f(x)}{|x|^\nu(\omega - x)}{\rm d}x 
\sim -2 
\omega^{2\lfloor m/2 \rfloor - m + 1} \displaystyle\int_0^a
x^{2 \lfloor (m-1)/2 \rfloor - \nu} g(x) \; {\rm d}x 
\end{split}
\end{equation}
in the same asymptotic limit, as long as the integral does not vanish as well.
\end{corollary}


\begin{theorem}
\label{theorem45}
Under all relevant conditions given in Theorem-\ref{theorem3}, the principal-value integral
\begin{equation}
\label{splitthewhole1withbcmainabspositivesgn}
\begin{split}
& \operatorname{PV}\!\!\!\int_{-a}^a
\dfrac{\textup{sgn}(x) f(x)}{|x|^\nu(\omega - x)} {\rm d}x \\
&\hspace{1cm}
= - \sum_{k=0}^\infty 
\omega^k \bbint{0}{a}
\dfrac{{\rm d}x}{x^{k+\nu+1}}
\left[
(-1)^{k}  f(-x) 
+ f(x)
\right] 
- \pi \cot(\pi\nu/2)
\dfrac{f(\omega)}{|\omega|^\nu}  
\end{split}
\end{equation}
holds for $0 < \nu < 1$. If $f(x)$ is even, then
\begin{equation}
\label{splitthewhole1withbcmainabspositivesgneven}
\begin{split}
& \operatorname{PV}\!\!\!\int_{-a}^a
\dfrac{\textup{sgn}(x) f(x)}{|x|^\nu(\omega - x)} {\rm d}x 
= - 2 \sum_{k=0}^\infty 
\omega^{2k} \bbint{0}{a}
\dfrac{f(x)}{x^{2k+\nu+1}} {\rm d}x
- \pi \cot(\pi\nu/2)
\dfrac{f(\omega)}{|\omega|^\nu} .
\end{split}
\end{equation}
If $f(x) = x^m g(x)$, then equation \eqref{splitthewhole1withbcmainabspositivesgn} becomes
\begin{equation}
\label{splitthewhole1withbcmainabspositivesgnxgx}
\begin{split}
\operatorname{PV}\!\!\! \int_{-a}^a
& \dfrac{\textup{sgn}(x) f(x)}{|x|^\nu(\omega - x)} {\rm d}x 
= -\sum_{k=0}^{m-1} 
\omega^k \int_{0}^{a} \dfrac{\left[ (-1)^{k+m} g(-x)+ g(x) \right] {\rm d}x}{x^{k-m+\nu+1}} \\
&- \sum_{k=0}^\infty \omega^{k+m} \bbint{0}{a} \dfrac{\left[(-1)^k g(-x) + g(x) \right] {\rm d}x}{x^{k+\nu+1}} 
- \pi \cot(\pi\nu/2) \dfrac{\omega^m g(\omega)}{|\omega|^\nu} ,
\end{split}
\end{equation}
and if $g(x)$ is even, then equation \eqref{splitthewhole1withbcmainabspositivesgnxgx} simplifies into
\begin{equation}
\label{splitthewhole1withbcmainabspositivesgnxgxeven}
\begin{split}
\operatorname{PV}\!\!\! \int_{-a}^a
& \dfrac{\textup{sgn}(x) f(x)}{|x|^\nu(\omega - x)} {\rm d}x 
= - 2 \sum_{k=0}^\infty \omega^{2k+m} \bbint{0}{a} \dfrac{g(x)}{x^{2k+\nu+1}} {\rm d}x
- \pi \cot(\pi\nu/2) \dfrac{\omega^m g(\omega)}{|\omega|^\nu}
\\
&\hspace{0.5cm}
- 2 \sum_{k=0}^{\lfloor (m-2)/2 \rfloor} 
\omega^{2k + 2\lfloor (m+1)/2 \rfloor - m} \int_{0}^{a} \dfrac{ g(x) \; {\rm d}x}{x^{2k-2\lfloor (m-2)/2 \rfloor +\nu- 1}}  .
\end{split}
\end{equation}
\end{theorem}

\begin{proof} We use the definitions for $\operatorname{sgn}(x)$ and  $|x|$ given in equations \eqref{signum} and \eqref{absolutevalue}, respectively, and proceed in the same manner as in Theorem-\ref{theorem43} to establish the results.
\end{proof}

\begin{corollary}
Under the same relevant conditions as in Theorem-\ref{theorem45}, equations \eqref{splitthewhole1withbcmainabspositivesgn} and \eqref{splitthewhole1withbcmainabspositivesgneven} have the dominant behavior
\begin{equation}
\label{splitthewhole1withbcmainabspositivesgncorollary}
\operatorname{PV}\!\!\!\int_{-a}^a
\dfrac{\textup{sgn}(x)f(x)}{|x|^\nu(\omega - x)}{\rm d}x 
\sim -\pi \cot(\pi\nu/2)  \dfrac{f(0)} {|\omega|^\nu} 
\end{equation}
as $\omega \rightarrow 0$.
When $f(x) = x^m g(x)$, equation \eqref{splitthewhole1withbcmainabspositivesgnxgx}
has the dominant behavior
\begin{equation}
\label{splitthewhole1withbcmainabspositivesgncorollaryb}
\begin{split}
\operatorname{PV}\!\!\! \int_{-a}^a
& \dfrac{\textup{sgn}(x) f(x)}{|x|^\nu(\omega - x)} {\rm d}x 
\sim - \bbint{0}{a} \dfrac{\left[ (-1)^{m} g(-x)+ g(x) \right] {\rm d}x}{x^{-m+\nu+1}}
\end{split}
\end{equation}
in the same asymptotic limit, provided that the finite-part integral does not vanish. And if $g(x)$ is even symmetric, equation \eqref{splitthewhole1withbcmainabspositivesgnxgxeven} has the dominant behavior
\begin{equation} 
\label{splitthewhole1withbcmainabspositivesgncorollaryc}
\operatorname{PV}\!\!\! \int_{-a}^a
\dfrac{\textup{sgn}(x) f(x)}{|x|^\nu(\omega - x)} {\rm d}x 
\sim \begin{cases}
-2 \omega\displaystyle\bbint{0}{a} x^{-\nu-1} g(x)  
\; {\rm d}x &\text{for} \;\;\; m = 1 , \\
\\
- 2  \omega^{2\lfloor (m+1)/2 \rfloor - m} \displaystyle\int_{0}^{a} \dfrac{ g(x) \; {\rm d}x}{x^{\nu -2\lfloor (m-2)/2 \rfloor - 1}} 
\; &\text{for} \;\;\; m \geq 2 
\end{cases}
\end{equation}
in the same asymptotic limit, as long as the finite-part or the convergent integral does not vanish as well.
\end{corollary}


\section{Special Reductions of Full Hilbert Transform Integrals}
\label{otherformhilbert}


\subsection{Case $\nu = 0$}

\begin{theorem}
\label{theorem51}
Let $f(x)$ be complex analytic in the interval $[0,a]$ for some positive $a<\infty$. If $f(0)\neq 0$, then for all $\omega\in(0,a)$ satisfying $\omega<\rho_0$
\begin{equation}
\label{KKnobc1a}
\begin{split}
\operatorname{PV}\!\!\!\int_0^a 
\dfrac{f(x)}{\omega^2-x^2} {\rm d}x
=&\; - \sum_{k=0}^\infty \omega^{2k}
\bbint{0}{a} \dfrac{f(x)}{x^{2k+2}} {\rm d}x 
+ \dfrac{1}{2 \omega } 
[f(\omega) - f(-\omega)] \ln\omega .
\end{split}
\end{equation}
If $f(x)$ is even while $f(0)\neq 0$, then the principal-value integral reduces to
\begin{equation}
\label{KKnobc1aeven}
\begin{split}
\operatorname{PV}\!\!\!\int_0^a 
\dfrac{f(x)}{\omega^2-x^2} {\rm d}x
=&\; - \sum_{k=0}^\infty \omega^{2k}
\bbint{0}{a} \dfrac{f(x)}{x^{2k+2}} {\rm d}x .
\end{split}
\end{equation}
And if $f(x) = x^m g(x)$ where $m$ is a positive integer and $g(0)\neq 0$, then
\begin{equation}
\label{KKnobc1xmgx}
\begin{split}
&\operatorname{PV}\!\!\! \int_0^a 
\dfrac{f(x)}{\omega^2 - x^2} {\rm d}x
= - \sum_{k=0}^{\lfloor (m-2)/2 \rfloor}
\omega^{2k} \displaystyle\int_0^a x^{m-2k-2} g(x) \; {\rm d}x \\
&\hspace{2mm} - 
\sum_{k=0}^\infty 
\omega^{2k+2\lfloor m/2 \rfloor} \bbint{0}{a} \dfrac{x^m g(x) \; {\rm d}x}{x^{2k+2\lfloor m/2 \rfloor +2}} 
+ \dfrac{\omega^{m-1}}{2} 
\left[ g(\omega) - (-1)^m g(-\omega) \right]\ln(\omega).
\end{split}
\end{equation}
The third term in equation \eqref{KKnobc1xmgx} vanishes when $g(\omega)$ is an even function and $m$ is an even integer. Equations \eqref{KKnobc1a}-\eqref{KKnobc1xmgx}
also hold for $a=\infty$ if the principal-value integrals exist in the limit $a\rightarrow\infty$. 
\end{theorem}

\begin{proof}
To prove (\ref{KKnobc1a}), we perform partial-fraction expansion on the kernel $(\omega^2-x^2)^{-1}$ and distribute the principal-value integral to yield
\begin{equation}
\label{KKnobc3}
\begin{split}
\operatorname{PV}\!\!\!\int_0^a
\dfrac{f(x)}{\omega^2 - x^2} {\rm d}x
&= \dfrac{1}{2\omega }
\operatorname{PV}\!\!\!\int_0^a
\dfrac{f(x)}{\omega-x} {\rm d}x
+ \dfrac{1}{2\omega }
\displaystyle\int_0^a
\dfrac{f(x)}{\omega+x} {\rm d}x ,
\end{split}
\end{equation}
which is again a Stieltjes and Hilbert transform decomposition of the given principal-value integral. Substituting equations (\ref{stieltjesandhilberttransformpair2}) and (\ref{mainequationprincipalvalueorder1onepole}) back into equation (\ref{KKnobc3}), we find that the finite-part integrals $\bbint{0}{a}x^{-2k-1} f(x) \;{\rm d}x$ for all non-negative $k$ cancel, leaving only the finite-part integrals $\bbint{0}{a}x^{-2k-2} f(x) \; {\rm d}x$ to contribute. Simplification of the resulting expression yields \eqref{KKnobc1a}.

We only impose the even-parity condition of $f(x)$ in proving equation \eqref{KKnobc1aeven}. After substituting $f(x) = f(-x)$ to equation \eqref{KKnobc1a}, we retrieve equation \eqref{KKnobc1aeven}.

In order to derive equation \eqref{KKnobc1xmgx}, we substitute $f(x)=x^m g(x)$ to equation \eqref{KKnobc1a}, segregate the terms containing convergent and finite-part integrals, and apply the definition of the floor function as given in equation \eqref{splitthewhole8c}. Upon applying $g(\omega)=g(-\omega)$, we see that the third term of equation \eqref{KKnobc1xmgx} vanishes when $m$ is an even integer.

Although the separate Stieltjes and Hilbert transforms on the left-hand side of \eqref{KKnobc3} would require the more stringent condition of integrability of $f(x)x^{-1}$ at infinity than the integrability of $f(x) x^{-2}$ at infinity required by the existence of the principal-value integral in the limit $a\rightarrow\infty$, equations \eqref{KKnobc1a}, \eqref{KKnobc1aeven}, and \eqref{KKnobc1xmgx} require only integrability of $f(x) x^{-2}$ at infinity which is satisfied when the principal-value integral exists in the limit $a\rightarrow\infty$. Under this condition, equations \eqref{KKnobc1a}, \eqref{KKnobc1aeven}, and \eqref{KKnobc1xmgx} hold for $a=\infty$.
\end{proof}

\begin{corollary}
Under the same relevant conditions as in Theorem-\ref{theorem51}, equation \eqref{KKnobc1a} has the dominant behavior 
\begin{equation}
\label{KKnobc1acorollary}
\operatorname{PV}\!\!\!\int_0^a 
\dfrac{ f(x)}{\omega^2-x^2} {\rm d}x
\sim f'(0) \ln\omega 
\end{equation}
as $\omega\rightarrow 0$, provided that $f'(0)\neq 0$. When $f(x)$ is an even function, then equation \eqref{KKnobc1aeven} has the dominant behavior 
\begin{equation}
\label{KKnobc1acorollary2}
\operatorname{PV}\!\!\!\int_0^a 
\dfrac{f(x)}{\omega^2-x^2} {\rm d}x
\sim -\bbint{0}{a} \frac{f(x)}{x^2}\,\mathrm{d}x  
\end{equation}
in the same asymptotic limit, provided that the finite-part integral does not vanish. And when $f(x) = x^m g(x)$, then equation \eqref{KKnobc1xmgx} has the dominant behavior 
\begin{equation}
\label{eq5point7corollary}
\operatorname{PV}\!\!\!\int_0^a 
\dfrac{f(x)}{\omega^2-x^2} {\rm d}x
\sim \begin{cases}
g(0) \ln\omega &\text{for} \;\;\;\;\; m = 1 , \\
\\
- \displaystyle\int_0^a x^{m-2} g(x) \; {\rm d}x &\text{for} \;\;\;\;\; m \geq 2 
\end{cases}
\end{equation}
in the same asymptotic limit, as long as the integral in the last line of \eqref{eq5point7corollary} does not vanish as well.
\end{corollary}

\begin{proof}
We derive $f'(0)$ in equation \eqref{KKnobc1acorollary} using the Taylor expansion of $f(\omega)$ about $\omega=0$.
\end{proof}


\begin{theorem}
\label{theorem52}
If all relevant conditions given in Theorem-\ref{theorem51} hold, in particular $f(0)\neq 0$, then
\begin{equation}
\label{KKnobc1b}
\begin{split}
\operatorname{PV}\!\!\!\int_0^a 
\dfrac{x f(x)}{\omega^2 - x^2} {\rm d}x
=&\; - \sum_{k=0}^\infty 
\omega^{2k}
\bbint{0}{a} \dfrac{f(x)}{x^{2k+1}} {\rm d}x 
+ \dfrac{1}{2} 
[f(\omega) + f(-\omega)] \ln\omega .
\end{split}
\end{equation}
If $f(x)$ is even, then equation \eqref{KKnobc1b} reduces to
\begin{equation}
\label{KKnobc1ceven}
\begin{split}
\operatorname{PV}\!\!\!\int_0^a 
\dfrac{x f(x)}{\omega^2 - x^2} {\rm d}x
=&\; - \sum_{k=0}^\infty 
\omega^{2k}
\bbint{0}{a} \dfrac{f(x)}{x^{2k+1}} {\rm d}x + 
f(\omega) \ln\omega .
\end{split}
\end{equation}
And if $f(x) = x^m g(x)$, then equation \eqref{KKnobc1b} simplifies into
\begin{equation}
\label{KKnobc1cevenxmgx}
\begin{split}
& \operatorname{PV}\!\!\! \int_0^a 
\dfrac{x f(x)}{\omega^2 - x^2} {\rm d}x
= - \sum_{k=0}^{\lfloor (m-1)/2 \rfloor}
\omega^{2k} \displaystyle\int_0^a x^{m-2k-1} g(x) \; {\rm d}x \\
&\hspace{2mm} - 
\sum_{k=0}^\infty 
\omega^{2k+2\lfloor (m+1)/2 \rfloor} \bbint{0}{a} \dfrac{x^m g(x) \; {\rm d}x}{x^{2k+2\lfloor (m+1)/2 \rfloor + 1}} 
+ \dfrac{\omega^m}{2} 
\left[ g(\omega) + (-1)^m g(-\omega)\right]
\ln(\omega) .
\end{split}
\end{equation}
When $g(\omega)$ is even and $m$ is an odd integer, the third term in \eqref{KKnobc1cevenxmgx} vanishes.
\end{theorem}

\begin{proof} We proceed in the same manner as in Theorem-\ref{theorem51}. The only difference here is that $x^{-1} f(x)$ must be integrable at infinity for the results to hold for $a=\infty$. 
\end{proof}

\begin{corollary}
Under the same conditions as in Theorem-\ref{theorem52}, equations \eqref{KKnobc1b} and \eqref{KKnobc1ceven} have the dominant behavior
\begin{equation}
\begin{split}
\operatorname{PV}\!\!\!\int_0^a 
\dfrac{x f(x)}{\omega^2 - x^2} {\rm d}x
\sim f(0)  \ln\omega 
\end{split}
\end{equation}
as $\omega \rightarrow 0$. Similarly, equation \eqref{KKnobc1cevenxmgx} has the dominant behavior
\begin{equation}
\begin{split}
& \operatorname{PV}\!\!\! \int_0^a 
\dfrac{x f(x)}{\omega^2 - x^2} {\rm d}x
\sim - \displaystyle\int_0^a x^{m-1} g(x) \; {\rm d}x 
\end{split}
\end{equation}
in the same asymptotic limit, provided that the integral does not vanish.
\end{corollary}


\subsection{Case $\nu \neq 0$}
\begin{theorem}
\label{theorem53}
Under all relevant conditions given in Theorem-\ref{theorem51}, the following equality holds,
\begin{equation}
\label{KKwithbc1a}
\begin{split}
\operatorname{PV}\!\!\!\int_0^a 
\dfrac{f(x) \; {\rm d}x } 
{x^{\nu}(\omega^2-x^2)} 
=- \sum_{k=0}^\infty \omega^{2k}
&\; \bbint{0}{a} \dfrac{f(x)}{x^{2k+\nu+2}} {\rm d}x \\ &\; +
\frac{\pi}{2\omega^{\nu+1}} 
\left( \dfrac{f(-\omega)}{\sin(\pi\nu)}
- \dfrac{f(\omega)}{\tan(\pi\nu)}\right) ,
\end{split}
\end{equation}
for $0<\nu<1$. If $f(x)$ is even, then 
\begin{equation}
\label{KKwithbc1aeven}
\begin{split}
\operatorname{PV}\!\!\!\int_0^a
\dfrac{f(x) \; {\rm d}x}{x^\nu (\omega^2-x^2)} 
=&\; - \sum_{k=0}^\infty \omega^{2k}
\bbint{0}{a} \dfrac{f(x)}{x^{2k+\nu+2}} {\rm d}x 
+ \dfrac{\pi}{2} 
\dfrac{f(\omega)}{\omega^{\nu+1}}
\tan(\pi\nu/2) .
\end{split}
\end{equation}
And if $f(x) = x^m g(x)$, then equation \eqref{KKwithbc1a} becomes
\begin{equation}
\label{KKwithbc1axmgx}
\begin{split}
& \operatorname{PV}\!\!\! \int_0^a 
\dfrac{f(x) \; {\rm d}x}{x^\nu (\omega^2 - x^2)} 
= - \sum_{k=0}^{\lfloor (m-2)/2 \rfloor}
\omega^{2k} \displaystyle\int_0^a x^{m-2k-\nu-2} g(x) \; {\rm d}x \\
&\hspace{2mm} -
\sum_{k=0}^\infty 
\omega^{2k+2\lfloor m/2 \rfloor } \bbint{0}{a} \dfrac{x^m g(x) \; {\rm d}x}{x^{2k+2\lfloor m/2 \rfloor+\nu+2}} 
- \dfrac{\pi \omega^{m-\nu-1}}{2}
\left[ \dfrac{g(\omega)}{\tan(\pi\nu)} - \dfrac{(-1)^m g(-\omega)}{\sin(\pi\nu)} \right] .
\end{split}
\end{equation}
\end{theorem}

\begin{proof} We proceed in the same manner as in Theorem-\ref{theorem51}, only this time the decomposition into Stieltjes and one-sided Hilbert transform requires the results \eqref{stieltjesandhilberttransformpair2b} and \eqref{mainequationprincipalvalueorder1onepolewithbc}. The results hold for $a=\infty$ when the principal value exists as $a\rightarrow\infty$ or when $f(x) x^{-2-\nu}$ is integrable at infinity.
\end{proof}

\begin{corollary}
Under the same relevant conditions as in Theorem-\ref{theorem53}, equations \eqref{KKwithbc1a} and \eqref{KKwithbc1aeven} have the dominant behavior
\begin{equation}
\operatorname{PV}\!\!\!\int_0^a 
\dfrac{f(x) \; {\rm d}x } 
{x^{\nu}(\omega^2-x^2)}
\sim \frac{\pi f(0) }{2\omega^{\nu+1}} 
\tan(\pi\nu/2) 
\end{equation}
as $\omega\rightarrow 0$. If $f(x) = x^m g(x)$, then equation \eqref{KKwithbc1axmgx} has the dominant behavior
\begin{equation}
\operatorname{PV}\!\!\!\int_0^a 
\dfrac{f(x) \; {\rm d}x  } 
{x^{\nu}(\omega^2-x^2)}
\sim
\begin{cases} 
- \dfrac{\pi}{2\omega^\nu} g(0) \cot(\pi\nu/2) &\text{for} \;\;\; m = 1 , \\
\\
- \displaystyle\int_{0}^{a} x^{m-\nu-2} g(x) \; {\rm d}x &\text{for} \;\;\; m \geq 2 
\end{cases}
\end{equation}
in the same asymptotic limit, provided that the integral in the last line does not vanish.
\end{corollary}


\begin{theorem}
\label{theorem54}
Under all relevant conditions given in Theorem-\ref{theorem53}, the following equality holds,
\begin{equation}
\label{KKwithbc1b}
\begin{split}
\operatorname{PV}\!\!\!\int_0^a 
& \dfrac{x^{1-\nu} f(x)}{\omega^2 - x^2} {\rm d}x \\
&= - \sum_{k=0}^\infty \omega^{2k}
\bbint{0}{a} \dfrac{f(x)}{x^{2k+\nu+1}} {\rm d}x
- \frac{\pi}{2\omega^\nu} 
\left( \dfrac{f(\omega)}{\tan(\pi\nu)}
+ \dfrac{f(-\omega)}{\sin(\pi\nu)} \right)  .
\end{split}
\end{equation}
If $f(x)$ is even, then
\begin{equation}
\label{KKwithbc1bceven}
\begin{split}
\operatorname{PV}\!\!\!\int_0^a 
& \dfrac{x^{1-\nu} f(x)}{\omega^2 - x^2} {\rm d}x 
= - \sum_{k=0}^\infty \omega^{2k}
\bbint{0}{a} \dfrac{f(x)}{x^{2k+\nu+1}} {\rm d}x
- \dfrac{\pi}{2}
\dfrac{f(\omega)}{\omega^\nu} 
\cot(\pi\nu/2) .
\end{split}
\end{equation}
And if $f(x) = x^m g(x)$, then equation \eqref{KKwithbc1b} simplifies into
\begin{equation}
\label{KKwithbc1bcxmgx}
\begin{split}
&\operatorname{PV}\!\!\! \int_0^a 
\dfrac{x^{1-\nu} f(x)}{\omega^2 - x^2} {\rm d}x
= - \sum_{k=0}^{\lfloor (m-1)/2 \rfloor}
\omega^{2k} \displaystyle\int_0^a x^{m-2k-\nu-1} g(x) \; {\rm d}x \\
&\hspace{2mm} -
\sum_{k=0}^\infty 
\omega^{2k+2\lfloor (m+1)/2\rfloor} \bbint{0}{a} \dfrac{x^m g(x) \; {\rm d}x}{x^{2k+2\lfloor (m+1)/2\rfloor+\nu+1}} 
- \dfrac{\pi \omega^{m-\nu}}{2}
\left[ \dfrac{g(\omega)}{\tan(\pi\nu)} + \dfrac{(-1)^m g(-\omega)}{\sin(\pi\nu)} \right] .
\end{split}
\end{equation}
\end{theorem}

\begin{proof} We proceed in the same manner as in Theorem-\ref{theorem53}. The results hold for $a=\infty$ when the principal value exists as $a\rightarrow\infty$ or when $f(x) x^{-1-\nu}$ is integrable at infinity.
\end{proof}

\begin{corollary}
Under the same conditions as in Theorem-\ref{theorem54}, equations \eqref{KKwithbc1b} and \eqref{KKwithbc1bceven} have the dominant behavior 
\begin{equation}
\begin{split}
\operatorname{PV}\!\!\!\int_0^a 
& \dfrac{x^{1-\nu} f(x)}{\omega^2 - x^2} {\rm d}x 
\sim - \dfrac{\pi}{2}
\dfrac{f(0)}{\omega^\nu} 
\cot(\pi\nu/2)
\end{split}
\end{equation}
as $\omega \rightarrow 0$. For the case when $f(x) = x^m g(x)$, equation \eqref{KKwithbc1bcxmgx} has the dominant behavior
\begin{equation}
\operatorname{PV}\!\!\!\int_0^a 
\dfrac{x^{1-\nu} f(x)}{\omega^2 - x^2} {\rm d}x 
\sim -\displaystyle\int_0^a x^{m-\nu-1} g(x) \; {\rm d}x
\end{equation}
in the same asymptotic limit, provided that the integral does not vanish.
\end{corollary}


\section{Examples}
\label{examplesection}


\subsection{Example 1}
\label{example1}

\hspace{\parindent} In this example, we use finite-part integration to reproduce the known full Hilbert transform,
\begin{equation}
\label{oddevenexample1eq10ex1}
\begin{split}
\operatorname{PV}\!\!\!\int_{-\infty}^\infty
\dfrac{e^{iax}}{\omega - x} {\rm d}x
= &\; -i \pi \; \textup{sgn}(a) e^{i a \omega} ,
\end{split}
\end{equation}
for all real $a\neq 0$. This is the Hilbert transform of the function $e^{i a x}$ whose complex extension is entire. The relevant equation is given by (\ref{splitthewhole6maintheorem}), and its application to the given Hilbert transform \eqref{oddevenexample1eq10ex1} yields
\begin{equation}
\label{oddevenexample1eq2ex1}
\begin{split}
\operatorname{PV}\!\!\!\int_{-\infty}^\infty
\dfrac{e^{iax}}{\omega - x} {\rm d}x 
= - 
\sum_{k=0}^\infty \omega^k
\bbint{0}{\infty} 
\dfrac{e^{iax}}{x^{k+1}} {\rm d}x 
+ 
\sum_{k=0}^\infty (-\omega)^k
\bbint{0}{\infty} 
\dfrac{e^{-iax}}{x^{k+1}} {\rm d}x ,
\end{split}
\end{equation}
where the two series expansions in \eqref{oddevenexample1eq2ex1} converge for all real $\omega$. In  Section-\ref{fpiexample1}, we derive the finite-part integral
\begin{equation}
\label{oddevenexample1eq2xxex1}
 \bbint{0}{\infty} \dfrac{e^{iax}}{x^{k+1}} {\rm d}x
= -\dfrac{(ia)^k}{k!}
\left[ 
\ln|a| - \dfrac{i\pi}{2}\operatorname{sgn}(a) - \psi(k+1)
\right], \;\; k=0, 1, 2, \cdots    
\end{equation}
for all real $a\neq 0$. Substituting equation \eqref{oddevenexample1eq2xxex1} back into \eqref{oddevenexample1eq2ex1}, we obtain the two infinite sums
\begin{equation}
\label{xx1ex1}
\begin{split}
\sum_{k=0}^{\infty} 
\omega^k \bbint{0}{\infty} \frac{e^{i a x}}{x^{k+1}}\,{\rm d}x 
&= -\sum_{k=0}^{\infty}\omega^k \frac{(i a )^k}{k!}\left(\ln|a|-i\frac{\pi}{2}\operatorname{sgn}(a) - \psi(k+1)\right)\\
&= \left(i\frac{\pi}{2}\operatorname{sgn}a-\ln|a|\right)e^{i a \omega} + \sum_{k=0}^{\infty} \omega^k \frac{(ia)^k}{k!} \psi(k+1),
\end{split}
\end{equation}
\begin{eqnarray}
\label{xx2ex1}
\sum_{k=0}^{\infty} (-\omega)^k \bbint{0}{\infty} \frac{e^{-i a x}}{x^{k+1}}\,{\rm d}x \!\!\! &=&\!\!\!-\sum_{k=0}^{\infty}(-1)^k\omega^k \frac{(-i a )^k}{k!}\left(\ln|a|-i\frac{\pi}{2}\operatorname{sgn}(-a) - \psi(k+1)\right)\nonumber\\
&=&\!\!\!\left(-i\frac{\pi}{2}\operatorname{sgn}a-\ln|a|\right)e^{i a \omega} + \sum_{k=0}^{\infty} \omega^k \frac{(ia)^k}{k!} \psi(k+1) .
\end{eqnarray}
Substituting equations \eqref{xx1ex1} and \eqref{xx2ex1} back into equation \eqref{oddevenexample1eq2ex1} leads to the known result \eqref{oddevenexample1eq10ex1}.


\subsection{Example 2}
\label{example2}

\hspace{\parindent} We obtain here the alternative method of evaluating
\begin{equation}
\label{oddevenexample1eq10}
\begin{split}
\operatorname{PV}\!\!\!\int_{-\infty}^\infty
\dfrac{\operatorname{sech}^2(ax)}{\omega - x} {\rm d}x 
\end{split}
\end{equation}
for all real $a$ and for $|\omega|<\rho_0=\pi/(2a)$. The integral in \eqref{oddevenexample1eq10} was considered in \cite{sechsquared}. Since $f(x)$ is even symmetric with $f(0)\neq 0$, then we use equation \eqref{splitthewhole6maintheoremeven} to evaluate equation \eqref{oddevenexample1eq10}. In doing so, we write the principal-value integral as
\begin{equation}
\label{oddevenexample1eq2}
\begin{split}
\operatorname{PV}\!\!\!\int_{-\infty}^\infty
\dfrac{\operatorname{sech}^2(ax)}{\omega - x} {\rm d}x 
= - 
\sum_{k=0}^\infty \omega^{2k+1}
\bbint{0}{\infty} 
\dfrac{\operatorname{sech}^2(ax)}{x^{2k+2}} {\rm d}x .
\end{split}
\end{equation}
Note that $\operatorname{sech}^2(ax)$ has singularity points at $z=\pm in \pi/(2a)$ for $n=1,3,5,7,\cdots$. It then follows that the series expansion in \eqref{oddevenexample1eq2} converges for any real $|\omega|<\pi/(2a)$. In Section-\ref{fpiexample2} we derive the finite-part integral
\begin{equation}
\label{oddevenexample1eq2xx}
 \bbint{0}{\infty} \dfrac{\operatorname{sech}^2(ax)}{x^{2k+2}} {\rm d}x
= (-1)^{k+1} \dfrac{2k+2}{\pi^{2k+2}}a^{2k+1}(2^{2k+3}-1)\zeta(2k+3)   
\end{equation}
for $k=0, 1, 2, \cdots $. If we substitute equation \eqref{oddevenexample1eq2xx} to \eqref{oddevenexample1eq2}, then we obtain 
\begin{equation}
\label{xx1}
\begin{split}
\operatorname{PV}\!\!\!\int_{-\infty}^\infty
\dfrac{\operatorname{sech}^2(ax)}{\omega - x} {\rm d}x
= \dfrac{4}{\pi}\sum_{k=0}^\infty (-1)^k \left(\dfrac{a\omega}{\pi}\right)^{2k+1} (k+1) (2^{2k+3}-1) \zeta(2k+3) ,
\end{split}
\end{equation}
which is the Hilbert transform for $\operatorname{sech}^2(ax)$ for $|\omega|<\pi/(2a)$.


\subsection{Example 3}
\label{example3}

\hspace{\parindent}We demonstrate how new series representations for $K_0(x)$ can be obtained using finite-part integration. Consider
\begin{equation}
\label{ex2eq1a}
\operatorname{PV}\!\!\!\int_{0}^\infty 
\dfrac{e^{-x}I_0(x)}{\omega-x} {\rm d}x
= - e^{-\omega}K_0(\omega), \;\;\;\;\;
\omega>0 ,
\end{equation}
which is given by equation (8H.3) in \cite[p.~504]{hilbertbookking_2009b}. Since the principal-value integral in \eqref{ex2eq1a} is one-sided, then we use equation \eqref{mainequationprincipalvalueorder1onepole} to rewrite the left-hand side of equation \eqref{ex2eq1a} as
\begin{equation}
\label{ex2eq1}
\displaystyle\int_{0}^\infty 
\dfrac{e^{-x}I_0(x)}{\omega-x} {\rm d}x
= -\sum_{k=0}^\infty \omega^k 
\bbint{0}{\infty }\dfrac{e^{-x}I_0(x)}{x^{k+1}} {\rm d}x + e^{-\omega}I_0(\omega)\ln(\omega) .
\end{equation}
The series in \eqref{ex2eq1} is convergent for all $\omega>0$ because $e^{-x}I_0(x)$ is entire in the complex plane. In Section-\ref{fpiexample3} we derive the finite-part integral
\begin{equation}
\label{ex2eq2}
\begin{split}
\bbint{0}{\infty} 
& \dfrac{e^{-ax}I_n(ax)}{x^{k+n+1}} {\rm d}x
= (-1)^{k+1} \dfrac{(2a)^{k+n}}{\sqrt{\pi}} \dfrac{\Gamma(k+n+1/2)}{k!(k+2n)!} \\
&\hspace{0.4cm}
\times \left[\ln(2a) + \psi(k+n+1/2)-\psi(k+1)-\psi(k+2n+1) \right] .
\end{split}
\end{equation}
And when we use equation \eqref{ex2eq2} to \eqref{ex2eq1}, we obtain  
\begin{equation}
\label{ex2eq2mellin4}
\begin{split}
e^{-x} 
& K_0(x)
= - e^{-x}I_0(x)\ln(2x) \\
&\hspace{0.2cm}
-\dfrac{1}{\sqrt{\pi}}
\sum_{k=0}^\infty
(-1)^k \dfrac{(2x)^k \Gamma(k+1/2)}{(k!)^2}[\psi(k+1/2)-2\psi(k+1)] 
\end{split}
\end{equation}
for $x>0$. 

We do similar steps leading to \eqref{ex2eq2mellin4} for
\begin{equation}
\label{ex2eq2mellin4b}
\displaystyle\int_{-\infty}^0 \dfrac{e^{x}I_0(x)}{\omega-x} {\rm d}x
= e^{\omega}K_0(\omega),
\;\;\;\omega>0 ,
\end{equation}
where \eqref{ex2eq2mellin4b} is given in equation (8H.4) in \cite[p.~504]{hilbertbookking_2009b}. However, we need to implement the change of variable $x\rightarrow -x$ in equation \eqref{ex2eq2mellin4b} so that it becomes
\begin{equation}
\label{ex2eq2mellin4c}
\displaystyle\int_0^{\infty} \dfrac{e^{-x}I_0(x)}{\omega + x} {\rm d}x
= e^{\omega}K_0(\omega) .
\end{equation}
This time, the integral involved in \eqref{ex2eq2mellin4b} is a Stieltjes transform. So when we use \eqref{stieltjesandhilberttransformpair2}, then we have
\begin{equation}
\label{ex2eq2mellin4d}
\displaystyle\int_0^{\infty} \dfrac{e^{-x}I_0(x)}{\omega + x} {\rm d}x
= \sum_{k=0}^\infty (-1)^k \omega^k 
\bbint{0}{\infty }\dfrac{e^{-x}I_0(x)}{x^{k+1}} {\rm d}x - e^{\omega}I_0(\omega)\ln(\omega) .
\end{equation}
And after we substitute equation \eqref{ex2eq2} into equation \eqref{ex2eq2mellin4d}, we obtain
\begin{equation}
\label{ex2eq2mellin5}
\begin{split}
e^{x} 
& K_0(x)
= - e^{x} I_0(x)\ln(2x) \\
&\hspace{0.2cm}
-\dfrac{1}{\sqrt{\pi}}
\sum_{k=0}^\infty
\dfrac{(2x)^k \Gamma(k+1/2)}{(k!)^2}[\psi(k+1/2)-2\psi(k+1)] 
\end{split}
\end{equation}
for $x>0$, which is identical to equation (93) in \cite[p.~15]{StieltjesFinitePartpapermain13}.

Equation \eqref{ex2eq2mellin4} is a new series representation for $K_0(x)$ that is derived using finite-part integration, as they do not appear in the table of series representations given in \cite{Handbook1, Handbook2, Handbook3, Handbook4, Handbook5} and cannot be evaluated using Mathematica software. In addition to that, we can obtain other new series representations from the linear combinations of equations \eqref{ex2eq2mellin4} and \eqref{ex2eq2mellin5}. Recall that $\sinh(x)$ and $\cosh(x)$ are defined as
\begin{equation}
\label{sinhcosh}
\sinh(x)
= \dfrac{e^x-e^{-x}}{2},\;
\cosh(x)=\dfrac{e^x+e^{-x}}{2}.
\end{equation}
Taking the linear combinations of equations \eqref{ex2eq2mellin4} and \eqref{ex2eq2mellin5} leads to
\begin{equation}
\label{ex2eq5}
\begin{split}
\sinh(x) 
& K_0(x)
= - \sinh(x)I_0(x) \ln(2x) \\
&\hspace{0.1cm} -\dfrac{1}{\sqrt{\pi}}
\sum_{k=0}^\infty \dfrac{(2x)^{2k+1}\Gamma(2k+3/2)}{((2k+1)!)^2}[\psi(2k+3/2)-2\psi(2k+2)] ,
\end{split}
\end{equation}
\begin{equation}
\label{ex2eq6}
\begin{split}
\cosh(x) 
& K_0(x)
= - \cosh(x)I_0(x) \ln(2x) \\
&\hspace{0.1cm} -\dfrac{1}{\sqrt{\pi}}
\sum_{k=0}^\infty \dfrac{(2x)^{2k}\Gamma(2k+1/2)}{((2k)!)^2}[\psi(2k+1/2)-2\psi(2k+1)] ,
\end{split}
\end{equation}
\begin{equation}
\label{ex2eq7}
\begin{split}
& K_0(x)
= - I_0(x)\ln(2x) \\
&\hspace{0.4cm}
-\dfrac{\cosh(x)}{\sqrt{\pi}}
\sum_{k=0}^\infty
\dfrac{(2x)^{2k} \Gamma(2k+1/2)}{((2k)!)^2}[\psi(2k+1/2)-2\psi(2k+1)] \\
&\hspace{0.4cm}
+\dfrac{\sinh(x)}{\sqrt{\pi}}
\sum_{k=0}^\infty
\dfrac{(2x)^{2k+1} \Gamma(2k+3/2)}{((2k+1)!)^2}[\psi(2k+3/2)-2\psi(2k+2)] .
\end{split}
\end{equation}

Equations \eqref{ex2eq2mellin4}, \eqref{ex2eq2mellin5}, \eqref{ex2eq5} and \eqref{ex2eq6} also lead to the new series representations for some Meijer-G functions. In particular, references \cite{InternetSourceMeijerGx1, InternetSourceMeijerGx2, InternetSourceMeijerGx3} provide the following equalities:
\begin{align}
\setlength\arraycolsep{1pt}
e^{-x}K_0(x)
&= \sqrt{\pi} G_{1,2}^{2,0}
\left( \begin{matrix} 2x  \end{matrix} \; \Bigg| \; 
\begin{matrix}
1/2 \\ 0,\; 0
\end{matrix}
\right) , 
\label{MeijerG1} \\
e^{x}K_0(x)
&= \dfrac{1}{\sqrt{\pi}} G_{1,2}^{2,1}
\left( \begin{matrix} 2x \end{matrix} \; \bigg| \; 
\begin{matrix}
1/2 \\ 0,\; 0
\end{matrix}
\right) , \\
\sinh(x)K_0(x)
&= \dfrac{1}{2\sqrt{2}} G_{2,4}^{2,2}
\left( \begin{matrix} x, \dfrac{1}{2} \end{matrix} \; \bigg| \; 
\begin{matrix}
1/4,\; 3/4 \\ 1/2,\; 1/2,\; 0,\; 0
\end{matrix}
\right) , \\
\cosh(x)K_0(x)
&= \dfrac{1}{2\sqrt{2}} G_{2,4}^{2,2}
\left( \begin{matrix} x, \dfrac{1}{2} \end{matrix} \; \bigg| 
\begin{matrix}
1/4,\; 3/4 \\ 0,\; 0,\; 1/2,\; 1/2
\end{matrix}
\right) .
\end{align}
When we compare equations \eqref{ex2eq2mellin4} and \eqref{MeijerG1}, we obtain
\begin{equation}
\begin{split}
G_{1,2}^{2,0}
\left( \begin{matrix} x  \end{matrix} \; \Bigg| \; 
\begin{matrix}
1/2 \\ 0,\; 0
\end{matrix}
\right)
&= -\dfrac{1}{\sqrt{\pi}} e^{-x/2}I_0(x/2)\ln(x) \\
&\hspace{-0.8cm} 
- \dfrac{1}{\pi}
\sum_{k=0}^\infty (-1)^k x^k \dfrac{ \Gamma(k+1/2)}{(k!)^2} [\psi(k+1/2)-2\psi(k+1)] .
\end{split}
\end{equation}
And we do the same technique for the other given Meijer-G functions to obtain
\begin{equation}
\begin{split}
G_{1,2}^{2,1}
\left( \begin{matrix} x \end{matrix} \; \bigg| \; 
\begin{matrix}
1/2 \\ 0,\; 0
\end{matrix}
\right) 
&= -\sqrt{\pi} e^{x/2}I_0(x/2)\ln(x) \\
&\hspace{-0.6cm} 
- \sum_{k=0}^\infty x^k \dfrac{ \Gamma(k+1/2)}{(k!)^2} [\psi(k+1/2)-2\psi(k+1)] ,
\end{split}
\end{equation}
\begin{equation}
\begin{split}
G_{2,4}^{2,2}
\left( \begin{matrix} x, \dfrac{1}{2} \end{matrix} \; \bigg| \; 
\begin{matrix}
1/4,\; 3/4 \\ 1/2,\; 1/2,\; 0,\; 0
\end{matrix}
\right)
&= -2\sqrt{2} \sinh(x)I_0(x)\ln(2x) \\
&\hspace{-2.9cm}
-2\sqrt{\dfrac{2}{\pi}}
\sum_{k=0}^\infty \dfrac{(2x)^{2k+1}\Gamma(2k+3/2)}{((2k+1)!)^2}[\psi(2k+3/2)-2\psi(2k+2)] ,
\end{split}
\end{equation}
\begin{equation}
\begin{split}
G_{2,4}^{2,2}
\left( \begin{matrix} x, \dfrac{1}{2} \end{matrix} \; \bigg| \; 
\begin{matrix}
1/4,\; 3/4 \\ 0,\; 0,\; 1/2,\; 1/2
\end{matrix}
\right)
&= -2\sqrt{2} \cosh(x)I_0(x)\ln(2x) \\
&\hspace{-2.9cm}
-2\sqrt{\dfrac{2}{\pi}}
\sum_{k=0}^\infty \dfrac{(2x)^{2k}\Gamma(2k+1/2)}{((2k)!)^2}[\psi(2k+1/2)-2\psi(2k+1)] .
\end{split}
\end{equation}


\section{Mellin Transform Method in Evaluating Finite-Part Integrals}
\label{mellintransformmethod}

\hspace{\parindent}In principle, the finite-part integrals can be obtained using the canonical method, that is, as prescribed by the definitions \eqref{canonical} and \eqref{introfpirelevantconcept4}. However, the definitions may not be convenient in extracting the finite-part of a divergent integral. A powerful method is established in \cite{FPIRegularized} and is based on the Mellin transform
\begin{equation}
    \mathcal{M}_a[f(x); 1-\lambda]=\int_0^a x^{-\lambda} f(x)\,\mathrm{d}x, \;\;\; 0<a\leq\infty . 
\end{equation}
At most, the transform exists in a bounded or unbounded strip along the imaginary axis. The relationship between the Mellin transform and the finite-part integral
\begin{equation}
    \bbint{0}{a}\frac{f(x)}{x^{\lambda}}\, \mathrm{d}x, \;\;\;\operatorname{Re}\lambda\geq  1, 
\end{equation}
is through the analytic continuation of the Mellin transform, which is denoted by
\begin{equation}
    \mathcal{M}_a^*[f(x);1-\lambda] .
\end{equation}
Under the condition that $f(x)$ is analytic at the origin, the Mellin transform has at most simple poles along the real axis.

In \cite{FPIRegularized} it is established that the value of the analytic continuation is equal to that of the finite-part integral if $\lambda_0$ is a regular point of the analytic continuation,
\begin{equation}
    \bbint{0}{a}\frac{f(x)}{x^{\lambda_0}}\, \mathrm{d}x=\mathcal{M}^*_a[f(x);1-\lambda_0].
\end{equation}
If $\lambda_0$ happens to be  simple pole of the analytic continuation, then the finite-part integral is equal to the regularized limit of $\mathcal{M}^*_a[f(x);1-\lambda]$ at $\lambda_0$,
\begin{equation}
    \bbint{0}{a}\frac{f(x)}{x^{\lambda_0}}\,\mathrm{d}x = \reglim{\lambda}{\lambda_0}\mathcal{M}^*_a[f(x);1-\lambda] .
\end{equation}
If $\lambda_0$ is an isolated singularity of $\mathcal{M}^*_a[f(x);1-\lambda]$, then the regularized limit of $\mathcal{M}^*_a[f(x);1-\lambda]$ at $\lambda_0$ is, by definition, the value of the regular part of the Laurent series expansion of $\mathcal{M}^*_a[f(x);1-\lambda]$ in a deleted neighborhood of $\lambda_0$. If $\mathcal{M}^*_a[f(x);1-\lambda]=F(\lambda)/G(\lambda)$, where $F(\lambda)$ and $G(\lambda)$ are both analytic at $\lambda_0$ with $F(\lambda_0)\neq 0$ and $G(\lambda_0)=0$, the regularized limit at $\lambda_0$ is given by
\begin{equation}
\label{mellinsimplepole}
    \reglim{\lambda}{\lambda_0}\frac{F(\lambda)}{G(\lambda)}=\frac{F'(\lambda_0)}{G'(\lambda_0)}-\frac{F(\lambda_0) G''(\lambda_0)}{2 (G'(\lambda_0))^2} ,
\end{equation}
provided $\lambda_0$ is a simple zero of $G(\lambda)$. If it happens that $\lambda_0$ is a removable singularity of $\mathcal{M}^*_a[f(x);1-\lambda]$, then the second term in equation \eqref{mellinsimplepole} vanishes and the regularized limit reduces to the Cauchy limit, which is given by
\begin{equation}
\label{defCauchyLimit}
    \lim_{\lambda\rightarrow \lambda_0} \frac{F(\lambda)}{G(\lambda)}=\frac{F'(\lambda_0)}{G'(\lambda_0)} .
\end{equation}

We point out that not all functions posses Mellin transform so that the finite-part integral cannot be computed using Mellin transform. For such cases, the finite-part can be extracted using the definition or the contour integral representation of the finite-part integral.


\subsection{Example 1}\label{fpiexample1}

\label{section71}
\begin{figure}[t!]
\centering\includegraphics[width=0.26\linewidth]{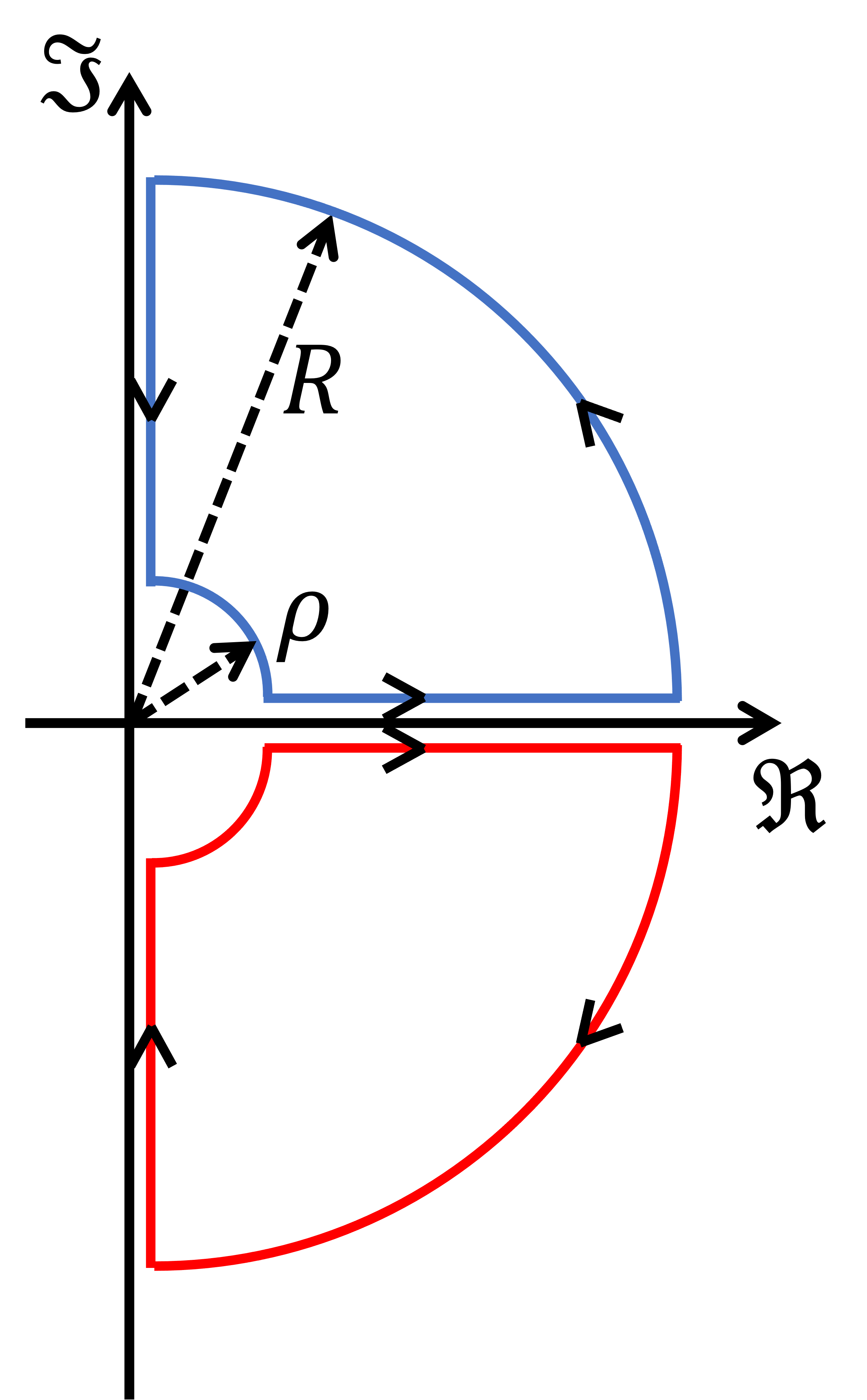}
\caption{The contours of integration that evaluate the Mellin transform for $e^{iax}$. The upper contour in this figure is used for $a>0$, whereas the lower contour is chosen for $a<0$.}
\label{contourmellinexponential}
\end{figure}

\hspace{\parindent} In this example, we use the Mellin transform and the canonical methods in evaluating $\bbint{0}{\infty} x^{-\lambda} e^{iax} \; {\rm d}x$ for all real $a\neq 0$ and for any real $\lambda \geq 1$. Let us begin with the Mellin transform method in extracting the finite-part. Using the contours shown in Figure-\ref{contourmellinexponential}, we obtain 
\begin{equation}
\label{ex1sec2eq10}
\mathcal{M}[e^{iax};1-\lambda]
= \int_{0}^{\infty} \dfrac{e^{iax}}{x^\lambda} {\rm d}x
= \dfrac{|a|^{\lambda-1}}{\Gamma(\lambda)} \dfrac{\pi}{\sin(\pi\lambda)} \exp\left( -\dfrac{i\pi(\lambda-1)}{2} \operatorname{sgn}(a) \right) 
\end{equation}
for $0<\operatorname{Re} \lambda<1$. The left-hand side of equation \eqref{ex1sec2eq10} extends analytically in the entire complex plane, then we have the desired analytic continuation
\begin{equation}
\label{mellinstarexample1}
\mathcal{M}^*[e^{iax};1-\lambda]
= \dfrac{|a|^{\lambda-1}}{\Gamma(\lambda)} \dfrac{\pi}{\sin(\pi\lambda)} \exp\left( -\dfrac{i\pi(\lambda-1)}{2} \operatorname{sgn}(a) \right) ,
\end{equation}
where \eqref{mellinstarexample1} is analytic for any non-integer $\lambda>1$ and has simple poles for any integer $\lambda \geq 1$. It then follows that the finite-part of $\bbint{0}{\infty} x^{-\lambda}e^{iax}\;{\rm d}x$ is
\begin{equation}
\label{ex1sec2eq10b}
\bbint{0}{\infty} \dfrac{e^{iax}}{x^\lambda} {\rm d}x
= \dfrac{|a|^{\lambda-1}}{\Gamma(\lambda)} \dfrac{\pi}{\sin(\pi\lambda)} \exp\left( -\dfrac{i\pi(\lambda-1)}{2} \operatorname{sgn}(a) \right)
\end{equation}
for all non-integer $\lambda > 1$. And we can simplify equation \eqref{ex1sec2eq10b} by letting $\lambda = m+\nu$, for $0<\nu<1$ and $m=1,2,3,\cdots$ to obtain
\begin{equation}
\label{fpiexpmellinwithnu}
\begin{split}
\bbint{0}{\infty} \dfrac{e^{iax}}{x^{m+\nu}} {\rm d}x
= |a|^{m+\nu-1} \Gamma(1-m-\nu)  \exp\left( -\dfrac{i\pi(m+\nu-1)}{2} \operatorname{sgn}(a) \right) .
\end{split}
\end{equation}

In order to evaluate $\bbint{0}{\infty} x^{-k-1} e^{iax} {\rm d}x$ for $k=0, 1, 2, \cdots$, we use the regularized limit because the analytic continuation \eqref{mellinstarexample1} has simple poles at any integer $\lambda\geq 1$. We begin by writing equation \eqref{ex1sec2eq10} as
\begin{equation}
\label{ex1sec2eq10cx}
\mathcal{M}^*[e^{iax};1-\lambda] = \dfrac{F(\lambda)}{G(\lambda)}
\end{equation}
where
\[
F(\lambda) = \pi \dfrac{|a|^{\lambda-1}}{\Gamma(\lambda)}  \exp\left( -\dfrac{i\pi(\lambda-1)}{2} \operatorname{sgn}(a) \right) , \;
G(\lambda) = \sin(\pi\lambda) .
\]
As described in equation \eqref{mellinsimplepole}, we can obtain $\bbint{0}{\infty} x^{-k-1} e^{iax} \;{\rm d}x$
through
\begin{equation}
\label{ex1sec2eq10c}
\begin{split}
\bbint{0}{\infty}
\dfrac{e^{iax}}{x^{k+1}} {\rm d}x 
&= \reglim{\lambda}{k+1}\mathcal{M}^*[e^{iax};1-\lambda] \\
&= \dfrac{F'(k+1)}{G'(k+1)}
- \dfrac{F(k+1)}{2} \dfrac{G''(k+1)}{(G'(k+1))^2} .
\end{split}
\end{equation}
Observe that the choice for $G(\lambda)$ leads to the vanishing of $G''(k+1)$ so that only the first term contributes. Executing equation \eqref{ex1sec2eq10c} leads to
\begin{equation}
\label{ex1sec2eq8}
\bbint{0}{\infty} \dfrac{e^{iax}}{x^{k+1}} {\rm d}x
= -\dfrac{(ia)^k}{k!}
\left[ 
\ln|a| - \dfrac{i\pi}{2}\operatorname{sgn}(a) - \psi(k+1)
\right]. 
\end{equation}

Another way of extracting the finite-part for $\bbint{0}{\infty} x^{-\lambda} e^{iax} \; {\rm d}x$ is by using the definition given by \eqref{canonical}. We show here in detail for $\lambda=m+\nu$ case, where $m=1,2,3,\cdots$ and $0<\nu<1$. Using the definition, we temporarily assign the lower limit of the divergent integral as $\varepsilon>0$ and the upper limit of integration to be $s>\varepsilon$ such as
\begin{equation}
\label{example2fullhilbertwithBC3}
\displaystyle\int_{\varepsilon}^{s} 
\dfrac{e^{iax}}{x^{m+\nu}} {\rm d}x 
\end{equation}
before we take the limit $s\rightarrow\infty$. If we use the series expansion of $e^{iax}$, then equation \eqref{example2fullhilbertwithBC3} becomes
\begin{equation}
\label{example2fullhilbertwithBC3b}
\begin{split}
\displaystyle\int_\varepsilon^s
\dfrac{e^{iax}}{x^{m+\nu}} dx
&= \sum_{n=0}^\infty 
\dfrac{(ia)^n}{n!} \displaystyle\int_\varepsilon^s
x^{n-m-\nu} {\rm d}x .
\end{split}
\end{equation}
There are two cases to consider in evaluating \eqref{example2fullhilbertwithBC3b}: $0 \leq n \leq m-1 $ and $n \geq m$. Using the two cases, we write the integral result as
\begin{equation}
\label{example2fullhilbertwithBC5}
\begin{split}
\displaystyle\int_\varepsilon^s
& \dfrac{e^{iax}}{x^{m+\nu}} {\rm d}x
= \sum_{n=0}^{m-1} 
\dfrac{(ia)^n}{n!} \dfrac{1}{n-m-\nu+1}
\dfrac{1}{s^{m+\nu-n-1}}
+ \sum_{n=m}^\infty
\dfrac{(ia)^n}{n!} \dfrac{s^{n-m-\nu+1}}{n-m-\nu+1}
 \\
& - \sum_{n=0}^{m-1} 
\dfrac{(ia)^n}{n!} \dfrac{1}{n-m-\nu+1}
\dfrac{1}{\varepsilon^{m+\nu-n-1}}
- \sum_{n=m}^\infty 
\dfrac{(ia)^n}{n!} \dfrac{\varepsilon^{n-m-\nu+1}}{n-m-\nu+1}
.
\end{split}
\end{equation}
Equation (\ref{example2fullhilbertwithBC5}) gives us the diverging term
\[
D_\varepsilon
= - \sum_{n=0}^{m-1} 
\dfrac{(ia)^n}{n!} \dfrac{1}{n-m-\nu+1}
\dfrac{1}{\varepsilon^{m+\nu-n-1}} .
\]
as $\varepsilon \rightarrow 0$. As a result, the finite-part of $\bbint{0}{s} x^{-m-\nu} e^{iax} \; {\rm d}x$ is
\begin{equation}
\label{example2fullhilbertwithBC6}
\begin{split}
\bbint{0}{s}
 \dfrac{e^{iax}}{x^{m+\nu}} {\rm d}x
= \sum_{n=0}^{m-1} 
\dfrac{(ia)^n}{n!} \dfrac{1}{n-m-\nu+1}
\dfrac{1}{s^{m+\nu-n-1}}
+ \sum_{n=0}^\infty
\dfrac{(ia)^{n+m}}{(n+m)!} \dfrac{s^{n-\nu+1}}{n-\nu+1} .
\end{split}
\end{equation}

We can further simplify equation (\ref{example2fullhilbertwithBC6}) by taking hypergeometric summation. We have the sum
\begin{equation}
\label{example2fullhilbertwithBC7}
\begin{split}
\sum_{n=0}^\infty
\dfrac{z^n}{(n-\nu+1)(n+m)!} 
= \dfrac{1}{(1-\nu) m!}
\times 
\pFq{2}{2}{1, 1-\nu}{1+m, 2-\nu}{z} .
\end{split}
\end{equation}
So when we take the limit of \eqref{example2fullhilbertwithBC6} as $s\rightarrow\infty$, we obtain
\begin{equation}
\label{example2fullhilbertwithBC7b}
\begin{split}
\bbint{0}{\infty}
\dfrac{e^{iax}}{x^{m+\nu}} {\rm d}x
= \dfrac{(ia)^m }{(1-\nu)m!}
\lim_{s\rightarrow\infty} 
s^{1-\nu} 
\pFq{2}{2}{1,1-\nu}{1+m, 2-\nu}{ias} .
\end{split}
\end{equation}
Reference \cite{InternetSource1} gives the asymptotic behavior
\[
\begin{split}
\pFq{2}{2}{a_1, a_2}{b_1, b_2}{z}
\approx &\;
\dfrac{\Gamma(a_2-a_1)\Gamma(b_1)\Gamma(b_2)}
{\Gamma(a_2)\Gamma(b_1-a_1)\Gamma(b_2-a_1)}
(-z)^{-a_1} (1 + \mathcal{O}(z^{-1})) \\
&\; +\dfrac{\Gamma(a_1-a_2)\Gamma(b_1)\Gamma(b_2)}
{\Gamma(a_1)\Gamma(b_1-a_2)\Gamma(b_2-a_2)}
(-z)^{-a_2} (1 + \mathcal{O}(z^{-1})) \\
&\; +\dfrac{\Gamma(b_1)\Gamma(b_2)}
{\Gamma(a_1)\Gamma(a_2)} 
z^{a_1+a_2-b_1-b_2} e^z 
(1 + \mathcal{O}(z^{-1})) 
\end{split}
\]
for $|z|\rightarrow\infty$. As a result, equation (\ref{example2fullhilbertwithBC6}) becomes
\begin{equation}
\label{example2fullhilbertwithBC9}
\begin{split}
\bbint{0}{\infty}
\dfrac{e^{iax}}{x^{m+\nu}} dx
= &\; (-1)^m (-ia)^{m+\nu-1} \dfrac{\Gamma(\nu)\Gamma(1-\nu)}
{\Gamma(m+\nu)} .
\end{split}
\end{equation}
Finally, we can use the reflection formula $\Gamma(\nu)\Gamma(1-\nu)$ for the Gamma function to retrieve equation \eqref{fpiexpmellinwithnu}.

We can perform the definition as shown above to reproduce equation \eqref{ex1sec2eq8}.


\subsection{Example 2}\label{fpiexample2}

\hspace{\parindent}We demonstrate in this example the Mellin transform and the contour integration methods to evaluate the finite-part integral $\bbint{0}{\infty} x^{-\lambda} \operatorname{sech}^2(ax) \; {\rm d}x$ for any real $\lambda \geq 1$. Again, we first extract the finite-parts using Mellin transform method. From \cite[p.~30]{TableOfMellinTransform}, we obtain
\begin{equation}
\label{fp1sech2c}
\begin{split}
\mathcal{M}[\operatorname{sech}^2(ax);1-\lambda]
&= \displaystyle\int_{0}^{\infty} x^{-\lambda} \operatorname{sech}^2(ax) \;{\rm d}x \\
&= 2 \lambda \dfrac{a^{\lambda-1}}{\pi^\lambda} (2^{1+\lambda}-1) \dfrac{\sin(\pi\lambda/2)}{\sin(\pi\lambda)} \zeta(1+\lambda) .
\end{split}
\end{equation}
for $\operatorname{Re} a>0$ and $\operatorname{Re} \lambda < 1$. 
Then we have the desired analytic continuation 
\begin{equation}
\label{mellinstarexample2}
\mathcal{M}^*[\operatorname{sech}^2(ax);1-\lambda]
= 2 \lambda \dfrac{a^{\lambda-1}}{\pi^\lambda} (2^{1+\lambda}-1) \dfrac{\sin(\pi\lambda/2)}{\sin(\pi\lambda)} \zeta(1+\lambda) ,
\end{equation}
where \eqref{mellinstarexample2} is well-defined for all non-integer $\lambda>1$, has simple poles for any odd integer $\lambda \geq 1$, and has removable singularities for any even integer $\lambda\geq 1$. It now follows that the finite-part for $\bbint{0}{\infty} x^{-\lambda} \operatorname{sech}^2(ax)\; {\rm d}x$ for all non-integer $\lambda > 1$ is
\begin{equation}
\label{fp1sech2d}
\bbint{0}{\infty} x^{-\lambda} \operatorname{sech}^2(ax) \;{\rm d}x
= 2 \lambda \dfrac{a^{\lambda-1}}{\pi^\lambda} (2^{1+\lambda}-1) \dfrac{\sin(\pi\lambda/2)}{\sin(\pi\lambda)} \zeta(1+\lambda).
\end{equation}
Furthermore, we can let $\lambda=m+\nu$ for $m=1,2,3,\cdots$ and $0<\nu<1$ on  \eqref{fp1sech2d} to arrive at the finite-part integral
\begin{equation}
\label{fp1sech3}
\begin{split}
\bbint{0}{\infty}  \dfrac{\operatorname{sech}^2(ax)}{x^{m+\nu}} \;{\rm d}x
&= 2 (-1)^m a^{m+\nu-1} \dfrac{m+\nu}{\pi^{m+\nu}} (2^{1+m+\nu} -1) \\
&\hspace{0.8cm}
\times \dfrac{\sin(\pi(m+\nu)/2)}{\sin(\pi\nu)}\zeta(m+\nu+1) .
\end{split}
\end{equation}

We can also determine the finite-part of $\bbint{0}{\infty} x^{-\lambda} \operatorname{sech}^2(ax)\; {\rm d}x$ for all integer values of $\lambda\geq 1$ from equation \eqref{mellinstarexample2} using regularized limit or Cauchy limit, depending on the nature of the singularity. First, let us write the right-hand side of \eqref{mellinstarexample2} as
\begin{equation}
\label{fp1sech3b}
\mathcal{M}^*[\operatorname{sech}^2(ax);1-\lambda]
= \dfrac{F(\lambda)}{G(\lambda)}
\end{equation}
where 
\[
F(\lambda) = 2 \lambda \dfrac{a^{\lambda-1}}{\pi^\lambda} (2^{1+\lambda}-1) \sin(\pi\lambda/2) \zeta(1+\lambda) ,\;
G(\lambda)=\sin(\pi\lambda) .
\]
For $k=0,1,2,\cdots$, equation \eqref{fp1sech3b} has simple poles at $\lambda = 2k+1$, so that the finite-part integral $\bbint{0}{\infty} x^{-2k-1} \operatorname{sech}^2(ax)\; {\rm d}x$ is equal to the regularized limit of the analytic continuation at $\lambda=2k+1$,
\begin{equation}
\label{fp1sech3c}
\begin{split}
\bbint{0}{\infty} \dfrac{\operatorname{sech}^2(ax)}{x^{2k+1}}{\rm d}x
&= 
\reglim{\lambda}{2k+1}\mathcal{M}^*[\operatorname{sech}^2(ax);1-\lambda] \\
&=\dfrac{F'(2k+1)}{G'(2k+1)}-\dfrac{F(2k+1)}{2} \dfrac{G''(2k+1)}{(G'(2k+1))^2} .
\end{split}
\end{equation}
Again observe that $G''(2k+1)$ vanishes, so that only the first term of \eqref{fp1sech3c} contributes. Implementing  equation \eqref{fp1sech3c} gives us
\begin{equation}
\label{fp1sech4}
\begin{split}
\bbint{0}{\infty}
& \dfrac{\operatorname{sech}^2(ax)}{x^{2k+1}} {\rm d}x
= 2(-1)^{k+1} \dfrac{a^{2k}}{\pi^{2k+2}}\bigg[ (2k+1) (2^{2k+2}-1)\zeta'(2k+2)  \\
&\hspace{-0.6cm}
\left.+ \zeta(2k+2)\left( 2^{2k+2}\left(1+(2k+1) \ln\left(\dfrac{2a}{\pi} \right) \right)
- \left(1+ (2k+1) \ln\left(\dfrac{a}{\pi}\right) \right)\right)\right] .
\end{split}
\end{equation}
On the other hand, for the 
same values of $k$, equation \eqref{fp1sech3b} has removable singularity at $\lambda = 2k+2$. And we can determine the finite-part for $\bbint{0}{\infty} x^{-2k-2} \operatorname{sech}^2(ax)\; {\rm d}x$ using Cauchy limit
\begin{equation}
\label{fp1sech4ax}
\bbint{0}{\infty} \dfrac{\operatorname{sech}^2(ax)}{x^{2k+2}}{\rm d}x
= 
\lim_{\lambda\rightarrow 2k+2}
\mathcal{M}[\operatorname{sech}^2(ax); 1-\lambda]
=\dfrac{F'(2k+2)}{G'(2k+2)} ,
\end{equation}
which leads to
\begin{equation}
\label{fp1sech5}
\begin{split}
\bbint{0}{\infty}
& \dfrac{\operatorname{sech}^2(ax)}{x^{2k+2}} {\rm d}x
= (-1)^{k+1} \dfrac{4k+4}{(2\pi)^{2k+2}} (2a)^{2k+1} \left( 2^{2k+3}-1\right) \zeta(2k+3) .
\end{split}
\end{equation}

\begin{figure}[t!]
\centering\includegraphics[width=0.9\linewidth]{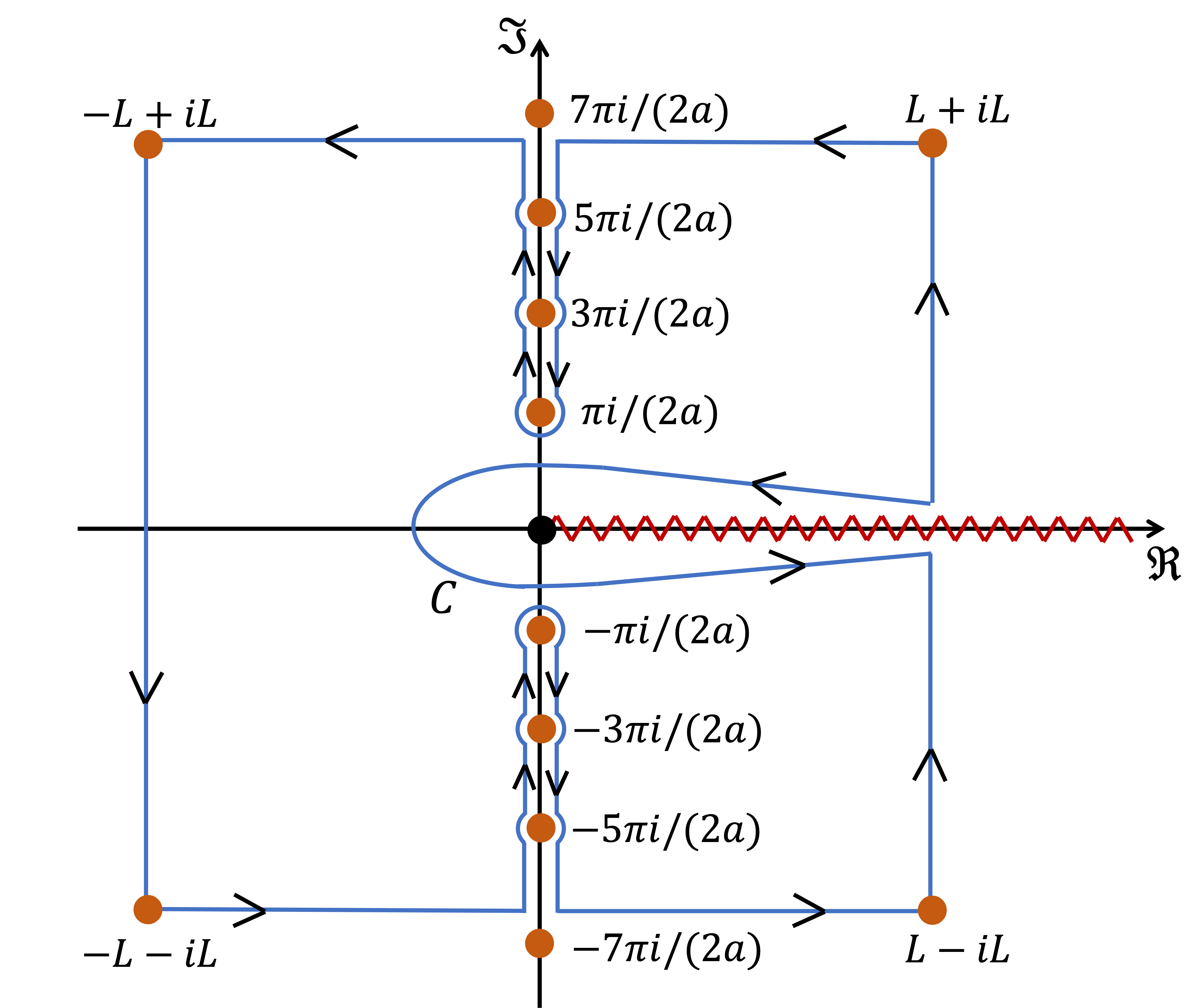}
\caption{The contour of integration that evaluates $\bbint{0}{\infty} x^{-\lambda}\operatorname{sech}^2(ax) \; {\rm d}x$ for any real $\lambda\geq 1$.}
\label{contourmellinsech}
\end{figure}

Another method of obtaining equations \eqref{fp1sech3}, \eqref{fp1sech4}, and \eqref{fp1sech5} is through the contour integral representations of the finite-part integral given by equations \eqref{principalvaluegeneralresult3b} and \eqref{principalvaluegeneralresult3withbc}. To derive the finite-part integral \eqref{fp1sech3}, we use the contour integral representation \eqref{principalvaluegeneralresult3withbc}. We do so by deforming the contour $C$ of integration into the square contour as shown in Figure-\ref{contourmellinsech}. The integrals along the sides of the square vanish in the limit $L\to\infty$, On the other hand, the opposing integrals along the imaginary axis cancel, effectively leaving only the integrals around each pole of $\operatorname{sech}^2(a z)$ to contribute. Then, by the residue theorem, we obtain
\begin{equation}
\label{fp1sech5b}
\bbint{0}{\infty} \dfrac{\operatorname{sech}^2(ax)}{x^{m+\nu}} {\rm d}x
= \dfrac{2\pi i}{1-e^{-2\pi i\nu}} \left[ \sum_{s=0}^\infty \underset{z=z_s}{\text{Res}}  \dfrac{\operatorname{sech}^2(az)}{z^{m+\nu}} 
+ \sum_{t=0}^\infty \underset{z=z_t}{\text{Res}}  \dfrac{\operatorname{sech}^2(az)}{z^{m+\nu}} \right] . 
\end{equation}
where the $z_s$'s and the $z_t$'s are the poles of $\operatorname{sech}^2(az)$ in the upper and lower complex plane, repectively, and are given by 
\[
\begin{split}
z_s &= \left( \dfrac{2s+1}{2a} \right)\pi i 
= \pi \left( \dfrac{2s+1}{2a} \right) e^{i\pi/2} , \\
z_t &= - \left( \dfrac{2t+1}{2a} \right)\pi i 
= \pi \left( \dfrac{2t+1}{2a} \right) e^{3i\pi/2} ,
\end{split}
\]
The sums of residues are evaluated by means of the identity
\begin{equation}
\label{zetaidentity}
\zeta(z)
= \dfrac{1}{1-2^{-z}} \sum_{n=0}^\infty \dfrac{1}{(2n+1)^z}
\end{equation}
for $\operatorname{Re}z>1$ \cite[p.~602]{Handbook1}. Substituting the poles back into equation \eqref{fp1sech5b}, we retrieve equation \eqref{fp1sech3} upon using equation \eqref{zetaidentity}.

We can perform the same steps shown above to derive \eqref{fp1sech4} and \eqref{fp1sech5}. But this time, we utilize equation \eqref{principalvaluegeneralresult3b} to implement the contour integral method in extracting the finite-parts for the remaining two equations.


\subsection{Example 3}\label{fpiexample3}

\hspace{\parindent} We use the Mellin transform method to obtain the finite-part integral 
\begin{equation}
\label{ex2eq2mellin1a}
\bbint{0}{\infty} 
\dfrac{e^{-ax} I_n(ax)}{x^{\lambda+n}} {\rm d}x ,
\end{equation}
where $a>0$, $n$ is any non-negative integer, and $\lambda\geq 1$. In equation \eqref{ex2eq2mellin1a}, we identify $f(x)=x^{-n}e^{-ax}I_n(ax)$ with $f(0)\neq 0$. From \cite[p.~191]{TableOfMellinTransform}, we have the Mellin transform
\begin{equation}
\label{ex2eq2mellin2}
\begin{split}
\mathcal{M}[x^{-n}e^{-ax}I_n(ax); 1-\lambda]
&= \int_{0}^{\infty} 
\dfrac{e^{-ax} I_n(ax)}{x^{\lambda+n}} {\rm d}x \\
&= \sqrt{\pi} \dfrac{(2a)^{\lambda+n-1}}{\sin(\pi\lambda)} \dfrac{\Gamma(\lambda+n-1/2)}{\Gamma(\lambda)\Gamma(\lambda+2n) } ,
\end{split}
\end{equation}
for $\text{Re} \; a>0$ and $-n+1/2<\operatorname{Re}\lambda< 1$. The right-hand side of \eqref{ex2eq2mellin2} is analytic for all non-integer $\lambda$ except at $\lambda=-n+1/2$, and extends the Mellin transform in the entire complex plan. Hence, the desired analytic continuation is given by
\begin{equation}
\label{ex2eq2mellin2bxx}
\mathcal{M}^*[x^{-n} e^{-ax}I_n(ax); 1-\lambda]
=\sqrt{\pi} \dfrac{(2a)^{\lambda+n-1}}{\sin(\pi\lambda)}  \dfrac{\Gamma(\lambda+n-1/2)}{\Gamma(\lambda)\Gamma(\lambda+2n) } ,
\end{equation}
which has simple poles at any integer $\lambda\geq 1$. Then for non-integer $\lambda>1$, we have the finite-part integral, 
\begin{equation}
\label{ex2eq2mellin2b}
\bbint{0}{\infty} 
\dfrac{e^{-ax} I_n(ax)}{x^{\lambda+n}} {\rm d}x
= \sqrt{\pi} \dfrac{(2a)^{\lambda+n-1}}{\sin(\pi\lambda)}  \dfrac{\Gamma(\lambda+n-1/2)}{\Gamma(\lambda)\Gamma(\lambda+2n) } .
\end{equation}

For integer $\lambda=k+1$, $k=0,1,2,\cdots$, the finite-part is obtained via the regularized limit. We rationalize the analytic continuation in the form
\begin{equation}
\label{ex2eq2mellin2bx}
\mathcal{M}^*[x^{-n} e^{-ax}I_n(ax); 1-\lambda]
= \dfrac{F(\lambda)}{G(\lambda)}
\end{equation}
where
\[
\begin{split}
F(\lambda)
= \sqrt{\pi} (2a)^{\lambda+n-1}  \dfrac{\Gamma(\lambda+n-1/2)}{\Gamma(\lambda)\Gamma(\lambda+2n)} , \;
G(\lambda)
=  \sin(\pi\lambda) .
\end{split}
\]
Then 
\begin{equation}
\label{ex2eq2mellin3}
\begin{split}
\bbint{0}{\infty} \dfrac{e^{-ax}I_n(ax)}{x^{k+n+1}} {\rm d}x
&= \reglim{\lambda}{k+1}
\mathcal{M}^*[x^{-n} e^{-ax}I_n(ax); 1-\lambda] \\
&= \dfrac{F'(k+1)}{G'(k+1)}
- F(k+1) \dfrac{G''(k+1)}{2(G'(k+1))^2} .
\end{split}
\end{equation}
Again $G''(k+1)$ vanishes, so that only the first term of \eqref{ex2eq2mellin3} contributes to the finite-part. Executing equation \eqref{ex2eq2mellin3} leads to
\begin{equation}
\label{ex2eq2xxx}
\begin{split}
\bbint{0}{\infty} 
& \dfrac{e^{-ax}I_n(ax)}{x^{k+n+1}} {\rm d}x
= (-1)^{k+1} \dfrac{(2a)^{k+n}}{\sqrt{\pi}} \dfrac{\Gamma(k+n+1/2)}{k!(k+2n)!} \\
&\hspace{0.4cm}
\times \left[\ln(2a) + \psi(k+n+1/2)-\psi(k+1)-\psi(k+2n+1) \right] .
\end{split}
\end{equation}


\begin{appendix}


\section{A Table of Hilbert Transforms}
\label{appendix3}
In this Appendix we list down some Hilbert transforms obtained by direct application of the results obtained in this paper.

\begin{enumerate}[label = A.\arabic*]

\item For $a>0$, $0<\nu<1$, and $-\infty<\omega<\infty$ \label{itemc12}
\[
\begin{split}
\operatorname{PV}\!\!\!\int_{-\infty}^\infty
& \dfrac{e^{iax}}{|x|^\nu(\omega - x)} {\rm d}x 
= \pi \tan(\pi\nu/2)
\; \textup{sgn}(\omega) \dfrac{e^{ia\omega}}{|\omega|^\nu} \\ 
&\hspace{0.5cm} -2\pi i (i\omega)^{-\nu} e^{ia\omega}
\sin(\pi\nu/2) \csc(\pi\nu)
\left[1 - \dfrac{\Gamma(\nu, ia\omega)}{\Gamma(\nu)} \right]
\end{split}
\]

\item For $a>0$, $0<\nu<1$, and $-\infty<\omega<\infty$ \label{itemc13}
\[
\begin{split}
\operatorname{PV}\!\!\!\int_{-\infty}^\infty
& \dfrac{\text{sgn}(x) \; e^{iax}}{|x|^\nu(\omega - x)} {\rm d}x 
= - \pi \cot(\pi\nu/2)
\dfrac{e^{ia\omega}}{|\omega|^\nu} \\
&\hspace{0.5cm} + 2 \pi (iw)^{-\nu} e^{ia\omega} 
\cos(\pi\nu/2) \csc(\pi\nu)
\left[1 - \dfrac{\Gamma(\nu, ia\omega)}{\Gamma(\nu)} \right]
\end{split}
\]

\item \label{itemc10} For $a>0$, $0<\nu<1$, and $\omega>0$
\[
\begin{split}
\operatorname{PV}\!\!\!\int_0^\infty 
\dfrac{\omega e^{-ax}}{x^\nu(\omega^2-x^2)} {\rm d}x
=&\; \dfrac{\pi}{2}\dfrac{1}{\omega^\nu}
\left[
\dfrac{e^{a\omega}}{\sin(\pi\nu)}
- \dfrac{e^{-a\omega}}{\tan(\pi\nu)}
\right] \\
&\hspace{-2.5cm}
-\dfrac{\pi}{\sin(\pi\nu)}
\dfrac{a^{\nu+1} \omega}{\Gamma(2+\nu)}
\pFq{1}{2}{1}{1+\nu/{2}, (3+\nu)/{2}}{\dfrac{a^2\omega^2}{4}} 
\end{split}
\]

\item \label{itemc11} For $a>0$, $0<\nu<1$, and $\omega > 0$
\[
\begin{split}
\operatorname{PV}\!\!\!\int_0^\infty 
\dfrac{x^{1-\nu} e^{-ax}}{\omega^2-x^2} {\rm d}x
=&\; - \dfrac{\pi}{2}\dfrac{1}{\omega^\nu}
\left[
\dfrac{e^{-a\omega}}{\tan(\pi\nu)}
+ \dfrac{e^{a\omega}}{\sin(\pi\nu)}
\right] \\
&\hspace{-2.5cm}
+ \dfrac{\pi}{\sin(\pi\nu)}
\dfrac{a^\nu}{\Gamma(1+\nu)}
\pFq{1}{2}{1}{(1+\nu)/2, 1+\nu/2}{\dfrac{a^2 \omega^2}{4}} 
\end{split}
\]

\item \label{itemc5} For $a>0$ and $\omega>0$
\[
\begin{split}
\operatorname{PV}\!\!\!\int_{0}^\infty 
\dfrac{\omega (J_0(ax))^2}{\omega^2-x^2} {\rm d}x
&= \dfrac{1}{2}\operatorname{PV}\!\!\!\int_{-\infty}^\infty 
\dfrac{(J_0(ax))^2}{\omega-x} {\rm d}x \\
&= 4 \dfrac{a\omega}{\pi} 
\pFq{2}{3}{1,1}{3/2,3/2,3/2}{-(a\omega)^2}
\end{split}
\]

\item \label{itemc6} For $a>0$ and $\omega>0$
\[
\begin{split}
\operatorname{PV}&\!\!\!\int_{0}^\infty 
\dfrac{x(J_0(ax))^2}{\omega^2-x^2} {\rm d}x
= (J_0(a\omega))^2\ln(a\omega) \\
&- \dfrac{1}{2\sqrt{\pi}} 
\sum_{k=0}^\infty 
(-1)^{k}(a\omega)^{2k}
\dfrac{\Gamma(k+1/2)}
{(k!)^3}
\left[3\psi(k+1) - \psi\left(k+\dfrac{1}{2}\right) 
\right] 
\end{split}
\]

\item \label{itemc7} For $a>0$, $0<\nu<1$, and $\omega>0$ 
\[
\begin{split}
\operatorname{PV}&\!\!\!\int_{0}^\infty 
\dfrac{\omega(J_0(ax))^2}
{x^\nu(\omega^2-x^2)} {\rm d}x
= \dfrac{1}{2}
\operatorname{PV}\!\!\!\int_{-\infty}^\infty 
\dfrac{(J_0(ax))^2}
{|x|^\nu(\omega-x)} {\rm d}x
= \dfrac{\pi}{2\omega^\nu} \tan(\pi\nu/2)
\left( J_0(a\omega) \right)^2 \\
&\hspace{-0.3cm} 
+ \dfrac{\sqrt{\pi}}{2} 
\dfrac{a^\nu(a\omega)}{\cos(\pi\nu/2)} 
\dfrac{\Gamma\left(1+\dfrac{\nu}{2}\right)}{\left(\Gamma\left(\dfrac{3+\nu}{2} \right)\right)^{3}} 
\times
\pFq{2}{3}{1, 1+\nu/{2}}{(3+\nu)/{2}, (3+\nu)/{2}, (3+\nu)/{2}}{-(a\omega)^2}
\end{split}
\]

\item \label{itemc8} For $a>0$, $0<\nu<1$, and $\omega>0$ 
\[
\begin{split}
& \operatorname{PV}\!\!\!\int_{0}^\infty 
\dfrac{x^{1-\nu}(J_0(ax))^2}
{\omega^2-x^2} {\rm d}x
= \dfrac{1}{2}
\operatorname{PV}\!\!\!\int_{-\infty}^\infty 
\dfrac{\textup{sgn}(x)(J_0(ax))^2}
{|x|^\nu(\omega-x)} {\rm d}x 
= -\dfrac{\pi}{2\omega^\nu} 
\dfrac{\left( J_0(a\omega) \right)^2}{\tan(\pi\nu/2)} \\
& \hspace{0.6cm} 
+ \dfrac{\sqrt{\pi}}{2} 
\dfrac{a^\nu}{\sin(\pi\nu/2)} 
\dfrac{\Gamma\left(\dfrac{1+\nu}{2}\right)}{\left(\Gamma\left(1+\dfrac{\nu}{2} \right)\right)^{3}}
\times 
\pFq{2}{3}{1,(1+\nu)/2}{1+\nu/2, 1+\nu/2, 1+\nu/2}{-(a\omega)^2}
\end{split}
\]

\item \label{itemc9} For $a>0$, $0<\nu<1$, and $-\infty<\omega<\infty$
\[
\begin{split}
\operatorname{PV}&\!\!\!\int_{-\infty}^\infty 
\dfrac{(J_0(ax))^2}
{x^\nu(\omega-x)} {\rm d}x
= -\pi i \dfrac{(J_0(a\omega))^2}{\omega^\nu} \\
& + \sqrt{\pi}a^{\nu+1} \omega e^{-i\pi\nu/2}
\dfrac{\Gamma\left(1+\dfrac{\nu}{2}\right)}
{\left(\Gamma\left(\dfrac{3+\nu}{2} \right)\right)^{3} }
\times \pFq{2}{3}{1, 1+\nu/{2}}{(3+\nu)/{2}, (3+\nu)/{2}, (3+\nu)/{2}}{-(a\omega)^2} \\
&+ i\sqrt{\pi}a^\nu e^{-i\pi\nu/2}
\dfrac{\Gamma\left(\dfrac{1+\nu}{2}\right)}
{\left(\Gamma\left(1+\dfrac{\nu}{2} \right)\right)^{3} }
\times \pFq{2}{3}{1,(1+\nu)/2}{1+\nu/2, 1+\nu/2, 1+\nu/2}{-(a\omega)^2}
\end{split}
\]

\item \label{itemc1} For $a>0$ and $0<\omega<a$
\[
\begin{split}
\operatorname{PV}\!\!\!\int_{0}^\infty
& \dfrac{x}{\omega^2 - x^2} 
\dfrac{{\rm d}x}{\sqrt{x^2+a^2}}  
= \dfrac{\ln(\omega/a)}
{\sqrt{a^2+\omega^2}} \\
&-\dfrac{1}{2a\sqrt{\pi}} 
\sum_{k=0}^\infty \left(- \dfrac{\omega^2}{a^2}\right)^k \dfrac{\Gamma(k+1/2)}{\Gamma(k+1)} \left[ \psi(k+1) - \psi(k+1/2) \right]
\end{split}
\]

\item \label{itemc2} For $0<\nu<1$, $a>0$, and $\omega>0$
\[
\begin{split}
\operatorname{PV}\!\!\!\int_{0}^\infty
&\; \dfrac{\omega}{x^\nu(\omega^2 - x^2)} 
\dfrac{{\rm d}x}{\sqrt{x^2+a^2}}  
= \dfrac{1}{2} 
\operatorname{PV}\!\!\!\int_{-\infty}^\infty
\dfrac{1}{|x|^\nu (\omega - x)}
\dfrac{{\rm d}x}{\sqrt{x^2+a^2}} \\
=&\; -\dfrac{\omega}{2\sqrt{\pi} a^{2+\nu}}
\Gamma\left(1+\dfrac{\nu}{2}\right) \Gamma\left(-\dfrac{1+\nu}{2}\right)
\pFq{2}{1}{1, 1+\nu/{2}}{(3+\nu)/2}{-\dfrac{\omega^2}{a^2} } \\
&\; + \dfrac{\pi}{2} \dfrac{\tan(\pi\nu/2)}
{\omega^\nu\sqrt{\omega^2+a^2}} 
\end{split}
\]

\item \label{itemc3} For $0<\nu<1$, $a>0$, and $\omega>0$
\[
\begin{split}
\operatorname{PV}\!\!\!\int_{0}^\infty
&\; \dfrac{x^{1-\nu}}{\omega^2 - x^2} 
\dfrac{{\rm d}x}{\sqrt{x^2+a^2}}  
= \dfrac{1}{2} 
\operatorname{PV}\!\!\!\int_{-\infty}^\infty
\dfrac{\textup{sgn}(x)}{|x|^\nu (\omega - x)}
\dfrac{{\rm d}x}{\sqrt{x^2+a^2}} \\
=&\; -\dfrac{1}{2\sqrt{\pi} a^{1+\nu}}
\Gamma\left(-\dfrac{\nu}{2}\right) \Gamma\left(\dfrac{1+\nu}{2}\right) 
\pFq{2}{1}{1,(1+\nu)/2}{1+\nu/2}{-\dfrac{\omega^2}{a^2}} \\
&\; - \dfrac{\pi}{2} \dfrac{\cot(\pi\nu/2)}
{\omega^\nu\sqrt{\omega^2+a^2}} 
\end{split}
\]

\item \label{itemc4} For $0<\nu<1$, $a>0$, and $-\infty<\omega<\infty$
\[
\begin{split}
\operatorname{PV}&\!\!\!\int_{-\infty}^\infty
\dfrac{1}{x^\nu (\omega - x)}
\dfrac{{\rm d}x}{\sqrt{x^2+a^2}} 
= -\dfrac{\pi i}{\omega^\nu}
\dfrac{1}{\sqrt{\omega^2+a^2}} \\
& -\dfrac{e^{-\pi i \nu/2}}
{a^{\nu+1}\sqrt{\pi}}
\left[
i \sin\left(\dfrac{\pi\nu}{2}\right)
\Gamma\left(-\dfrac{\nu}{2}\right)
\Gamma\left(\dfrac{1+\nu}{2}\right)
\pFq{2}{1}{1,(1+\nu)/2}{1+\nu/2}{-\dfrac{\omega^2}{a^2}}
\right. \\
&\left. 
+ \dfrac{\omega}{a}
\cos\left(\dfrac{\pi\nu}{2}\right) 
\Gamma\left(-\dfrac{(\nu+1)}{2}\right)
\Gamma\left(\dfrac{2+\nu}{2}\right)
\pFq{2}{1}{1, 1+\nu/{2}}{(3+\nu)/2}{-\dfrac{\omega^2}{a^2} }
\right] 
\end{split}
\]

\item \label{itemc22} For $c>0$, $0 < \nu < 1$, and $\omega>0$
\[
\begin{split}
\operatorname{PV}\!\!\!\int_0^\infty 
&\; \dfrac{{\rm d}x}{x^\nu(\omega - x)} 
\dfrac{1}{c^3+x^3}
= - \dfrac{\pi}{\omega^\nu} 
\dfrac{\cot(\pi\nu)}{\omega^3+c^3} \\
& + \dfrac{\pi}{3c^{2+\nu}} \dfrac{1}{\omega^3+c^3}
\left[ c^2 \csc\left(\dfrac{\pi\nu}{3}\right)
+ c \omega  
\sec\left( \dfrac{\pi\nu}{3} - \dfrac{\pi}{6}\right)
+ \omega^2
\sec\left( \dfrac{\pi\nu}{3} + \dfrac{\pi}{6}\right)
\right] 
\end{split}
\]

\item \label{itemc23} For $c>0$, $0 < \nu < 1$, and $\omega\neq c > 0$
\[
\begin{split}
\operatorname{PV}\!\!\!\int_0^\infty 
&\; \dfrac{\omega}{x^\nu(\omega^2 - x^2)} 
\dfrac{{\rm d}x}{c^3+x^3}
= \dfrac{\pi}{2\omega^\nu} 
\left[ \dfrac{\csc(\pi\nu)}{c^3 - \omega^3} 
- \dfrac{\cot(\pi\nu)}{c^3 + \omega^3} 
\right]\\
& - \dfrac{\pi}{3c^{\nu+2}} \dfrac{\omega}{c^6 - \omega^6}
\left[ c^2 \omega^2 \csc\left( \dfrac{\pi\nu}{3}\right)
+ \omega^4 
\sec\left( \dfrac{\pi\nu}{3} + \dfrac{\pi}{6}\right)
- c^4 
\sec\left( \dfrac{\pi\nu}{3} - \dfrac{\pi}{6}\right)
\right] 
\end{split}
\]

\item \label{itemc24} For $c>0$, $0 < \nu < 1$, and $\omega\neq c >0$,
\[
\begin{split}
\operatorname{PV}\!\!\!\int_0^\infty 
&\; \dfrac{x^{1-\nu}}{\omega^2 - x^2} 
\dfrac{{\rm d}x}{c^3+x^3}
= - \dfrac{\pi}{2\omega^\nu} 
\left[ \dfrac{\csc(\pi\nu)}{c^3 - \omega^3} 
 + \dfrac{\cot(\pi\nu)}{c^3 + \omega^3} 
\right]\\
& + \dfrac{\pi}{3c^{\nu+1}} \dfrac{1}{c^6 - \omega^6}
\left[ c^4 \csc\left( \dfrac{\pi\nu}{3}\right)
+ c^2 \omega^2 \sec\left( \dfrac{\pi\nu}{3} + \dfrac{\pi}{6}\right)
- \omega^4 
\sec\left( \dfrac{\pi\nu}{3} - \dfrac{\pi}{6}\right)
\right] 
\end{split}
\]

\item \label{itemc14} For $a > 0$, $0 < \nu < 1$, and $0< \omega < c$
\[
\begin{split}
\operatorname{PV}\!\!\!\int_{0}^\infty
\dfrac{\omega}{x^\nu(\omega^2 - x^2)} 
\dfrac{e^{-ax}}{x+c} {\rm d}x 
=&\; -\dfrac{\pi e^{ac}}{\sin(\pi\nu)} 
\dfrac{1}{c^{\nu+1}}
\dfrac{\omega}{c}
\sum_{k=0}^\infty 
\left(\dfrac{\omega}{c}\right)^{2k}
\dfrac{\Gamma(2k+2+\nu,ac)}
{\Gamma(2k+2+\nu)} \\
& + \dfrac{\pi}{2\omega^\nu} 
\left[ 
\dfrac{1}{\sin(\pi\nu)}
\dfrac{e^{a\omega}}{c - \omega}
- \dfrac{1}{\tan(\pi\nu)}
\dfrac{e^{-a\omega}}{\omega+c}
\right] 
\end{split}
\]

\item \label{itemc15} For $a > 0$, $0 < \nu < 1$, and $0< \omega < c$
\[
\begin{split}
\operatorname{PV}\!\!\!\int_{0}^\infty
\dfrac{x^{1-\nu}}{\omega^2 - x^2} 
\dfrac{e^{-ax}}{x+c} {\rm d}x 
=&\; \dfrac{\pi e^{ac}}{\sin(\pi\nu)} 
\dfrac{1}{c^{\nu+1}}
\sum_{k=0}^\infty 
\left(\dfrac{\omega}{c}\right)^{2k}
\dfrac{\Gamma(2k+1+\nu,ac)}
{\Gamma(2k+1+\nu)} \\
& - \dfrac{\pi}{2\omega^\nu} 
\left[ 
\dfrac{1}{\sin(\pi\nu)}
\dfrac{e^{a\omega}}{c - \omega}
+ \dfrac{1}{\tan(\pi\nu)}
\dfrac{e^{-a\omega}}{\omega+c}
\right] 
\end{split}
\]

\item \label{itemc16} For $a > 0$ and $0< \omega < c$
\[
\begin{split}
\operatorname{PV}\!\!\!\int_{0}^\infty
& \dfrac{\omega}{\omega^2 - x^2} 
\dfrac{e^{-ax}}{x+c} {\rm d}x 
= e^{ac} \dfrac{\omega}{c^2} 
\sum_{k=0}^\infty
\left( \dfrac{\omega}{c} \right)^{2k}
\dfrac{\gamma(2k+2,ac)}
{(2k+1)!}
[\ln(a)-\psi(2k+2)] \\
&- a^2 \omega e^{ac}
\sum_{k=0}^\infty 
\dfrac{a^{2k} \omega^{2k}}{(2k+1)! (2k+2)^2} 
\times \pFq{2}{2}{2k+2, 2k+2}{2k+3, 2k+3}{-ac} \\
&+ \dfrac{\omega e^{ac}}{c^2-\omega^2} \ln(c)
+ \dfrac{1}{2}
\left[ \dfrac{e^{-a\omega}}{c+\omega}
- \dfrac{e^{a\omega}}{c-\omega}\right] \ln(\omega)
\end{split}
\]

\item \label{itemc17} For $a > 0$ and $0< \omega < c$
\[
\begin{split}
\operatorname{PV}\!\!\!\int_{0}^\infty
&\dfrac{x}{\omega^2 - x^2} 
\dfrac{e^{-ax}}{x+c} {\rm d}x 
= -\dfrac{e^{ac}}{c} 
\sum_{k=0}^\infty
\left( \dfrac{\omega}{c} \right)^{2k}
\dfrac{\gamma(2k+1,ac)}
{(2k)!}
[\ln(a)-\psi(2k+1)] \\
&+ a e^{ac}
\sum_{k=0}^\infty 
\dfrac{a^{2k} \omega^{2k}}{(2k)! (2k+1)^2} \times 
\pFq{2}{2}{2k+1, 2k+1}{2k+2, 2k+2}{-ac} \\
&- \dfrac{c e^{ac}}{c^2-\omega^2} \ln(c)
+ \dfrac{1}{2}
\left[ \dfrac{e^{-a\omega}}{c+\omega}
+ \dfrac{e^{a\omega}}{c-\omega}\right] \ln(\omega)
\end{split}
\]

\item \label{itemc18} For $s>0$, $\mu > 0$, and $0 < \omega < s$
\[
\begin{split}
\operatorname{PV}&\!\!\!\int_0^\infty 
\dfrac{{\rm d}x}{(\omega - x)} 
\dfrac{1}{(s+x)^\mu}
= - \dfrac{1}{s^\mu\Gamma(\mu)}
\sum_{n=0}^\infty 
\left(-\dfrac{\omega}{s}\right)^n
\dfrac{\Gamma(n+\mu)}{\Gamma(n+1)} \psi(n+1) \\
&+ \dfrac{1}{s^\mu\Gamma(\mu)}
\sum_{n=0}^\infty 
\left(-\dfrac{\omega}{s}\right)^n
\dfrac{\Gamma(n+\mu)}{\Gamma(n+1)} \psi(n+\mu) 
+ \dfrac{\ln(\omega/s)}{(s+\omega)^\mu}
\end{split}
\]

\item \label{itemc19} For $s>0$, $\mu > 0$, $0< \nu <1$, and $\omega>0$
\[
\begin{split}
\operatorname{PV}\!\!\!\int_0^\infty 
& \dfrac{{\rm d}x}{x^\nu(\omega - x)} 
\dfrac{1}{(s+x)^\mu}
= -\dfrac{\pi}{\omega^\nu}
\dfrac{\cot(\pi\nu)}{(s+\omega)^\mu} \\
&+ \dfrac{\pi}{s^{\nu+\mu}}
\dfrac{1}{\sin(\pi\nu)}
\dfrac{\Gamma(\nu+\mu)}
{\Gamma(\mu)\Gamma(1+\nu)} 
\times 
\pFq{2}{1}{1,\mu+\nu}{1+\nu}{-\dfrac{\omega}{s}}
\end{split}
\]

\item \label{itemc20} For $s>0$, $\mu>0$, $0<\nu<1$, and $\omega\neq s > 0$ 
\[
\begin{split}
\operatorname{PV}&\!\!\!\int_0^\infty 
\dfrac{\omega}{x^\nu(\omega^2 - x^2)} 
\dfrac{{\rm d}x}{(s+x)^\mu}
= \dfrac{\pi}{2\omega^\nu} 
\left[
\dfrac{\csc(\pi\nu)}{(s-\omega)^\mu} 
- \dfrac{\cot(\pi\nu)}{(s+\omega)^\mu}
\right] \\
&\hspace{-0.4cm}
- \dfrac{\pi}{s^{\nu+\mu}}
\dfrac{(\omega/s)}{\sin(\pi\nu)}
\dfrac{\Gamma(1+\mu+\nu)}
{\Gamma(\mu)\Gamma(2+\nu)} 
\pFq{3}{2}{1,(1+\mu+\nu)/{2}, 1 + (\mu+\nu)/{2}}{1+\nu/{2}, (3+\nu)/{2}}{\dfrac{\omega^2}{s^2}}
\end{split}
\]

\item \label{itemc21} For $s>0$, $\mu>0$, $0<\nu<1$, and $\omega \neq s > 0$
\[
\begin{split}
\operatorname{PV}&\!\!\!\int_0^\infty 
\dfrac{x}{x^\nu(\omega^2 - x^2)} 
\dfrac{{\rm d}x}{(s+x)^\mu}
= -\dfrac{\pi}{2\omega^\nu} 
\left[
\dfrac{\csc(\pi\nu)}{(s-\omega)^\mu} 
+ \dfrac{\cot(\pi\nu)}{(s+\omega)^\mu}
\right] \\
&+\dfrac{\pi}{s^{\nu+\mu}}
\dfrac{1}{\sin(\pi\nu)}
\dfrac{\Gamma(\mu+\nu)}
{\Gamma(\mu)\Gamma(1+\nu)} 
\pFq{3}{2}{1,(\mu+\nu)/{2}, (1+\mu+\nu)/{2}}{(1+\nu)/{2}, 1+\nu/{2}}{\dfrac{\omega^2}{s^2}}
\end{split}
\]

\item \label{itemc32} For $a>0$, $0 < \nu < 1$, and $0 < \omega < \pi/a$
\[
\begin{split}
\operatorname{PV}\!\!\!\int_0^\infty 
\dfrac{{\rm d}x}{x^{\nu}(\omega-x)}
& \dfrac{1}{e^{ax}+1}
= -\dfrac{\pi}{\omega^\nu} 
\dfrac{\cot(\pi\nu)}{1+e^{a\omega}} \\
& -\pi a^\nu 
\sum_{k=0}^\infty 
(a\omega)^k \dfrac{(2^{k+\nu+1}-1)}{\sin(\pi(k+\nu))}
\dfrac{\zeta(-k-\nu)}{\Gamma(k+\nu+1)} 
\end{split}
\]

\item \label{itemc31} For $a>0$ and $0 < \omega < \pi/a$
\[
\begin{split}
\operatorname{PV}&\!\!\!\int_0^\infty 
\dfrac{{\rm d}x}{(e^{ax}+1)(\omega-x)}
= \dfrac{\gamma}{2} +
\dfrac{1}{2} \ln\left( \dfrac{2a}{\pi} \right) 
+ \sum_{n=1}^\infty (a \omega)^{2n}
\dfrac{(2^{2n+1}-1)}{\Gamma(2n+1)}  \zeta'(-2n) \\
&\; - \dfrac{1}{a \omega} 
\sum_{n=1}^\infty \dfrac{(a\omega)^{2n}}{\Gamma(2n)} 
\left[\ln(a) - 2^{2n}\ln(2a) + (2^{2n}-1)\psi(2n) \right]
\zeta(1-2n) \\
&\;  - \dfrac{1}{a \omega} 
\sum_{n=1}^\infty \dfrac{(a\omega)^{2n}}{\Gamma(2n)} 
(2^{2n}-1)\zeta'(1-2n) 
+ \dfrac{\ln(\omega)}{e^{a\omega}+1} 
\end{split}
\]

\item \label{itemBessel1} For $a > 2$, $\nu\not\in \mathbb{Z}$, and $\omega > 0$
\[
\begin{split}
\operatorname{PV}&\!\!\!\int_{-\infty}^\infty 
\dfrac{\cos(ax)J_{\nu}(x)J_{-\nu}(x)}{\omega - x} {\rm d}x
=2\times \operatorname{PV}\!\!\!\int_{0}^\infty 
\dfrac{\omega\cos(ax)J_{\nu}(x)J_{-\nu}(x)}{\omega^2 - x^2} {\rm d}x
\\
&= \pi \operatorname{sinc}(\pi\nu) \sum_{k=0}^\infty (-1)^k 
\dfrac{(a\omega)^{2k+1}}{(2k+1)!}  
\pFq{3}{2}{1/2,  -k-1/2, -k}{1+\nu,1-\nu}{\dfrac{4}{a^2}}
\end{split}
\]

\item \label{itemBessel1x} For $a > 2$, $\nu\not\in \mathbb{Z}$, and $\omega > 0$
\[
\begin{split}
\operatorname{PV}&\!\!\!\int_{-\infty}^\infty 
\dfrac{\operatorname{sgn}(x)\cos(ax)J_{\nu}(x)J_{-\nu}(x)}{\omega - x} {\rm d}x
= 2\times \operatorname{PV}\!\!\!\int_{0}^\infty 
\dfrac{x\cos(ax)J_{\nu}(x)J_{-\nu}(x)}{\omega^2 - x^2} {\rm d}x
\\
&= 2\cos(a\omega)J_\nu(\omega)J_{-\nu}(\omega)\ln(\omega) \\
&\hspace{0.3cm}
+ 2\operatorname{sinc}(\pi\nu)
\sum_{k=0}^\infty (-1)^k \dfrac{(a\omega)^{2k}}{(2k)!}
\ln(a) \times \pFq{3}{2}{1/2,-k,-k+1/2}{1+\nu,1-\nu}{\dfrac{4}{a^2}}
\\ 
&\hspace{0.3cm}
-\dfrac{2}{\pi} \dfrac{1}{a^2}
\sum_{k=0}^\infty \dfrac{(-1)^k}{k+1} \dfrac{(2\omega)^{2k}}{(2k)!} 
\dfrac{\Gamma(k+1/2)\Gamma(k+3/2)}{\Gamma(2+k+\nu)\Gamma(2+k-\nu)} 
\pFq{4}{3}{1, 1, {3}/{2}, k+{3}/{2}}{2+k, 2+k+\nu, 2+k-\nu}{\dfrac{4}{a^2}} \\
&\hspace{0.3cm}
- 2\operatorname{sinc}(\pi\nu) \sum_{k,l=0}^\infty (-1)^{k+l} \dfrac{(a\omega)^{2k+2l}}{(2k+2l)!}
\dfrac{(1/2)_l(-k-l)_l(-k-l+1/2)_l}{(1+\nu)_l(1-\nu)_l} \psi(2k+1) \dfrac{(4/a^2)^l}{l!}  
\end{split}
\]

\item \label{itemBessel2} For $a > 2$, $\nu\not\in \mathbb{Z}$, and $\omega > 0$
\[
\begin{split}
\operatorname{PV}&\!\!\!\int_{-\infty}^\infty 
\dfrac{\sin(ax)J_{\nu}(x)J_{-\nu}(x)}{\omega - x} {\rm d}x
= 2\times  \operatorname{PV}\!\!\!\int_{0}^\infty 
\dfrac{x\sin(ax)J_{\nu}(x)J_{-\nu}(x)}{\omega^2 - x^2} {\rm d}x
\\
&= - \pi \operatorname{sinc}(\pi\nu) \sum_{k=0}^\infty (-1)^k 
\dfrac{(a\omega)^{2k}}{(2k)!}  
\pFq{3}{2}{1/2,  -k+1/2, -k}{1+\nu,1-\nu}{\dfrac{4}{a^2}}
\end{split}
\]

\item \label{itemBessel2x} For $a > 2$, $\nu\not\in \mathbb{Z}$, and $\omega > 0$
\[
\begin{split}
\operatorname{PV}&\!\!\!\int_{-\infty}^\infty 
\dfrac{\operatorname{sgn}(x)\sin(ax)J_{\nu}(x)J_{-\nu}(x)}{\omega - x} {\rm d}x
= 2\times\operatorname{PV}\!\!\!\int_{0}^\infty 
\dfrac{\omega\sin(ax)J_{\nu}(x)J_{-\nu}(x)}{\omega^2 - x^2} {\rm d}x
\\
&= 2\sin(a\omega)J_\nu(\omega)J_{-\nu}(\omega)\ln(\omega) \\
&\hspace{0.3cm}
+ 2\operatorname{sinc}(\pi\nu)
\sum_{k=0}^\infty (-1)^k \dfrac{(a\omega)^{2k+1}}{(2k+1)!}
\ln(a) \times 
\pFq{3}{2}{1/{2},-k,-k-1/{2}}{1+\nu,1-\nu}{\dfrac{4}{a^2}}
\\ 
&\hspace{0.3cm}
+\dfrac{2}{\pi} \dfrac{1}{a}
\sum_{k=0}^\infty \dfrac{(-1)^k}{k+1} \dfrac{(2\omega)^{2k+1}}{(2k+1)!} 
\dfrac{(\Gamma(k+3/2))^2}{\Gamma(2+k+\nu)\Gamma(2+k-\nu)} \pFq{4}{3}{1/{2}, 1, 1, k+3/{2}}{2+k, 2+k+\nu, 2+k-\nu}{\dfrac{4}{a^2}} \\
&\hspace{0.3cm}
-2\operatorname{sinc}(\pi\nu) \sum_{k,l=0}^\infty (-1)^{k+l} \dfrac{(a\omega)^{2k+2l+1}}{(2k+2l+1)!}
\dfrac{(1/2)_l(-k-l)_l(-k-l-1/2)_l}{(1+\nu)_l(1-\nu)_l} \psi(2k+2) \dfrac{(4/a^2)^l}{l!} 
\end{split}
\]

\item \label{itemBessel3} For $a > 2$, $n\in 0, \mathbb{Z}^+$, and $\omega > 0$
\[
\begin{split}
\operatorname{PV}&\!\!\!\int_{-\infty}^\infty 
\dfrac{\cos(ax) (J_{n}(x))^2}{\omega - x} {\rm d}x
= 2 \times \operatorname{PV}\!\!\!\int_0^\infty 
\dfrac{\omega \cos(ax) (J_{n}(x))^2}{\omega^2 - x^2} {\rm d}x
\\
&= \pi \dfrac{(\omega/2)^{2n}}{(n!)^2} \sum_{k=0}^\infty (-1)^k 
\dfrac{(a\omega)^{2k+1}}{(2k+1)!}  
\pFq{3}{2}{n+1/2,  -k-1/2, -k}{1+n,1+2n}{\dfrac{4}{a^2}}
\end{split}
\]

\item \label{itemBessel3x} For $a > 2$, $n\in 0, \mathbb{Z}^+$, and $\omega > 0$
\[
\begin{split}
\operatorname{PV}&\!\!\!\int_{-\infty}^\infty 
\dfrac{\operatorname{sgn}(x) \cos(ax) (J_{n}(x))^2}{\omega - x} {\rm d}x
= 2\times \operatorname{PV}\!\!\!\int_{0}^\infty 
\dfrac{x \cos(ax) (J_{n}(x))^2}{\omega^2 - x^2} {\rm d}x
\\
&\hspace{-0.3cm}
= 2 \cos(a\omega)(J_n(\omega))^2\ln(\omega) \\
&\hspace{-0.1cm}
- 2 \left( \dfrac{1}{2a}\right)^{2n} \dfrac{(-1)^n}{(n!)^2} \sum_{k=0}^{n-1} (-1)^k (a\omega)^{2k} (2n-2k-1)!
\times \pFq{3}{2}{n+1/2,n-k,n-k+1/2}{n+1,2n+1}{\dfrac{4}{a^2}} \\
&\hspace{-0.1cm}
+ 2 \dfrac{(\omega/2)^{2n}}{(n!)^2} 
\sum_{k=0}^\infty (-1)^k \dfrac{(a\omega)^{2k}}{(2k)!} 
\ln(a) \times 
\pFq{3}{2}{n+1/{2}, -k, -k+1/{2}}{n+1,2n+1}{\dfrac{4}{a^2}} \\
&\hspace{-0.1cm} 
- \dfrac{2a^{-2}\omega^{2n}}{\sqrt{\pi}} \sum_{k=0}^\infty 
\dfrac{(-1)^k \omega^{2k}  \times \Gamma(k+n+3/2)}{(k+1)!(k+n+1)!(k+2n+1)!} 
\times 
\pFq{4}{3}{1,1,3/{2},k+n+3/{2}}{2+k, 2+k+n,2+k+2n}{\dfrac{4}{a^2}} \\
&\hspace{-0.1cm} 
- 2 \dfrac{(\omega/2)^{2n}}{(n!)^2} 
\sum_{k,l=0}^\infty (-1)^{k+l} \dfrac{(a\omega)^{2k+2l}}{(2k+2l)!} \dfrac{(n+1/2)_l (-k-l)_l (-k-l+1/2)_l }{(n+1)_l (2n+1)_l} \psi(2k+1) \dfrac{(4/a^2)^l}{l!} 
\end{split}
\]

\item \label{itemBessel4} For $a > 2$, $n\in 0, \mathbb{Z}^+$, and $\omega > 0$
\[
\begin{split}
\operatorname{PV}&\!\!\!\int_{-\infty}^\infty 
\dfrac{\sin(ax) (J_{n}(x))^2}{\omega - x} {\rm d}x
= 2\times \operatorname{PV}\!\!\!\int_{0}^\infty 
\dfrac{x \sin(ax) (J_{n}(x))^2}{\omega^2 - x^2} {\rm d}x
\\
&= - \pi \dfrac{(\omega/2)^{2n}}{(n!)^2} \sum_{k=0}^\infty (-1)^k 
\dfrac{(a\omega)^{2k}}{(2k)!}  
\pFq{3}{2}{n+1/2,  -k+1/2, -k}{1+n,1+2n}{\dfrac{4}{a^2}}
\end{split}
\]

\item \label{itemBessel4x} For $a > 2$, $n\in 0, \mathbb{Z}^+$, and $\omega > 0$
\[
\begin{split}
\operatorname{PV}&\!\!\!\int_{-\infty}^\infty 
\dfrac{\operatorname{sgn}(x) \sin(ax) (J_{n}(x))^2}{\omega - x} {\rm d}x
= 2\times \operatorname{PV}\!\!\!\int_{0}^\infty 
\dfrac{\omega \sin(ax) (J_{n}(x))^2}{\omega^2 - x^2} {\rm d}x
\\
&\hspace{-0.2cm}= 2 \sin(a\omega)(J_n(\omega))^2\ln(\omega) \\
&\hspace{0.1cm}+ 2 \left( \dfrac{1}{2a}\right)^{2n} \dfrac{(-1)^n}{(n!)^2} \sum_{k=0}^{n-1} (-1)^k (a\omega)^{2k+1} (2n-2k-2)!
\times \pFq{3}{2}{n+1/2,n-k,n-k-1/2}{n+1,2n+1}{\dfrac{4}{a^2}} \\
&\hspace{0.1cm}+ 2 \dfrac{(\omega/2)^{2n}}{(n!)^2} 
\sum_{k=0}^\infty (-1)^k \dfrac{(a\omega)^{2k+1}}{(2k+1)!} 
\ln(a) \times 
\pFq{3}{2}{n+1/{2}, -k, -k-1/{2}}{n+1,2n+1}{\dfrac{4}{a^2}} \\
&\hspace{0.1cm} + \dfrac{2a^{-1}\omega^{2n}}{\sqrt{\pi}} \sum_{k=0}^\infty 
\dfrac{(-1)^k \omega^{2k+1}  \times
\Gamma(k+n+3/2)}{(k+1)!(k+n+1)!(k+2n+1)!} \times
\pFq{4}{3}{1/{2},1,1,k+n+3/{2}}{2+k, 2+k+n,2+k+2n}{\dfrac{4}{a^2}} \\
&\hspace{0.1cm} 
- 2 \dfrac{(\omega/2)^{2n}}{(n!)^2} \sum_{k,l=0}^{\infty} (-1)^{k+l} \dfrac{(a\omega)^{2k+2l+1}}{(2k+2l+1)!} \dfrac{(n+1/2)_l (-k-l)_l (-k-l-1/2)_l }{(n+1)_l (2n+1)_l} \psi(2k+2) \dfrac{(4/a^2)^l}{l!} \\
\end{split}
\]

\item \label{itemc25} For $a > 0$ and $-\infty < \omega < \infty$
\[
\begin{split}
\operatorname{PV}&\!\!\!\int_{-\infty}^\infty 
\dfrac{\textup{Ai}(ax)}{\omega - x} {\rm d}x
= - \dfrac{a^2\omega^2}{2}
\pFq{1}{2}{1}{4/3, 5/3}{\dfrac{a^3\omega^3}{9} } \\
&+ \dfrac{\pi}{3\sqrt{3}} 
\dfrac{(a\omega\sqrt{a\omega})^{2/3}
I_{-1/3}\left(\dfrac{2}{3} a^{3/2} \omega^{3/2}\right)
+ (a\omega) I_{1/3}\left(\dfrac{2}{3} a^{3/2} \omega^{3/2} \right)}
{(a\omega\sqrt{a\omega})^{1/3}} 
\end{split}
\]

\item \label{itemc26} For $a > 0$ and $\omega > 0$
\[
\begin{split}
& \operatorname{PV}\!\!\!\int_0^\infty 
\dfrac{\textup{Ai}(ax)}{\omega - x} {\rm d}x
= \textup{Ai}(a\omega)\ln(\omega)
- \dfrac{a^2\omega^2}{6}
\pFq{1}{2}{1}{4/3, 5/3}{\dfrac{a^3\omega^3}{9} } \\
&\hspace{0.5cm}
+ \dfrac{(a\omega)^{1/2}}{9}
\dfrac{\pi}{\sqrt{3}} 
\left[ I_{-1/3}\left(\dfrac{2}{3} a^{3/2} \omega^{3/2}\right)
+ I_{1/3}\left(\dfrac{2}{3} a^{3/2} \omega^{3/2} \right) 
\right]
\\
&\hspace{0.5cm}
+ \dfrac{(a\omega)^{1/2}}{9}
\ln\left( \dfrac{a^3}{9}\right)
\left[ I_{-1/3}\left(\dfrac{2}{3} a^{3/2} \omega^{3/2}\right)
- I_{1/3}\left(\dfrac{2}{3} a^{3/2} \omega^{3/2} \right) 
\right] \\
&\hspace{0.5cm}
+ \sum_{n=1}^\infty 
\left( \dfrac{a^3 \omega^3}{9} \right)^n
\left[ \dfrac{3^{-1/3}}{(a\omega)^2} 
\dfrac{\psi(n)+\psi(n+1/3)}{\Gamma(n)\Gamma(n+1/3)} - \dfrac{3^{1/3}}{(a\omega)^3} \dfrac{\psi(n)+\psi(n-1/3)}{\Gamma(n)\Gamma(n-1/3)}
\right]
\end{split}
\]

\item \label{itemc27} For $a > 0$ and $\omega > 0$
\[
\begin{split}
\operatorname{PV}\!\!\!\int_0^\infty 
& \dfrac{\textup{Ai}(ax)}{x^\nu(\omega - x)} {\rm d}x
= -\dfrac{\pi}{\tan(\pi\nu)}
\dfrac{\textup{Ai}(a\omega)}{\omega^\nu} \\
& \hspace{-0.7cm} 
- \dfrac{2}{3^{5/2+2\nu/3}}
\dfrac{ \pi a^{\nu+2} \omega^2}
{1+2\cos\left( \dfrac{2\pi\nu}{3} \right)} 
\dfrac{
\pFq{1}{2}{1}{(4+\nu)/{3}, (5+\nu)/{3}}{\dfrac{a^3\omega^3}{9} }
}{\Gamma\left(\dfrac{4+\nu}{3}\right)\Gamma\left(\dfrac{5+\nu}{3}\right)} \\
& \hspace{-0.7cm} 
- \dfrac{2}{3^{11/6+2\nu/3}}
\dfrac{\pi a^{\nu+1}\omega}
{1+\sqrt{3}\sin\left( \dfrac{2\pi\nu}{3}\right) - \cos\left( \dfrac{2\pi\nu}{3} \right)} 
\dfrac{
\pFq{1}{2}{1}{1+\nu/{3}, (4+\nu)/{3}}{\dfrac{a^3\omega^3}{9} }
}{\Gamma\left(\dfrac{3+\nu}{3}\right)
\Gamma\left(\dfrac{4+\nu}{3}\right)}
\\
&\hspace{-0.7cm} 
- \dfrac{2}{3^{7/6+2\nu/3}} 
\dfrac{\pi a^\nu}
{1-\sqrt{3}\sin\left( \dfrac{2\pi\nu}{3}\right) - \cos\left( \dfrac{2\pi\nu}{3} \right)}
\dfrac{
\pFq{1}{2}{1}{(2+\nu)/3, 1+\nu/3}{\dfrac{a^3\omega^3}{9}}
}{\Gamma\left(\dfrac{2+\nu}{3}\right)
\Gamma\left(\dfrac{3+\nu}{3}\right)}
\end{split}
\]

\item \label{itemc28} For $a > 0$ and $-\infty<\omega<\infty$
\[
\begin{split}
\operatorname{PV}\!\!\!\int_{-\infty}^\infty 
\dfrac{\textup{Ai}(ax)}{x^\nu(\omega - x)} {\rm d}x
= - & \pi i
\dfrac{\textup{Ai}(a\omega)}{\omega^\nu} \\
&\hspace{-3.2cm}
- \dfrac{2\pi a^{\nu+2} \omega^2}{3^{5/2+2\nu/3} } 
\dfrac{1 + 2e^{-\pi i \nu}\cos\left(\dfrac{\pi\nu}{3}\right)}
{1+2\cos\left(\dfrac{2\pi\nu}{3}\right)} \dfrac{
\pFq{1}{2}{1}{(4+\nu)/3, (5+\nu)/3}{\dfrac{a^3\omega^3}{9}}
}
{\Gamma\left(\dfrac{4+\nu}{3}\right)
\Gamma\left(\dfrac{5+\nu}{3}\right)} \\
&\hspace{-3.2cm}
- \dfrac{2\pi a^{\nu+1} \omega}{3^{11/6+2\nu/3} }
\dfrac{1
- e^{-\pi i \nu}
\left[\sqrt{3} \sin\left(\dfrac{\pi\nu}{3}\right) 
+ \cos\left(\dfrac{\pi\nu}{3}\right) \right]}
{1 + \sqrt{3}\sin\left( \dfrac{2\pi\nu}{3}\right) -\cos\left( \dfrac{2\pi\nu}{3}\right)} 
\dfrac{
\pFq{1}{2}{1}{1+\nu/3,(4+\nu)/3}{\dfrac{a^3\omega^3}{9}}
}{ \Gamma\left(\dfrac{3+\nu}{3}\right)
\Gamma\left(\dfrac{4+\nu}{3}\right)} 
\\
&\hspace{-3.2cm}
- \dfrac{2 \pi a^\nu}{3^{7/6+2\nu/3} }  
\dfrac{1
+ e^{-\pi i \nu}
\left[\sqrt{3} \sin\left(\dfrac{\pi\nu}{3}\right) 
- \cos\left(\dfrac{\pi\nu}{3}\right) \right]}
{1 - \sqrt{3}\sin\left( \dfrac{2\pi\nu}{3}\right) -\cos\left( \dfrac{2\pi\nu}{3}\right)} 
\dfrac{
\pFq{1}{2}{1}{(2+\nu)/3,1+\nu/3}{\dfrac{a^3\omega^3}{9}}
}{\Gamma\left(\dfrac{2+\nu}{3}\right)
\Gamma\left(\dfrac{3+\nu}{3}\right)}
\end{split}
\]

\item \label{itemc29} For $a>0$ and $\omega>0$
\[
\begin{split}
\operatorname{PV}\!\!\!\int_0^\infty 
& \dfrac{\text{Ai}(ax) \text{Ai}'(ax)}
{\omega-x} {\rm d}x
= \text{Ai}(a\omega) \text{Ai}'(a\omega) \ln(\omega)
- \dfrac{\ln(12^{1/3} a)}{\pi \sqrt{12}}
\pFq{1}{2}{1/2}{1/3, 2/3}{\dfrac{4}{9} a^3\omega^3} \\
&\hspace{-0.8cm}
+ \dfrac{\ln(12^{1/3} a)}{\sqrt{12\pi}}  
\left[ 
\dfrac{12^{1/3} a\omega}{\Gamma(1/6)}
\pFq{1}{2}{5/6}{2/3, 4/3}{\dfrac{4}{9}a^3\omega^3}
- \dfrac{18^{1/3} a^2\omega^2}{\Gamma(-1/6)}
\pFq{1}{2}{7/6}{4/3, 5/3}{\dfrac{4}{9} a^3\omega^3}
\right] \\
&\hspace{-0.8cm} 
- \dfrac{1}{3\sqrt{12\pi}} 
\sum_{k=0}^\infty 
\dfrac{(-12^{1/3}a\omega)^k \psi(1/2-k/3)}
{\Gamma(k+1)\Gamma(1/2-k/3)} + \dfrac{1}{\sqrt{12\pi}} 
\sum_{k=0}^\infty 
\dfrac{(-12^{1/3}a\omega)^k \psi(k+1)}
{\Gamma(k+1)\Gamma(1/2-k/3)} 
\end{split}
\]

\item \label{itemc30} For $a>0$, $\omega>0$, and $0<\nu<1$ 
\[
\begin{split}
\operatorname{PV}&\!\!\!\int_0^\infty 
\dfrac{\text{Ai}(ax) \text{Ai}'(ax)}
{x^\nu(\omega-x)} {\rm d}x
= - \dfrac{\pi}{\tan(\pi\nu)}
\dfrac{\text{Ai}(a\omega) \text{Ai}'(a\omega)}{\omega^\nu} \\
&\hspace{-0.4cm}
-\dfrac{12^{\nu/3}a^\nu}{\sqrt{12\pi}} 
\dfrac{\cos(\pi\nu/3)}{\sin(\pi\nu)}
\dfrac{\Gamma\left(\dfrac{1}{2}
+\dfrac{\nu}{3}\right)}{\Gamma(1+\nu)}
\pFq{2}{3}{1,(3+2\nu)/6}{(1+\nu)/3, (2+\nu)/3, 1+\nu/3}{\dfrac{4}{9} a^3\omega^3} \\
&\hspace{-0.4cm}
+12^{1/3}a\omega\dfrac{12^{\nu/3}a^\nu}{\sqrt{12\pi}} 
\dfrac{\cos\left(\dfrac{\pi}{3}(1+\nu)\right)}{\sin(\pi\nu)}
\dfrac{\Gamma\left(\dfrac{5}{6}
+\dfrac{\nu}{3}\right)}{\Gamma(2+\nu)}
\pFq{2}{3}{1, (5+2\nu)/6}{(2+\nu)/3, 1+\nu/3, (4+\nu)/3}{\dfrac{4}{9} a^3\omega^3} \\
&\hspace{-0.4cm}
+ \left(12^{1/3}a\omega\right)^2
\dfrac{12^{\nu/3}a^\nu}{\sqrt{12\pi}} 
\dfrac{\cos\left(\dfrac{\pi}{3}(1-\nu)\right)}{\sin(\pi\nu)}
\dfrac{\Gamma\left(\dfrac{7}{6}
+\dfrac{\nu}{3}\right)}{\Gamma(3+\nu)}
\pFq{2}{3}{1, (7+2\nu)/6}{1+\nu/3, (4+\nu)/3, (5+\nu)/3}{\dfrac{4}{9} a^3\omega^3}
\end{split}
\]

\end{enumerate}


\newpage
\section{Table of Finite-part Integrals}
\label{appendix4}
In this Appendix, we list down the relevant finite-part integrals used in obtaining the Hilbert transform in the preceding Appendix. They were extracted by means of the analytic continuation of the Mellin transform.

\begin{enumerate}[label = B.\arabic*]

\item Used in items (\ref{itemc10}) and (\ref{itemc11})
\[
\begin{split}
\bbint{0}{\infty} \dfrac{e^{-ax}}{x^{m+\nu}} {\rm d}x
= \dfrac{(-1)^m a^{m+\nu-1}\pi}{\sin(\pi\nu)\Gamma(m+\nu)} \;\;\;\;\;\;\text{for}\;\; a > 0, \; 0 < \nu < 1,\; m=1, 2, 3, \cdots
\end{split}
\]

\item Used in items (\ref{itemc5}), (\ref{itemc7}), (\ref{itemc8}), and (\ref{itemc9})
\[
\begin{split}
\bbint{0}{\infty} \dfrac{(J_0(ax))^2}{x^\lambda} {\rm d}x
= \dfrac{\sqrt{\pi}a^{\lambda-1}\Gamma(\lambda/2)}{2\cos(\pi\lambda/2)(\Gamma((\lambda+1)/2))^3}
\;\;\;\;\;\text{for}\;\; a>0, \; \lambda > 1, \; \lambda\neq 1, 3, 5, \cdots
\end{split}
\]

\item Used in item (\ref{itemc6})
\[
\begin{split}
\bbint{0}{\infty} \dfrac{(J_0(ax))^2}{x^{2n+1}} {\rm d}x
=&\; (-1)^n \dfrac{a^{2n}(1/2)_n}{(n!)^3} \left(\dfrac{3}{2}\psi(n+1) - \dfrac{1}{2}\psi(n+1/2) - \ln a \right) \\ &\;\;\;\;\text{for}\;\; a > 0,\; n = 0, 1, 2, \cdots 
\end{split}
\]

\item Used in item (\ref{itemc1})
\[
\begin{split}
\bbint{0}{\infty} \dfrac{{\rm d}x}{x^{2k+1} \sqrt{x^2 + a^2}}
= & \dfrac{(-1)^k}{2 \sqrt{\pi} a^{2k+1}}
\dfrac{\Gamma(k+1/2)}{\Gamma(k+1)}
[2 \ln a + \psi(k+1) - \psi(k+1/2)]
\\
& \;\;\;\text{for}\;\; a > 0, \; k =0, 1, 2, \cdots 
\end{split}
\]

\item Used in items (\ref{itemc2})-(\ref{itemc4})
\[
\begin{split}
\bbint{0}{\infty} \dfrac{{\rm d}x}{x^{m+\nu} \sqrt{x^2 + a^2}}
= & \dfrac{a^{-(m+\nu)}}{2\sqrt{\pi}} \Gamma\left(\dfrac{1}{2} - \dfrac{m+\nu}{2}\right)
\Gamma\left(\dfrac{m+\nu}{2}\right) \\
&\;\;\;\text{for}\;\; a > 0, \; m=1, 2, 3, \cdots
\; 0 < \nu < 1
\end{split}
\]

\item Used in items (\ref{itemc22})-(\ref{itemc24})
\[
\begin{split}
\bbint{0}{\infty} \dfrac{{\rm d}x}{x^{m+\nu} (x^3+c^3)}
=&\; \dfrac{\pi}{3c^{m+\nu+2}\cos(\pi(2m+2\nu+1)/6)} \\
& \;\text{for}\;\; c > 0,
\; 0 < \nu < 1,\;
m=1, 2, 3, \cdots
\end{split}
\]

\item Used in items (\ref{itemc14}) and (\ref{itemc15})
\[
\begin{split}
\bbint{0}{\infty} \dfrac{e^{-ax}}{x^{m+\nu}(x+c)} {\rm d}x
= \dfrac{(-1)^m\pi e^{ac}\Gamma(m+\nu,ac)}{c^{m+\nu}\sin(\pi\nu)\Gamma(m+\nu)} 
\;\;\;\;\;\text{for}\;\; a,c > 0, \; 0 < \nu < 1,\; m=1, 2, 3, \cdots
\end{split}
\]

\item Used in items (\ref{itemc16}) and (\ref{itemc17})
\[
\begin{split}
\bbint{0}{\infty} \dfrac{e^{-ax}}{x^{n+1}(x+c)} {\rm d}x
= (-1)^n e^{ac}
& \left[ 
\dfrac{\ln c}{c^{n+1}} 
+ \dfrac{1}{n! c^{n+1}}(\ln a - \psi(n+1)) \gamma(1+n, ac) \right. \\
&\hspace{0.4cm}
\left. - \dfrac{a^{n+1}}{n!(n+1)^2} \times
\pFq{2}{2}{n+1, n+1}{n+2, n+2}{-ac}
\right] \\
& \;\text{for}\;\; a > 0, \; c > 0,\; n = 0, 1, 2, \cdots 
\end{split}
\]

\item Used in item (\ref{itemc18})
\[
\begin{split}
\bbint{0}{\infty} \dfrac{{\rm d}x}{x^{n+1} (s+x)^\mu}
=&\; \dfrac{(-1)^n\Gamma(n+\mu)}{s^{n+\mu} n! \Gamma(\mu)}
\left[\ln(s)+\psi(n+1)-\psi(n+\mu)\right] \\
& \;\text{for}\;\; s > 0,
\; \mu > 0,\;
n = 0, 1, 2, \cdots
\end{split}
\]

\item Used in items (\ref{itemc19})-(\ref{itemc21})
\[
\begin{split}
\bbint{0}{\infty} \dfrac{{\rm d}x}{x^{m+\nu} (s+x)^\mu}
=&\; \dfrac{(-1)^m \pi \Gamma(m + \nu + \mu - 1)}{s^{m + \mu + \nu - 1} \sin(\pi\nu) \Gamma(\mu) \Gamma(m+\nu)} \\
& \;\text{for}\;\; s > 0,
\; 0 < \nu < 1, \; \mu > 0,\;
m=1, 2, 3, \cdots
\end{split}
\]

\item Used in item (\ref{itemc32})
\[
\begin{split}
\bbint{0}{\infty} \dfrac{{\rm d}x}{x^{m+\nu}} \dfrac{1}{e^{ax}+1}
=&\; -(-1)^m \pi 
\dfrac{(2^{m+\nu}-1) a^{m+\nu-1}}
{\sin(\pi\nu)}
\dfrac{\zeta(1-m-\nu)}
{\Gamma(m+\nu)} \\
& \;\text{for}\;\; a > 0,
\; 0 < \nu < 1,\;
m=1, 2, 3, \cdots
\end{split}
\]

\item Used in item (\ref{itemc31})
\[
\begin{split}
\bbint{0}{\infty} \dfrac{{\rm d}x}{x} \dfrac{1}{e^{ax}+1}
=&\; \ln\left( \sqrt{\dfrac{\pi}{2a}}\right) 
- \dfrac{\gamma}{2}
\;\;\;\;\;\;\text{for}\;\; a > 0
\end{split}
\]

\item Used in item (\ref{itemc31})
\[
\begin{split}
\bbint{0}{\infty} \dfrac{{\rm d}x}{x^{2n}} \dfrac{1}{e^{ax}+1}
=&\; \dfrac{a^{2n-1}B_{2n}}{(2n)!} \left(2^{2n} \ln 2 - (2^{2n}-1)
(\psi(2n) - \ln a ) \right) \\
&+ \dfrac{a^{2n-1}(2^{2n}-1)}{(2n-1)!}\zeta'(1-2n)
\;\;\;\;\;\text{for}\;\; a > 0,
\; n=1, 2, 3, \cdots
\end{split}
\]

\item Used in item (\ref{itemc31})
\[
\begin{split}
\bbint{0}{\infty} \dfrac{{\rm d}x}{x^{2n+1}} \dfrac{1}{e^{ax}+1}
=&\; -\dfrac{(2^{2n+1}-1)a^{2n}}{(2n)!} \zeta'(-2n) 
\;\;\;\;\text{for}\;\; a > 0,
\; n=1, 2, 3, \cdots
\end{split}
\]

\item Used in item (\ref{itemBessel1})
\[
\begin{split}
\setlength\arraycolsep{1pt}
\bbint{0}{\infty} 
&\dfrac{\cos(ax)J_{\nu}(bx)J_{-\nu}(bx)}{x^{2k+2}} {\rm d}x 
= -\dfrac{\pi}{2} \operatorname{sinc}(\pi\nu) \dfrac{(-1)^{k} a^{2k+1}}{(2k+1)!} \\
&\hspace{-0.4cm}
\times \pFq{3}{2}{1/2,  -k-1/2, -k}{1+\nu,1-\nu}{\dfrac{4b^2}{a^2}}
\;\;\;\;\text{for}\;\; a>2b,\; \nu \neq 0, 1, 2, \cdots,
\; k=0, 1, 2, \cdots
\end{split}
\]

\item Used in item (\ref{itemBessel1x})
\[
\begin{split}
&\bbint{0}{\infty} \dfrac{\cos(ax)J_\nu(bx)J_{-\nu}(bx)}{x^{2k+1}} {\rm d}x = -\operatorname{sinc}(\pi\nu) \dfrac{(-1)^k a^{2k}}{(2k)!} \\
&\hspace{0.4cm} 
\times \left[ \ln(a) \times \pFq{3}{2}{1/2,-k,-k+1/2}{1+\nu,1-\nu}{\dfrac{4b^2}{a^2}}
\right. \\ &\hspace{1.2cm}
\left. - \sum_{l=0}^k \dfrac{(1/2)_l(-k)_l(-k+1/2)_l}{(1+\nu)_l(1-\nu)_l} \psi(2k-2l+1) \dfrac{(4b^2/a^2)^l}{l!} \right] \\
&\hspace{1.0cm} 
+ (-1)^k \dfrac{2^{2k}}{\pi} \dfrac{b^{2k+2}}{k+1} \dfrac{a^{-2}}{(2k)!} \dfrac{\Gamma(k+1/2)\Gamma(k+3/2)}{\Gamma(2+k+\nu)\Gamma(2+k-\nu)} 
\pFq{4}{3}{1, 1, {3}/{2}, k+{3}/{2}}{2+k, 2+k+\nu, 2+k-\nu}{\dfrac{4b^2}{a^2}} \\
&\hspace{5.3cm} 
\;\;\;\;\text{for}\;\; a>2b,\; \nu \neq 0, 1, 2, \cdots,
\; k=0, 1, 2, \cdots
\end{split}
\]

\item Used in item (\ref{itemBessel2})
\[
\begin{split}
\setlength\arraycolsep{1pt}
\bbint{0}{\infty} &\dfrac{\sin(ax)J_{\nu}(bx)J_{-\nu}(bx)}{x^{2k+3}} {\rm d}x 
= - \dfrac{\pi}{2} \operatorname{sinc}(\pi\nu)\dfrac{(-1)^{k}a^{2k+2}}{(2k+2)!} \\
&\hspace{-0.4cm}
\times \pFq{3}{2}{1/{2}, -k-{1}/{2}, -k-1}{1+\nu,1-\nu}{\dfrac{4b^2}{a^2}}
\;\;\;\;\text{for}\;\; a>2b,\; \nu \neq 0, 1, 2, \cdots,
\; k=0, 1, 2, \cdots
\end{split}
\]

\item Used in item (\ref{itemBessel2x})
\[
\begin{split}
&\bbint{0}{\infty} \dfrac{\sin(ax)J_\nu(bx)J_{-\nu}(bx)}{x^{2k+2}} {\rm d}x = -\operatorname{sinc}(\pi\nu) \dfrac{(-1)^k a^{2k+1}}{(2k+1)!} \\
&\hspace{0.4cm} 
\times \left[ \ln(a) \times 
\pFq{3}{2}{1/{2},-k,-k-1/{2}}{1+\nu,1-\nu}{\dfrac{4b^2}{a^2}}
\right. \\ &\hspace{1.2cm}
\left. -\sum_{l=0}^k \dfrac{(1/2)_l(-k)_l(-k-1/2)_l}{(1+\nu)_l(1-\nu)_l} \psi(2k-2l+2) \dfrac{(4b^2/a^2)^l}{l!} \right] \\
&\hspace{0.4cm} 
- (-1)^k \dfrac{2^{2k+1}}{\pi} \dfrac{b^{2k+2}}{k+1} \dfrac{a^{-1}}{(2k+1)!} \dfrac{(\Gamma(k+3/2))^2}{\Gamma(2+k+\nu)\Gamma(2+k-\nu)} \pFq{4}{3}{1/{2}, 1, 1, k+3/{2}}{2+k, 2+k+\nu, 2+k-\nu}{\dfrac{4b^2}{a^2}} \\
&\hspace{5.3cm} 
\;\;\;\;\text{for}\;\; a>2b,\; \nu \neq 0, 1, 2, \cdots,
\; k=0, 1, 2, \cdots
\end{split}
\]

\item Used in item (\ref{itemBessel3})
\[
\begin{split}
\setlength\arraycolsep{1pt}
\bbint{0}{\infty} &
\dfrac{\cos(ax)(J_n(bx))^2}{x^{2k+2n+2}} {\rm d}x
= -\dfrac{\pi b^{2n}}{2^{2n+1} (n!)^2} \dfrac{(-1)^{k} a^{2k+1}}{(2k+1)!} \\
&\hspace{-0.4cm}
\times \pFq{3}{2}{n + {1}/{2}, -k-{1}/{2}, -k}{1+n, 1+2n}{\dfrac{4b^2}{a^2}}
\;\;\;\;\text{for}\;\; a>2b,\; n = 0, 1, 2, \cdots,
\; k=0, 1, 2, \cdots
\end{split}
\]

\item Used in item (\ref{itemBessel3x})
\[
\begin{split}
&\bbint{0}{\infty} \dfrac{\cos(ax)(J_n(bx))^2}{x^{2k+2n+1}} {\rm d}x
= - (-1)^{k}\dfrac{(b/2)^{2n}}{(n!)^2} \dfrac{a^{2k}}{(2k)!} \\ 
&\hspace{0.2cm}\times \left[ \ln(a) \times 
\pFq{3}{2}{n+1/{2}, -k, -k+1/{2}}{n+1,2n+1}{\dfrac{4b^2}{a^2}}
\right. \\
&\hspace{0.8cm} \left. - \sum_{l=0}^{k} \dfrac{(n+1/2)_l (-k)_l (-k+1/2)_l }{(n+1)_l (2n+1)_l} \psi(2k-2l+1) \dfrac{(4b^2/a^2)^l}{l!} \right] \\
&\hspace{0.6cm}+ \dfrac{b^{2k+2n+2}}{a^2} \dfrac{(-1)^{k} }{\sqrt{\pi}} 
\dfrac{\Gamma(k+n+3/2)}{(k+1)!(k+n+1)!(k+2n+1)!} 
\pFq{4}{3}{1,1,3/{2},k+n+3/{2}}{2+k, 2+k+n,2+k+2n}{\dfrac{4b^2}{a^2}} \\
&\hspace{5.3cm} 
\;\;\;\;\text{for}\;\; a>2b,\; n = 0, 1, 2, \cdots,
\; k=0, 1, 2, \cdots
\end{split}
\]

\item Used in item (\ref{itemBessel4})
\[
\begin{split}
\setlength\arraycolsep{1pt}
\bbint{0}{\infty} &
\dfrac{\sin(ax)(J_n(bx))^2}{x^{2k+2n+3}} {\rm d}x
= -\dfrac{\pi b^{2n}}{2^{2n+1}(n!)^2} \dfrac{(-1)^{k} a^{2k+2}}{(2k+2)!} \\
&\hspace{-0.4cm} \times 
\pFq{3}{2}{n + {1}/{2}, -k-{1}/{2}, -k-1}{1+n, 1+2n}{\dfrac{4b^2}{a^2}}
\;\;\;\;\text{for}\;\; a>2b,\; n = 0, 1, 2, \cdots,
\; k=0, 1, 2, \cdots
\end{split}
\]

\item Used in item (\ref{itemBessel4x})
\[
\begin{split}
&\bbint{0}{\infty} \dfrac{\sin(ax)(J_n(bx))^2}{x^{2k+2n+2}} {\rm d}x
= - (-1)^{k}\dfrac{(b/2)^{2n}}{(n!)^2} \dfrac{a^{2k+1}}{(2k+1)!} \\ 
&\hspace{0.2cm}\times \left[ \ln(a) \times 
\pFq{3}{2}{n+1/{2}, -k, -k-1/{2}}{n+1,2n+1}{\dfrac{4b^2}{a^2}}
\right. \\
&\hspace{0.8cm} \left. - \sum_{l=0}^{k} \dfrac{(n+1/2)_l (-k)_l (-k-1/2)_l }{(n+1)_l (2n+1)_l} \psi(2k-2l+2) \dfrac{(4b^2/a^2)^l}{l!} \right] \\
&\hspace{0.6cm} - \dfrac{b^{2k+2n+2}}{a} \dfrac{(-1)^{k} }{\sqrt{\pi}} 
\dfrac{\Gamma(k+n+3/2)}{(k+1)!(k+n+1)!(k+2n+1)!} 
\pFq{4}{3}{1/{2},1,1,k+n+3/{2}}{2+k, 2+k+n,2+k+2n}{\dfrac{4b^2}{a^2}} \\
&\hspace{5.3cm} 
\;\;\;\;\text{for}\;\; a>2b,\; n = 0, 1, 2, \cdots,
\; k=0, 1, 2, \cdots
\end{split}
\]

\item Used in item (\ref{itemc25})
\[
\begin{split}
\bbint{0}{\infty} \dfrac{\text{Ai}(-ax)}{x^{3n}} {\rm d}x
=&\; 2\times (-1)^n \times 3^{n-1} a^{3n-1} \dfrac{n!}{(3n)!} 
\;\;\;\;\;\;\;\text{for}\;\; a > 0,
\; n=1, 2, 3, \cdots
\end{split}
\]

\item Used in item (\ref{itemc25})
\[
\begin{split}
\bbint{0}{\infty} & \dfrac{\text{Ai}(-ax)}{x^{3n-1}} {\rm d}x
= (-1)^n \dfrac{3^{-2n-1/3} a^{3n-2}}{\Gamma(n)\Gamma(n+1/3)}\\
&\times \left[
\ln\left(\dfrac{a^3}{9}\right)
+ \dfrac{2\sqrt{3}}{3}\pi
- \psi(n) - \psi(n+1/3)
\right] 
\;\;\;\text{for}\;\; a > 0,
\; n=1, 2, 3, \cdots
\end{split}
\]

\item Used in item (\ref{itemc25})
\[
\begin{split}
\bbint{0}{\infty} & \dfrac{\text{Ai}(-ax)}{x^{3n-2}} {\rm d}x
= (-1)^n \dfrac{3^{-2n+1/3} a^{3n-3}}{\Gamma(n)\Gamma(n-1/3)}\\
&\times \left[
\ln\left(\dfrac{a^3}{9}\right)
- \dfrac{2\sqrt{3}}{3}\pi
- \psi(n) - \psi(n-1/3)
\right] 
\;\;\;\text{for}\;\; a > 0,
\; n=1, 2, 3, \cdots
\end{split}
\]

\item Used in items (\ref{itemc25}) and (\ref{itemc26})
\[
\begin{split}
\bbint{0}{\infty} \dfrac{\text{Ai}(ax)}{x^{3n}} {\rm d}x
=&\; 3^{n-1} a^{3n-1} \dfrac{n!}{(3n)!}
\;\;\;\text{for}\;\; a > 0,
\; n=1, 2, 3, \cdots
\end{split}
\]

\item Used in items (\ref{itemc25}) and (\ref{itemc26})
\[
\begin{split}
\bbint{0}{\infty} \dfrac{\text{Ai}(ax)}{x^{3n-1}} {\rm d}x
=&\; \dfrac{3^{-2n-1/3} a^{3n-2}}{\Gamma(n)\Gamma(n+1/3)}
\left[
\ln\left(\dfrac{a^3}{9}\right)
- \dfrac{\sqrt{3}}{3}\pi
- \psi(n) - \psi(n+1/3)
\right] \\
&\;\;\;\text{for}\;\; a > 0,
\; n=1, 2, 3, \cdots
\end{split}
\]

\item Used in items (\ref{itemc25}) and (\ref{itemc26})
\[
\begin{split}
\bbint{0}{\infty} \dfrac{\text{Ai}(ax)}{x^{3n-2}} {\rm d}x
=&\; - \dfrac{3^{-2n+1/3} a^{3n-3}}{\Gamma(n)\Gamma(n-1/3)}
\left[
\ln\left(\dfrac{a^3}{9}\right)
+ \dfrac{\sqrt{3}}{3}\pi
- \psi(n) - \psi(n-1/3)
\right] \\
&\;\;\;\text{for}\;\; a > 0,
\; n=1, 2, 3, \cdots
\end{split}
\]

\item Used in items (\ref{itemc27}) and (\ref{itemc28})
\[
\begin{split}
\bbint{0}{\infty} \dfrac{\text{Ai}(ax)}{x^{m+\nu}} {\rm d}x
=&\; \pi \dfrac{3^{-(3+4m+4v)/6}}{2 a^{1-m-\nu}}
\dfrac{\csc\left(\dfrac{\pi}{3}(1-m-\nu)\right)
\csc\left(\dfrac{\pi}{3}(1+m+\nu)\right)}{
\Gamma\left(\dfrac{1+m+\nu}{3} \right)
\Gamma\left(\dfrac{2+m+\nu}{3} \right)}\\
& \;\text{for}\;\; a > 0,
\; 0 < \nu < 1,\;
m=1, 2, 3, \cdots
\end{split}
\]

\item Used in item (\ref{itemc28})
\[
\begin{split}
\bbint{0}{\infty} 
& \dfrac{\text{Ai}(-ax)}{x^{m+\nu}} {\rm d}x
= \pi
\dfrac{3^{-(3+4m+4v)/6}}{ a^{1-m-\nu}} \cos\left(\dfrac{\pi}{3}(m+\nu)\right)\\
&\times \dfrac{\csc\left(\dfrac{\pi}{3}(1-m-\nu)\right)
\csc\left(\dfrac{\pi}{3}(1+m+\nu)\right)}{
\Gamma\left(\dfrac{1+m+\nu}{3} \right)
\Gamma\left(\dfrac{2+m+\nu}{3} \right)} 
\;\;\;\;\;\; \text{for}\;\; a > 0,
\; 0 < \nu < 1,\;
m=1, 2, 3, \cdots 
\end{split}
\]

\item Used in item (\ref{itemc29})
\[
\begin{split}
\bbint{0}{\infty} 
& \dfrac{\text{Ai}(ax)\text{Ai}'(ax)}{x^{n+1}} {\rm d}x
= (-1)^n  \dfrac{2^{-(3-2n)/3} 3^{-(9-2n)/6} a^{n}\Gamma\left((3+2n)/6\right)}
{\pi^{3/2}\Gamma(n+1)} \\
&\times
\left( 
\left(\ln(12a^3) 
+ \psi\left( \dfrac{3+2n}{6}\right)
-3\psi(n+1)\right)
\cos\left( \dfrac{n\pi}{3}\right)
-\pi\sin\left(\dfrac{n\pi}{3}\right) \right) \\
&\;\;\;\;\;\;\;\;\;\;\;\;\;\;\;\;\;\;\;\;\;\;\;\;\;\;\;\;\;\;\;\;\;\;\;
\text{for}\;\; a > 0,
\; n = 0, 1, 2, \cdots 
\end{split}
\]

\item Used in item (\ref{itemc30})
\[
\begin{split}
\bbint{0}{\infty} \dfrac{\text{Ai}(ax)\text{Ai}'(ax)}{x^{m+\nu}} {\rm d}x
=&\; -12^{-(5-2m-2\nu))/6}
\pi^{-1/2} a^{m+\nu-1} 
\dfrac{\csc(\pi(m+\nu))}{\Gamma(m+\nu)}\\
&\times 
\Gamma\left(\dfrac{1+2m+2\nu}{6}\right)
\sin\left( \dfrac{\pi}{6}(1+2m+2\nu)\right)\\
& \;\text{for}\;\; a > 0,
\; 0 < \nu < 1,\;
m=1, 2, 3, \cdots
\end{split}
\]

\end{enumerate}

\end{appendix}


\section*{Data Availability}
The data that support the findings of this study are available from the corresponding author upon reasonable request.

\section*{Acknowledgement} This work was funded by the UP-System Enhanced Creative Work and Research Grant (ECWRG 2019-05-R).

\end{document}